\documentclass[12pt]{amsart}
\usepackage{amssymb,amsthm,amsmath,amscd,rsfs}
\usepackage{latexsym}

\newcommand{\NN}{{\mathbb{N}}}
\newcommand{\ZZ}{{\mathbb{Z}}}
\newcommand{\QQ}{{\mathbb{Q}}}
\newcommand{\RR}{{\mathbb{R}}}

\newcommand{\KK}{{\mathbb{K}}}
\newcommand{\GG}{{\mathbb{G}}}
\newcommand{\TT}{{\mathbb{T}}}

\newcommand{\XXX}{{\mathscr{X}}}
\newcommand{\YYY}{{\mathscr{Y}}}
\newcommand{\ZZZ}{{\mathscr{Z}}}

\newcommand{\UUU}{{\mathscr{U}}}

\newcommand{\BBB}{{\mathscr{B}}}
\newcommand{\EEE}{{\mathscr{E}}}
\newcommand{\AAA}{{\mathscr{A}}}
\newcommand{\TTT}{{\mathscr{T}}}
\newcommand{\CCC}{{\mathscr{C}}}
\newcommand{\DDD}{{\mathscr{D}}}

\newcommand{\SSS}{{\mathscr{S}}}
\newcommand{\PPP}{{\mathscr{P}}}

\newcommand{\OO}{{\mathcal{O}}}
\newcommand{\rank}{\operatorname{rk}}

\newcommand{\ch}{\operatorname{char}}

\newcommand{\bfK}{{\mathbf{K}}}

\newcommand{\Ker}{\operatorname{Ker}}

\newcommand{\Image}{\operatorname{Image}}

\newcommand{\Spec}{\operatorname{Spec}}
\newcommand{\Spf}{\operatorname{Spf}}
\newcommand{\Gal}{\operatorname{Gal}}

\newcommand{\cherncl}{{c}}

\newcommand{\ord}{\operatorname{ord}}

\newcommand{\Tr}{\operatorname{Tr}}

\newcommand{\aff}{\mathrm{aff}}

\newcommand{\an}{\mathrm{an}}

\newcommand{\val}{\operatorname{val}}
\newcommand{\Val}{\operatorname{Val}}
\newcommand{\str}{\operatorname{str}}
\newcommand{\initial}{\operatorname{in}}
\newcommand{\relin}{\operatorname{relin}}
\newcommand{\red}{\operatorname{red}}
\newcommand{\ndr}{\operatorname{nd-rk}}

\newcommand{\XX}{\mathrm{\mathcal{X}}}

\title
[Strict supports of canonical measures and applications to the GBC.]
{Strict supports of canonical measures and applications to the
geometric Bogomolov conjecture}
\author
{Kazuhiko Yamaki}
\date{November 19, 2015 (Version 6.2).}
\address
{Institute for Liberal Arts and Sciences,
Kyoto University, Kyoto, 606-8501, Japan}
\email{yamaki.kazuhiko.6r@kyoto-u.ac.jp}


\begin{document}

\theoremstyle{plain}
\newtheorem{Theorem}{Theorem}[section]
\newtheorem{Lemma}[Theorem]{Lemma}
\newtheorem{Proposition}[Theorem]{Proposition}
\newtheorem{Corollary}[Theorem]{Corollary}
\newtheorem{Main-Theorem}[Theorem]{Main Theorem}
\newtheorem{Theorem-Definition}[Theorem]{Theorem-Definition}
\theoremstyle{definition}
\newtheorem{Definition}[Theorem]{Definition}
\newtheorem{Remark}[Theorem]{Remark}
\newtheorem{Conjecture}[Theorem]{Conjecture}
\newtheorem{Claim}{Claim}
\newtheorem{Example}[Theorem]{Example}
\newtheorem{Key Fact}[Theorem]{Key Fact}
\newtheorem{ack}{Acknowledgments}       \renewcommand{\theack}{}
\newtheorem*{n-c}{Notation and convention}      
\newtheorem{citeTheorem}[Theorem]{Theorem}
\newtheorem{citeProposition}[Theorem]{Proposition}

\newtheorem{Step}{Step}

\renewcommand{\theTheorem}{\arabic{section}.\arabic{Theorem}}
\renewcommand{\theClaim}{\arabic{section}.\arabic{Theorem}.\arabic{Claim}}
\renewcommand{\theequation}{\arabic{section}.\arabic{Theorem}.\arabic{Claim}}

\def\Pf{\trivlist\item[\hskip\labelsep\textit{Proof.}]}
\def\endPf{\strut\hfill\framebox(6,6){}\endtrivlist}

\def\Pfo{\trivlist\item[\hskip\labelsep\textit{Proof of Proposition~\ref{ch-of-hyp}.}]}
\def\endPfo{\strut\hfill\framebox(6,6){}\endtrivlist}

\maketitle


\begin{abstract}
The Bogomolov conjecture claims that a closed subvariety
containing a dense subset of small points is a special kind of
subvarieties.
In the arithmetic setting
over number fields,
the Bogomolov conjecture for abelian varieties
has already been established as a theorem of Ullmo and Zhang,
but in the geometric setting over function field,
it has not
yet been solved completely.
There are only some partial results known
such as the
totally degenerate case
due to Gubler
and our recent work generalizing Gubler's result.

The key in establishing the previous results on
the Bogomolov conjecture
is the equidistribution method due to Szpiro--Ullmo--Zhang
with respect to
the canonical measures.
In this paper, we exhibit the limit of this method,
making an important contribution to the geometric version of the conjecture.
%
In fact,
by the crucial investigation of
the support of the canonical measure on a subvariety,
we show that 
the conjecture in full generality holds
if the conjecture holds for abelian varieties
which 
have anywhere good reduction.
As a consequence,
we establish a partial answer that
generalizes our previous result.
\end{abstract}

\section*{Introduction}


\renewcommand{\theTheorem}{\Alph{Theorem}}

\subsection{Background and our main results} \label{introtheorems}

The target of this paper is
the geometric Bogomolov conjecture,
which is 
the geometric version
of the theorem of Ullmo and Zhang.
In this subsection, we review the conjecture
with its background,
and then
we state our main results of this paper.

Let $K$ be a number field, 
or the function field 
of a normal projective variety
over an algebraically closed base field $k$.
We fix an algebraic closure $\overline{K}$ of $K$.
Let $A$ be an  abelian variety over $\overline{K}$ and
let $L$ be an ample line bundle on $A$.
Assume that $L$ is even,
i.e., $[-1]^{\ast} L = L$
for the endomorphism
$[-1] : A \to A$ give by $[-1] (a) = -a$.
Then the 
canonical height
function $\hat{h}_{L}$,
also called the N\'eron--Tate height,
associated with $L$
is a 
semi-positive definite quadratic form on 
$A 
\left(
\overline{K}
\right)$.
It is  well known that
$\hat{h}_{L} (x) = 0$ if $x$ is a torsion point.
Let $X$ be a closed subvariety of 
$A$.
We
put 
\[
X (\epsilon ; L)
:=
\left\{
x \in X 
\left(
\overline{K}
\right)
\left|
\hat{h}_{L} (x) \leq \epsilon
\right.
\right\}
\]
for
a real number $\epsilon > 0$.
Then the Bogomolov conjecture for abelian varieties
claims 
that $X (\epsilon ; L)$ is not Zariski-dense in $X$
for sufficiently small $\epsilon > 0$
unless $X$ is ``special'', such as a torsion subvariety
for example.
Note that if $X$ is a torsion subvariety, then $X (\epsilon ; L)$
is dense in $X$.

In the arithmetic case, the Bogomolov conjecture
has been established
by Ullmo for curves inside their Jacobians and by Zhang in general:

\begin{citeTheorem}[\cite{zhang2} and \cite{ullmo},
arithmetic version of
the Bogomolov conjecture]
\label{ABC}
Let $K$ be a number field.
Then,
$X (\epsilon ; L)$
is Zariski dense in $X$ for any $\epsilon > 0$
if and only if $X$ is a torsion subvariety.\footnote{The essential part
is the ``only if'' part.}
\end{citeTheorem}

The
``Bogomolov conjecture''
over a function field was widely open,
while the Bogomolov conjecture over number fields
was established.
In fact,
Moriwaki established in \cite{moriwaki5}
an arithmetic version of the Bogomolov conjecture
over a field $K$ \emph{finitely generated over $\QQ$},
but
it is a different problem because the
height in his setting is different from
the classical height over a function field, as mentioned in the introduction
of \cite{yamaki5}.

Nine years after Theorem~\ref{ABC},
Gubler made a breakthrough.
In \cite{gubler2},
he established 
the following theorem,
the same statement
as Theorem~\ref{ABC} 
for abelian varieties
totally degenerate at some place over a function field:

\begin{citeTheorem} [Theorem~1.1 in \cite{gubler2}] \label{thm:gubler}
Let $K$ be a function field.
Let $A$ be an abelian variety over $\overline{K}$.
Assume that $A$ is totally degenerate at some place.
Let $X \subset A$ be a closed subvariety.
If $X ( \epsilon ; L)$ is dense in $X$
for any $\epsilon > 0$, then
$X$ is a torsion subvariety.
\end{citeTheorem}

Then the next important problem was to generalize Theorem~\ref{thm:gubler}
for any abelian variety, but it was not a trivial task.
In fact,
Theorem~\ref{thm:gubler} does not holds for any abelian variety,
because in general an abelian variety has
subvarieties which are ``constant over $k$'', 
and this kind of subvarieties have
a dense subset of points of
height $0$.

Therefore
we needed to define some notion of subvarieties
which should be a counterpart
to the torsion subvarieties.
In \cite{yamaki5},
we introduced 
the notion of 
special subvarieties:
a closed subvariety
$X$ of $A$ is said to be \emph{special} if
there exist
an abelian variety $B$ over $k$, a closed subvariety $X' \subset B$,
a homomorphism $\phi : B_{\overline{K}} \to A$,
an abelian subvariety $A' \subset A$,
and a torsion point $\tau \in A ( \overline{K} )$ such that
$X = A' + \phi ( X'_{\overline{K}}) + \tau$,
where 
$B_{\overline{K}} = B \times_{\Spec k} \Spec \overline{K}$,
$X'_{\overline{K}} = X' \times_{\Spec k} \Spec \overline{K}$
(cf. Lemma~\ref{traceintro}~(1)).
Note that if $X$ is a special subvariety, then $X (\epsilon ; L)$ is dense
in $X$ 
for any $\epsilon > 0$
(cf. \cite[Corollary~2.8]{yamaki5}).

Using the notion of special subvarieties as a counterpart,
we formulated in \cite{yamaki5}
the following conjecture,
called 
the \emph{geometric Bogomolov conjecture}.

\begin{Conjecture} [cf. Conjecture 2.9 in \cite{yamaki5} and
Conjecture~\ref{GBCforAV}] \label{intGBCforAV}
Let $K$ be a function field.
Let $X$ be an irreducible closed subvariety of $A$.
Then,
if
$X (\epsilon ; L)$
is Zariski dense in $X$ 
for any $\epsilon > 0$,
then
$X$ is a special subvariety.
\end{Conjecture}


Conjecture~\ref{intGBCforAV} 
is still an open problem.
The best partial solution known so far
is the following,
where $A_v$ is the Berkovich analytic space associated to $A$ 
over $v$, and $b (A_v)$ is the abelian rank of $A_v$
(cf. \S~\ref{RtM}):

\begin{citeTheorem} [Corollary 5.6 in \cite{yamaki5}]
\label{b=1}
Let 
$A$ be an abelian variety over $\overline{K}$.
Suppose that there exists a place 
$v$ such that 
$b (A_{v}) \leq 1$.
Then the geometric Bogomolov conjecture holds for $A$.
\end{citeTheorem}

Theorem~\ref{b=1}
generalizes Theorem~\ref{thm:gubler}.
Indeed,
if
$A_v$ is totally degenerate, then $b(A_v) = 0 \leq 1$
and any special subvariety of $A$
is a torsion subvariety.

In this paper, we make
a major contribution
to
the geometric Bogomolov conjecture,
exhibiting the limit
of the methods of Ullmo, Zhang,
Gubler, and us on the Bogomolov conjecture.
All the results  explained so far
have been proved by a common method, i.e., the equidistribution method
due to Szpiro--Ullmo--Zhang
\cite{SUZ}.
We will resume more detailed description 
of this method in the next subsection.

We are now ready to state our main results.
Let $K$ be a function field.
For an abelian variety $A$ over $\overline{K}$, 
it can be seen that there exists a unique
maximal nowhere degenerate\footnote{An
abelian variety $B$ over $\overline{K}$
is said to be \emph{nowhere degenerate}
if $B$ has good reduction at any place.
}
abelian subvariety $\mathfrak{m}$ of $A$ (cf. \S~\ref{ndrc}).

\begin{Theorem} [Corollary~\ref{maintheorem3}]
\label{maintheorem3int}
Let $A$ be an abelian variety 
over $\overline{K}$
and let $\mathfrak{m}$ be 
the
maximal nowhere degenerate abelian subvariety of $A$.
Then,
the geometric Bogomolov conjecture holds for $A$
if and only if it holds for  $\mathfrak{m}$.
\end{Theorem}

As a consequence
of Theorem~\ref{maintheorem3int},
the geometric Bogomolov conjecture 
for abelian varieties is reduced to the conjecture
for those
without places of degeneration (cf. Conjecture~\ref{GBCforAVwithtoutdeg}).
This theorem 
also shows that the conjecture holds for simple abelian varieties
which are degenerate at some place
(cf. Remark~\ref{rem:generalization}).

Theorem~\ref{maintheorem3int} generalizes
all the important previous results concerning the
geometric Bogomolov conjecture.
To see that,
we define the nowhere-degeneracy rank of $A$ to be
$\ndr ( A ) := \dim \mathfrak{m}$
(cf. Definition~\ref{ndr}).
Since
the geometric Bogomolov conjecture holds for all elliptic curves,
the following theorem follows 
from Theorem~\ref{maintheorem3int}.\footnote{In
this paper, we 
give a direct proof of
Theorem~\ref{maintheoremint} not via
Theorem~\ref{maintheorem3int}.}

\begin{Theorem} [cf. Corollary~\ref{GBCfornondeggeq1}]
\label{maintheoremint}
Let $A$ be an abelian variety over $\overline{K}$.
Assume that 
$\ndr ( A) \leq 1$.
Then the geometric Bogomolov conjecture holds for $A$.
\end{Theorem}


We remark
that
$b (A_{v}) \geq \ndr (A)$ for any $v \in 
M_{\overline{K}}$.
Thus
Theorem~\ref{maintheoremint} generalizes
Theorem~\ref{b=1}
(cf. Remark~\ref{rem:generalization}) and hence Theorem~\ref{thm:gubler}.

This paper also discusses 
another version of the Bogomolov conjecture,
called
the Bogomolov conjecture \emph{for curves} over function fields.
This conjecture
claims that
the embedded curve in its Jacobian should have
only a finite number of small points unless it is isotrivial
(cf.
Conjecture~\ref{GBCforcurves}
for the precise statement).
It is still an open problem
with partial solutions,
while
the arithmetic version
is established
by
Ullmo
in \cite{ullmo}.
Under the assumption that
$K$ is the function field of a curve
over a field of characteristic zero,
Cinkir established in \cite{cinkir} an affirmative answer
to Conjecture~\ref{EGBCforcurves},
an effective version of 
this conjecture.
There is not a satisfactory answer
to the conjecture
in positive characteristic
except for
some partial answers.
We refer to \S~\ref{GBCforC} for more details.
In this paper, 
applying our arguments
on Conjecture~\ref{intGBCforAV}
to this problem, 
we 
make nontrivial contribution
including the case of positive characteristic.

\subsection{Idea of the proof} \label{strategy}


Our basic strategy
for the geometric Bogomolov conjecture
is based on the 
non-archimedean analogue of the proof of
Theorem~\ref{ABC}, as is done in \cite{gubler2} and \cite{yamaki5}.
Zhang's proof 
of Theorem~\ref{ABC}
relies on the equidistribution theorem of Szpiro--Ullmo--Zhang
with respect to the canonical measure
over an archimedean place.
Although we do not have an archimedean place over function fields,
we can define
the canonical measure on
the analytic space over a non-archimedean
place,
and it
plays a crucial role as a counterpart to the canonical measure over an
archimedean place.
Let $K$ be a function field.
Let 
$\KK = 
\overline{K}_{v}$ be the completion of $\overline{K}$ 
with respect to a place $v$ of $\overline{K}$ (cf. \S~\ref{NT}).
Let
$A$ be an abelian variety over $\KK$, $X \subset A$ a closed subvariety
of dimension $d$
and let $L$ be an even ample line bundle on $A$.
By the theorem of square, we have $[n]^{\ast} L = L^{n^{2}}$ for any $n \in \ZZ$.
Suppose that all of them can be defined over $\overline{K}$.
Then it is known that there exists a canonical metric
on $L$ defined by 
the condition $[n]^{\ast} || \cdot || = || \cdot ||^{n^{2}}$ via the above
identification.
Let $\overline{L}$ be the line bundle $L$ with a canonical metric.
Then
we have a measure 
$\cherncl_{1} (\overline{L}|_{X})^{\wedge d}$ 
on the associated Berkovich analytic space $X^{\an}$,
which was originally
introduced by Chambert-Loir
in \cite{chambert-loir}.
It
is a semipositive Borel measure,
and thus
the
probability measure
\[
\mu_{X^{\an}, \overline{L}} :=
\frac{1}{\deg_{L} X} \cherncl_{1} (\overline{L}|_{X})^{\wedge d}
\]
on $X^{\an}$ is defined.
This is
also called a canonical measure.


%


With respect to the canonical measures
over non-archimedean analytic spaces,
an argument analogous to Zhang's proof
works in some extent.
We refer to
the introduction of 
\cite{yamaki5} for details of the argument.
However,
if we wish to
establish some results concerning the geometric Bogomolov conjecture
by the same way,
we should need
some particular information
on canonical measures.
In the setting of Theorem~\ref{b=1}, actually,
we focused
on the minimal dimension
of the components of the support of the tropicalized canonical measure
(see \cite{yamaki5} for details).
In the setting of 
Theorems~\ref{maintheoremint} 
and 
\ref{maintheorem3int},
we will introduce strict supports of a canonical measure. 
Their study will be crucial in this paper and leads to the required information.

Our main results Theorem~\ref{maintheoremint} and Theorem~\ref{maintheorem3int}
are consequences of
Theorem~\ref{maintheorem1}.
We say that a closed subvariety of an abelian variety is tropically trivial
if the canonical tropicalization map of the abelian variety
contracts the subvariety
to a point for every place of the function field.
Theorem~\ref{maintheorem1}
says that, if the closed subvariety $X$ of $A$ has dense small points\footnote{
We can show that the property
``$X (\epsilon ; L)$ is dense in $X$ for any $\epsilon > 0$''
does not depend on an even ample line bundle $L$.
We say $X$ has dense small points if, 
for some (and hence any) $L$, $X (\epsilon ; L)$
is dense for any $\epsilon > 0$
(cf. \cite[Definition 2.2]{yamaki5}).},
then $X/G_X$ is tropically trivial,
where
$G_{X}$ is the stabilizer of $X$,
i.e.,
$G_{X} = \{ a \in A \mid a + X \subset X \}$.

The proof of Theorem~\ref{maintheorem1}
follows
the same strategy
as Zhang's proof.
In fact, the strategy
works well
by virtue of
Proposition~\ref{nondeg-supp2},
which 
gives us the
crucial information about the canonical measure.
In the remaining of this subsection,
we
describe what Proposition~\ref{nondeg-supp2} says,
and 
how it is used to prove
Theorem~\ref{maintheorem1}.
Let $A$ be an abelian variety over $\KK$ and let
$X \subset A$ be a closed subvariety.
Let
$\XXX'$ be a strictly semistable proper formal scheme with 
Raynaud generic fiber $X'$
viewed as a Berkovich analytic space (cf. \S~\ref{preliminary}),
and let $f : X' \to A^{\an}$ be a 
generically finite morphism such that $f (X') = X^{\an}$
with some technical assumptions.
Let $\overline{L}$ be an even ample line bundle on $A$
with a canonical metric.
Then we have the probability measure
$\mu_{X', f^{\ast} \overline{L}}$ on $X'$,
which has
the property that
$f_{\ast} \mu_{X', f^{\ast} \overline{L}} = \mu_{X^{\an} , L}$.

Let $S (\XXX')$ be the skeleton of $\XXX'$
(cf. \S~\ref{ssfss}).
It is a simplicial set and a subspace of $X'$;
we have a canonical simplex $\Delta_{S}$ for each
stratum $S$ of the special fiber of $\XXX'$
and $S (\XXX') = \bigcup_{S} \Delta_{S}$.
In
\cite{gubler4}, Gubler defined the 
notion of {non-degenerate} canonical simplices with respect
to $f$,
and
showed that
$\mu_{X', f^{\ast} \overline{L}}$
is a finite sum of the Lebesgue measures on the
non-degenerate canonical simplices.
It follows that
the support
$S_{X^{\an}}$ 
of $\mu_{X^{\an} , L}$
coincides with the image of the union of the non-degenerate
canonical simplices by $f$.
Further, he showed
that
$S_{X^{\an}}$ 
has a unique piecewise linear structure
such that the restriction of $f$ to each non-degenerate
canonical simplex is a piecewise linear map.

We fix a sufficiently refined polytopal decomposition of 
$S_{X^{\an}}$.\footnote{It should be a subdivisional
one with the terminology in \S~\ref{ssca}.}
A polytope $\sigma$ 
in $S_{X^{\an}}$ is called a \emph{strict support} of $\mu_{X^{\an}}$
if $\mu_{X^{\an}} - \epsilon \delta_{\sigma}$ is semipositive
for small $\epsilon > 0$
(cf. Definition~\ref{defSS}),
where $\delta_{\sigma}$ 
is the push-out to $S_{X^{\an}}$
of the Lebesgue measure on $\sigma$.
Suppose that $\sigma$ is a 
strict support
of $\mu_{X^{\an},L}$.
Then we
see that there is
a non-degenerate canonical simplex $\Delta_{S}$ with
$\sigma \subset f ( \Delta_{S} )$ and $\dim \sigma = \dim f ( \Delta_{S} )$
(cf. Lemma~\ref{nondeg-supp}).
Furthermore,
Proposition~\ref{nondeg-supp2}
with Lemma~\ref{rem:assumptionofprop}
shows that,
for
\emph{any} stratum
$\Delta_{S}$  
of $S( \XXX' )$
with $\sigma \subset
f( \Delta_{S} ) $ and 
$\dim \sigma = \dim f( \Delta_{S} )$,
we have
$\dim f ( \Delta_{S} ) = \dim \Delta_{S}$.\footnote{Precisely,
we should write $\overline{f}_{\aff} ( \Delta_{S} )$ instead of 
$f ( \Delta_{S} )$, etc.,
but we use this notion in this subsection even though it is not precise;
we are just sketching the idea of the
proof.}

Let us give an outline of
the proof of Theorem~\ref{maintheorem1} now.
To argue by contradiction,
we suppose that
$X / G_{X}$ is not tropically trivial
but
$X \subset A$ has dense small points. 
Replacing $X$ with $X / G_{X}$,
we may further assume that $X$ has trivial stabilizer.
Then, it follows from the
tropical non-triviality of
$X$ that 
$S_{X_{v}}$ 
has positive dimension for some place $v$,
where $X_v$ is the associated Berkovich space at $v$.
Thus
$\mu_{X_{v}}$ has a strict support of positive dimension.

Consider the homomorphism 
$A^{N} \to A^{N-1}$ 
defined by 
$(x_{1} , \ldots , x_{N}) \mapsto (x_{2} - x_{1} , \ldots , x_{N} - x_{N-1})$.
We set $Z := X^{N}$ and let $Y$ be the image of $Z$ by this homomorphism.
Let
$\alpha : Z \to Y$ 
be the restriction of this homomorphism.
By the triviality of $G_X$,
there exists $N \in \NN$ such
that $\alpha$ is generically finite.
Note that 
$\alpha$ contracts the
diagonal of $Z = X^{N}$ to a point.
Since
$\mu_{X_{v}}$ has a strict support of positive dimension
and since
$\mu_{Z_{v}}$ is the product of $N$ copies of
$\mu_{X_{v}}$, 
it follows that
there is a
strict support $\sigma$ of $\mu_{Z_{v}}$ 
with $\dim \alpha (\sigma) < \dim \sigma$.

Since $Z$ has dense small points,
we obtain $\alpha_{\ast} ( \mu_{ Z_{v}})
= \mu_{Y_{v}}$ by the equidistribution theorem.
Since $\sigma$ is a strict support of $\mu_{Z_{v}}$,
it follows
that $\alpha ( \sigma )$ is a strict support of $\mu_{Y_{v}}$.

By de Jong's alteration theorem, we have
a proper strictly semistable formal scheme $\ZZZ'$
with a generically finite surjective morphism
$g : (\ZZZ')^{\an} \to Z_{v}$.
Since $\alpha$ is
a generically finite surjective morphism,
the morphism $h := \alpha \circ g$ 
is also generically finite surjective morphism.
Since $\sigma$ is a strict support, 
there exists a 
non-degenerate canonical simplex $\Delta_{S}$
of $S (\ZZZ')$
with respect to $g$ such that
$\sigma \subset g (\Delta_{S})$ and $\dim \sigma = \dim g ( \Delta_{S} )$.
Then we have
$\alpha ( \sigma ) \subset h (\Delta_{S})$ and
$\dim \alpha ( \sigma ) = \dim h (\Delta_{S})$.
Since $\alpha (\sigma)$ is a strict support of $\mu_{Y_{v}}$,
it follows form
Proposition~\ref{nondeg-supp2}
that
$\dim h ( \Delta_{S}) = \dim \Delta_{S}$.
On the other hand,
the inequality $\dim \alpha (\sigma) < \dim \sigma $
tells us
$\dim h ( \Delta_{S})
=
\dim \alpha ( g  ( \Delta_{S}) ) < \dim \Delta_{S}$.
This is a contradiction. Thus 
Theorem~\ref{maintheorem1} follows.

Finally in this subsection, we notice a limit of this strategy.
As mentioned in the previous subsection,
Conjecture~\ref{intGBCforAV} can be reduced to
Conjecture~\ref{GBCforAVwithtoutdeg}, namely,
the geometric Bogomolov conjecture 
\emph{for nowhere degenerate abelian varieties}.
However,
the strategy used here is useless 
for 
such abelian varieties
because the support of the canonical measure
of a closed subvariety is 
the Dirac measure 
of a finite set.
This fact suggests that
our theorems \ref{maintheoremint} and \ref{maintheorem3int}
are the ultimate which
can be reached with the strategy based on
equidistribution theorems. 

\subsection{Organization}
This article consists of eight sections.
In \S~\ref{preliminary},
we recall
some basic facts on non-archimedean geometry.
In 
\S~\ref{torus-equivalence},
we introduce an affinoid-torus-action on a formal fiber
of the Berkovich analytic space associated to a
strictly semistable formal scheme,
and we show that this action is compatible with the reduction
(cf. Lemma~\ref{torus-action-ff}).
In \S~\ref{RaynaudextMumford},
we recall the 
Raynaud extension of an abelian variety
and
its tropicalization.
We also recall Mumford models
of the Raynaud extension.
In \S~\ref{non-degstrata2},
we
establish
that a morphism from a strictly semistable formal scheme 
to a Mumford model
of an abelian variety induces a torus-equivariant morphism
between their strata (cf. Proposition~\ref{equivariantdiagram})
with the use of Lemma~\ref{torus-action-ff}.
Furthermore,
we show a key
Lemma~\ref{mainproposition},
which will be crucially used to show that a stratum
over a strict support does not collapse.
We introduce the notion of strict support in \S~\ref{strictsupport},
and prove Proposition~\ref{nondeg-supp2} which was explained in the previous
subsection.
In 
\S~\ref{application},
using Proposition~\ref{nondeg-supp2},
we show that
if the closed subvariety $X$ of $A$ has dense small points,
then $X/G_X$ is tropically trivial
(cf. Theorem~\ref{maintheorem1}).
In \S~\ref{mainresults},
using Theorem~\ref{maintheorem1},
we establish results concerning the geometric Bogomolov conjecture
including
the main theorems mentioned in \S~\ref{introtheorems}.
Some non-trivial remarks on the conjecture for curves are made
in the last section.

\subsection*{Acknowledgments}
I would like to express my gratitude to Professor Matthew Baker and
Professor Amaury Thuillier.
We had a conversation
on 
the relationship between the
tropical
canonical measures and the initial degenerations
in the Bellairs workshop 2011 at Barbados,
which motivated me to establish
Lemma~\ref{mainproposition}.
I also thank organizers of that workshop,
especially Professor Xander Faber, for inviting me.
I would like to thank Professor Shu Kawaguchi
for his constructive comments
on a preliminary version of this paper.
My deepest appreciation goes to Professor Walter Gubler,
who read the first version carefully
and gave me a lot of valuable comments.
Finally, I would thank the referees,
who gave me a lot of helpful comments and suggestions.
This work was partly supported by
KAKENHI(21740012).

\renewcommand{\theTheorem}{\arabic{section}.\arabic{Theorem}}

\section{Preliminary} \label{preliminary}

\subsection{Convention and terminology} \label{conv-term}
When we write $\KK$,
it is an algebraically closed field 
which is complete with respect to a
non-trivial
non-archimedean absolute value $| \cdot | : \KK \to \RR_{\geq 0}$.
We put 
$
\KK^{\circ} := 
\{
a \in \KK \mid | a | \leq 1
\}
$,
the ring of integers of $\KK$,
$
\KK^{\circ \circ} := 
\{
a \in \KK \mid | a | < 1
\}
$,
the maximal ideal of the valuation ring $\KK^{\circ}$,
and
$\tilde{\KK} := \KK^{\circ} / \KK^{\circ \circ}$,
the residue field.
We put $\Gamma := \{ -\log |a| \mid a \in \KK^{\times} \}$,
the value group of $\KK$.

We also fix the notation used in convex geometry.
See \cite[\S~6.1 and Appendix~A]{gubler1}
for more details.
A polytope $\Delta$ of $\RR^{n}$
is said to be \emph{$\Gamma$-rational}
if it can be given as an intersection of
subsets of form $\{ \mathbf{u} \mid \mathbf{m} \cdot \mathbf{u} \geq c \}$ 
(closed half-space)
for some $\mathbf{m} \in \ZZ^{n}$ and $c \in \Gamma$.
When $\Gamma = \QQ$, a $\Gamma$-rational polytope is called
a \emph{rational} polytope.
A \emph{closed face} of $\Delta$ is a polytope
which is $\Delta$
itself or is of form $H \cap \Delta$ 
where $H$ is the boundary of a closed half-space
of $\RR^{n}$ containing $\Delta$.
An \emph{open face} of $\Delta$ means
the relative interior of any closed face of $\Delta$.

Let $\Omega$ be a subset of $\RR^{n}$.
A 
\emph{$\Gamma$-rational
polytopal decomposition $\CCC$ of $\Omega$} is a locally finite family of 
$\Gamma$-rational polytopes of $\RR^{n}$
such that
\begin{itemize}
\item
$\bigcup_{\Delta \in \CCC} \Delta = \Omega$, 
\item
for any $\Delta , \Delta' \in \CCC$, $\Delta \cap \Delta'$
is an empty set or
a face of $\Delta$ and $\Delta'$,
and
\item
for any $\Delta \in \CCC$, any closed face of $\Delta$ is in $\CCC$.
\end{itemize}

Let $\CCC$ be a $\Gamma$-rational polytopal decomposition of $\RR^{n}$.
Let $\Lambda$ be a lattice of $\RR^{n}$.
The polytopal decomposition $\CCC$ is said to be 
\emph{$\Lambda$-periodic}
if, for any $\Delta \in \CCC$,
we have $\lambda + \Delta \in \CCC$ for any $\lambda \in \Lambda$,
and the restriction of the quotient map $\RR^{n} \to \RR^{n} / \Lambda$ to
$\Delta$ is a homeomorphism to its image.
For a $\Lambda$-periodic 
$\Gamma$-rational polytopal decomposition $\CCC$ of $\RR^{n}$,
we set
$\overline{\CCC} := \{ \overline{\Delta} \mid \Delta \in \CCC \}$,
where $\overline{\Delta}$ is the image of $\Delta$ by
the quotient map $\RR^{n} \to \RR^{n} / \Lambda$.
Such a $\overline{\CCC}$ is called a 
\emph{$\Gamma$-rational polytopal decomposition of $\RR^{n} / \Lambda$}.

We mean by an algebraic variety over a field an
irreducible, reduced, and
separated scheme of finite type over the field unless
otherwise specified.

\subsection{Berkovich analytic spaces} \label{bersp}

We recall some notions and properties on
analytic spaces associated to
admissible formal schemes and those associated to
algebraic varieties,
as far as we use later.
For details, we refer to
Berkovich's original papers
\cite{berkovich1,berkovich2,berkovich3,berkovich4}
or Gubler's
expositions in \cite{gubler1,gubler4}.

Let $\KK \langle x_{1} , \ldots , x_n \rangle$
be the Tate algebra over $\KK$, i.e., 
the completion of the polynomial ring 
$\KK [ x_{1} , \ldots , x_n ]$
with respect to the Gauss norm.
A $\KK$-algebra $\mathfrak{A}$ isomorphic to 
$\KK \langle x_{1} , \ldots , x_n \rangle / I$
for some ideal $I$ of $\KK \langle x_{1} , \ldots , x_n \rangle$
is called a \emph{$\KK$-affinoid algebra}.
Let $\mathrm{Max} (\mathfrak{A})$ be the maximal spectrum
of $\mathfrak{A}$.
Let $|\cdot|_{\sup} : \mathfrak{A} \to \RR$
denote the supremum semi-norm over $\mathrm{Max} (\mathfrak{A})$.
The Berkovich spectrum of $\mathfrak{A}$
is
the set of multiplicative seminorms 
$\mathfrak{A}$
bounded with $|\cdot|_{\sup}$.
A (Berkovich) analytic space over $\KK$
is given by an atlas of Berkovich spectrums of $\KK^{\circ}$-affinoid
algebras.

A $\KK^{\circ}$-algebra 
is called an
 \emph{admissible $\KK^{\circ}$-algebra} if it does not have any
$\KK^{\circ}$-torsions and it is isomorphic to
$\KK^{\circ} \langle x_{1} , \ldots , x_{n} \rangle / I$
for some $n \in \NN$ and for some ideal
$I$ of $\KK^{\circ} \langle x_{1} , \ldots , x_{n} \rangle$.
Note that an admissible $\KK^{\circ}$-algebra is flat over $\KK^{\circ}$.
The formal spectrum 
of an admissible $\KK^{\circ}$-algebra
is called an \emph{affine admissible formal scheme}.
A formal scheme over $\KK^{\circ}$ is called an 
\emph{admissible formal scheme}\footnote{When we say
a formal scheme,
it always means
an admissible formal scheme in this article.}
if it has a locally finite open atlas of affine admissible
formal schemes.
Note that an admissible formal scheme
is flat over $\KK^{\circ}$.
For an admissible formal scheme $\XXX$,
we write $\tilde{\XXX} := \XXX \times_{\Spf \KK^{\circ}} 
\Spec \tilde{\KK}$ for the special fiber.
For a morphism $\varphi : \XXX \to \YYY$
of admissible formal schemes,
we write $\tilde{\varphi} : \tilde{\XXX}
\to \tilde{\YYY}$ for the induced morphism
between their special fibers.

Let $\XXX$ be an
admissible formal scheme over $\KK^{\circ}$.
Then we can associate a Berkovich analytic space $\XXX^{\an}$,
called the \emph{(Raynaud) generic fiber} of $\XXX$.
Further, we have a map
$
\mathrm{red}_{\XXX} :
\XXX^{\an} \to \tilde{\XXX}
$,
called the \emph{reduction map}.
It is known that
the reduction map is surjective.
Let $Z$ be a dense open subset of an irreducible component of $\tilde{\XXX}$
with the generic point $\eta_{Z} \in Z$.
Then there is a point $\xi_{Z} \in \XXX^{\an}$
with $\mathrm{red}_{\XXX} (\xi_{Z}) = \eta_{Z}$.
If the special fiber $\tilde{\mathscr{X}}$ is reduced, then
such a point $\xi_{Z}$ is unique, and we refer to it as 
\emph{the point corresponding
to (the generic point of) $Z$}.


We can associate
an analytic space to an
algebraic variety $X$ over $\KK$
as well,
and
we write
$X^{\an}$
for
the analytic space associated to 
$X$.
There is a 
natural inclusion
$X (\KK) \subset X^{\an}$.
We should remind the relationship between
the analytic space associated to
an algebraic variety and that done to an admissible formal scheme.
Let $\mathbf{X}$ 
be a scheme
flat and of finite type
over $\KK^{\circ}$ with the generic fiber $X$.
Let $\widehat{\mathbf{X}}$ be the formal completion with respect to
a principal open 
ideal of $\KK^{\circ}$ contained in $\KK^{\circ \circ}$.
Then 
$\widehat{\mathbf{X}}$ is an admissible formal scheme
and
$\widehat{\mathbf{X}}^{\an}$ is an analytic subdomain of
$X^{\an}$.
If
$\mathbf{X}$ is proper over $\KK^{\circ}$,
we have $\widehat{\mathbf{X}}^{\an} = X^{\an}$.

For a given analytic space $X$, an admissible formal scheme
having $X$ as
the generic fiber
is called a
\emph{formal model} of $X$.
Note that a formal model is flat over $\KK^{\circ}$ by definition.
Let $Y$ be a closed 
analytic subvariety of
$X$
and let
$\XXX$ 
be a formal model of $X$.
Then there
exists a unique admissible formal subscheme
$\YYY \subset \XXX$ with $\YYY^{\an} = Y$.
We call $\YYY$ the \emph{closure} of $Y$ in $\XXX$.

\subsection{Tori} \label{tori}

Let $\GG_{m}^{n}$ denote the split torus of rank $n$
over $\KK$.
Let
$x_{1} , \ldots , x_{n}$ denote the standard coordinates of $\GG_{m}^{n}$
unless otherwise noted, so that
$
\GG_{m}^{n} = \Spec \KK [x_{1}^{\pm} , \ldots , x_{n}^{\pm}]
$.
Let $(\GG_{m}^{n})^{\an}$ be the analytic space associated to $\GG_{m}^{n}$.
We
set
\[
(\GG_{m}^{n})_{1}^{\mathrm{f-sch}} :=
\Spf \KK^{\circ} [x_{1}^{\pm} , \ldots , x_{n}^{\pm}]
,
\]
the formal torus over $\KK^{\circ}$,
writing $(\GG_{m}^{n})^{\an}_{1}$ for the generic fiber of 
$(\GG_{m}^{n})_{1}^{\mathrm{f-sch}}$.
We call $(\GG_{m}^{n})_{1}^{\mathrm{f-sch}}$ the \emph{canonical model} of
$(\GG_{m}^{n})^{\an}_{1}$.
The analytic space
$(\GG_{m}^{n})^{\an}_{1}$ is 
an analytic subgroup of $(\GG_{m}^{n})^{\an}$
as well as
an affinoid subdomain of $(\GG_{m}^{n})^{\an}$.
The special fiber of $(\GG_{m}^{n})_{1}^{\mathrm{f-sch}}$
is an algebraic torus 
$(
\GG_{m}^{n})_{\tilde{\KK}}
=
\Spec \tilde{\KK} [ x_{1}^{\pm} , \ldots , x_{n}^{\pm} ]
$,
which
we call
the \emph{canonical reduction}
of $(\GG_{m}^{n})^{\an}_{1} $.

Each element $p \in (\GG_{m}^{n})^{\an}$ is regarded as a seminorm
on $\KK [ x_{1}^{\pm} , \ldots , x_{n}^{\pm} ]$.
We define a map $\val : (\GG_{m}^{n})^{\an} \to \RR^{n}$,
called the \emph{valuation map},
by
\[
p \mapsto ( - \log  p (x_{1}) , \ldots , - \log p (x_{n}))
.
\]

Let $\Delta$ be a $\Gamma$-rational polytope of $\RR^{n}$.
Then $U_{\Delta} := \val^{-1} ( \Delta)$ is an analytic subdomain of $(\GG_{m}^{n})^{\an}$.
In fact, it is the Berkovich spectrum of the affinoid algebra
\[
\KK \langle U_{\Delta} \rangle
:=
\left\{
\left.
\sum_{\mathbf{m} \in \ZZ^{n}}
a_{\mathbf{m}} x_{1}^{m_{1}} \cdots x_{n}^{m_{n}} 
\right|
\lim_{|\mathbf{m}| \to \infty}
v (a_{\mathbf{m}}) + \mathbf{m} \cdot \mathbf{u} = \infty
\text{ for any $\mathbf{u} \in \Delta$}
\right\}
.
\]
Note that $\val^{-1} ( \mathbf{0}) = (\GG_{m}^{n})^{\an}_{1}$.
We refer to  \cite[4.3]{gubler2} for more details.


\section{Torus-action on the formal fiber} \label{torus-equivalence}

An admissible formal scheme $\XXX'$ is called
a \emph{strictly semistable} formal scheme if any point of $\XXX'$
has an open neighborhood $\UUU'$
and an \'etale morphism
\addtocounter{Claim}{1}
\begin{align} \label{ssetale}
\psi :
\UUU' \to
\SSS :=
\Spf \KK^{\circ} 
\langle
x'_{0} , \ldots , x'_{d}
\rangle
/
(x'_{0} \cdots x'_{r} - \pi)
,
\end{align}
where $\pi \in \KK^{\circ \circ} \setminus \{ 0 \}$.

Throughout this section,
let $\XXX'$ be a strictly semistable
formal scheme over $\KK^{\circ}$ with generic fiber $X' = (\XXX')^{\an}$.
Let $\red_{\XXX'} : X' \to \tilde{\XXX'}$ be the
reduction map.
For a closed point $\tilde{p} \in \tilde{\XXX'}$,
we put
$X'_{+} ( \tilde{p} ) := (\red_{\XXX'})^{- 1} ( \tilde{p} )$
and call it the \emph{formal fiber} over $\tilde{p}$.
It
is an open subdomain of $X'$.
See \cite[2.8]{gubler1} for more details.



The purpose of this section is to
define a torus-action on $X'_{+} ( \tilde{p} )$
and show
that it is compatible with the reduction
as Lemma~\ref{torus-action-ff} claims.
This lemma
will be used in
the proof of Proposition~\ref{equivariantdiagram}.

\subsection{Strictly semistable formal schemes, their skeletons,
and subdivision} \label{ssfss}

We begin by recalling the notion of stratification of a 
reduced separated scheme $Y$ of finite type
over a field.
We start with $Y^{(0)} := Y$.
For each $r \in \NN$,
let $Y^{(r)} \subset Y^{(r - 1)}$ be the complement of the set of 
normal points in $Y^{(r - 1)}$.
Since the set of normal points is open and dense,
we obtain a chain of closed subsets:
\[
Y = Y^{(0)} \supsetneq Y^{(1)} \supsetneq \cdots \supsetneq 
Y^{(s)} \supsetneq Y^{(s+1)} = \emptyset
.
\]
The irreducible components of $Y^{(r)} \setminus Y^{(r + 1)}$
($0 \leq r \leq s$)
are called the \emph{strata of $Y$},
and the set of strata is denoted by $\str (Y)$.

We take an \'etale morphism as in (\ref{ssetale}).
Putting
$\SSS_{1} := \Spf \KK^{\circ} 
\langle
x'_{0} , \ldots , x'_{r}
\rangle
/
(x'_{0} \cdots x'_{r} - \pi)$
and
$\SSS_{2} := \Spf \KK^{\circ} 
\langle
x'_{r+1} , \ldots , x'_{d}
\rangle$, we have
$
\SSS = \SSS_{1} \times \SSS_{2}$.
Let $\tilde{o}$
denote the point of $\tilde{\SSS_{1}} = 
\Spec \tilde{\KK} [x'_{0} ,  \ldots , x'_{r} ] / (x'_{0}   \cdots  x'_{r} )$ 
defined by
$x'_{0} = \cdots = x'_{r} = 0$.

The following proposition is due to Gubler.

\begin{citeProposition}[Proposition~5.2 of 
\cite{gubler4}] \label{gubler4Proposition5.2}
Any formal open covering of $\XXX'$
admits a refinement $\{ \UUU' \}$
by formal open subsets $\UUU'$ with \'etale morphisms
as in (\ref{ssetale})
which have the following properties.
\begin{enumerate}
\renewcommand{\labelenumi}{(\alph{enumi})}
\item
Any $\UUU'$ is a formal affine open subscheme of $\XXX'$.
\item
There exists a distinguished stratum $S$ of $\tilde{\XXX'}$
associated to $\UUU'$ such that for any stratum $T$ of $\tilde{\XXX'}$,
we have $S \subset \overline{T}$ if and only if
$\tilde{\UUU'} \cap \overline{T} \neq \emptyset$,
where $\overline{T}$ is the closure of $T$ in $\tilde{\XXX'}$.
\item
The subset
$\tilde{\psi}^{-1} ( \{ \tilde{o} \} \times \tilde{\SSS_{2}})$
is the stratum of $\tilde{\UUU'}$
which is equal to $\tilde{\UUU'} \cap S$
for the distinguished stratum $S$ from (b).
\item
Any stratum of $\tilde{\XXX'}$
is the distinguished stratum of a suitable $\UUU'$.
\end{enumerate}
\end{citeProposition}

We can define a subspace $S (\XXX')$ of $X' = (\XXX')^{\an}$
called the
\emph{skeleton}.
It has a canonical structure of an abstract
simplicial set which reflects the incidence relations
between the strata of $\tilde{\XXX'}$.
We briefly recall some properties of skeletons here,
and refer to \cite{berkovich4} or \cite[5.3]{gubler4} for more details.

First, we recall the skeletons
$S(\SSS_{1})$ and $S(\SSS)$.
We set $\GG_{m}^{r} := \Spec \KK [(x'_{1})^{\pm} , \ldots , (x'_{r})^{\pm}]$,
the algebraic torus with the standard coordinates $x'_{1} , \ldots , x'_{r}$,
and consider the associated analytic group $(\GG_{m}^{r})^{\an}$.
Let $\val' : (\GG_{m}^{r})^{\an} \to \RR^{r}$
denote the valuation map with respect to the
coordinates $x'_{1} , \ldots , x'_{r}$
(cf. \S~\ref{tori}).\footnote{We
write here $\val'$ instead of $\val$ to emphasize
that it is the valuation map with respect to the 
coordinates $x'_{1} , \ldots , x'_{r}$.}
We regard $\SSS_{1}^{\an}$ as a rational subdomain of $(\GG_{m}^{r})^{\an}$
by omitting $x'_{0}$.
We put $\Delta' :=
\{
(u'_{1} , \ldots , u'_{r})
\in \RR_{\geq 0}^{r}
\mid
u'_{1} + \cdots + u'_{r} \leq v (\pi)
\}$.
Then $\val'$ induces a map
$\SSS_{1}^{\an} \to \Delta'$ by restriction.
It is known that
this map restricts to an isomorphism $S ( \SSS_{1} ) \cong \Delta'$,
and that
the first projection $\SSS \to \SSS_{1}$ restricts to 
an isomorphism $S (\SSS ) \cong S (\SSS_{1} )$.
Thus we have identifications $S (\SSS ) = S (\SSS_{1} ) = \Delta'$.

Recall that the vertices of $\Delta'$
correspond to the irreducible components of 
$\tilde{\SSS}$.
In fact,
let $v_{0}$ denote the origin in $\Delta'$
and
let $v_{i}$ denote the vertex of
$\Delta'$
whose $i$-th coordinate is
$v (\pi)$ ($i = 0 , \ldots , r$).
The points $v_{0} , v_{1} , \ldots , v_{r}$ are the vertices of $\Delta'$.
Let $\xi_{i} \in S (\SSS)$ be the point corresponding to
the irreducible component
of $\Spec \tilde{\KK} [ x_{0}' , \ldots , x_{d}' ] / (x_{0}' \cdots x_{r}')$
defined by $x_{i} = 0$ for each $i = 0 , \ldots , r$
(cf. \S~\ref{bersp}).
Then we have $\xi_{i} = v_{i}$ via the identification
$S (\SSS ) = \Delta'$.

Let $S$ be a
stratum of $\tilde{\XXX'}$ of codimension $r$.
We
take a formal affine open subset $\UUU'$
of $\XXX'$ such that $S$ is the distinguished stratum of $\UUU'$
in Proposition~\ref{gubler4Proposition5.2},
and an \'etale morphism as in (\ref{ssetale}).
The skeleton $S (\UUU')$ of $( \UUU')^{\an}$
is a subset of $S (\XXX')$, and 
it is known that $\psi$ in (\ref{ssetale})
induces an isomorphism
$S (\UUU) \cong S(\SSS)$ between skeletons.
Thus
\addtocounter{Claim}{1}
\begin{align} \label{skelton-simplex}
S(\UUU') \cong S (\SSS) = S(\SSS_{1}) = \Delta'
\end{align}
and, in particular,
the subset $S(\UUU')$ of $S (\XXX')$
is a simplex. 
It is also known that
$S(\UUU')$ depends only on $S$ and not
on the choice of $\UUU'$,
so that
we write $\Delta_{S} = S(\UUU')$ and call it the \emph{canonical simplex}
corresponding to $S$.
The canonical simplices $\{ \Delta_{S} \}_{S \in \str ( \tilde{\XXX'})}$
cover $S(\XXX')$,
which gives
a canonical structure of an abstract
simplicial set to the skeleton $S(\XXX')$.

We have a continuous map $\Val : X' \to S (\XXX')$
which restricts to the identity on $S(\XXX')$.
If $S$ is a distinguished stratum of $\tilde{\XXX'}$ 
associated to $\UUU'$ in the sense of Proposition~\ref{gubler4Proposition5.2},
then the restriction $(\UUU')^{\an} \to \Delta_{S}$
of $\Val$ to $(\UUU')^{\an}$
is described as follows:
regarding $\SSS_{1}^{\an}$ as a rational subdomain of $(\GG_{m}^{r})^{\an}$
and
using the identification 
\addtocounter{Claim}{1}
\begin{align} \label{canonicalsimplexidentification}
\Delta_{S} =
\Delta' 
\end{align}
given by
(\ref{skelton-simplex}),
we can describe $\Val$ as $\Val (p) = 
\val' ( \psi^{\an} (p)) 
$
for $p \in (\UUU')^{\an}$.

Let $\DDD$ be a $\Gamma$-rational subdivision of the skeleton $S ( \XXX' )$.
This means that $\DDD$ is a family of polytopes,
each contained in a canonical simplex, such that
$\{ \Delta \in \DDD \mid \Delta \subset \Delta_{S}\}$
is a 
$\Gamma$-rational 
polytopal decomposition of $\Delta_{S}$ for any stratum $S$
of $\tilde{\XXX'}$
(cf. \cite[5.4]{gubler4}).

\begin{Remark} \label{rem:Proposition5.5gubler4}
By
\cite[Proposition~5.5]{gubler4},
we have
a unique formal model $\XXX''$ of $(\XXX')^{\an}$ associated to $\DDD$.
Note that there is a morphism $\iota' : \XXX'' \to \XXX'$
extending the identity on the generic fiber $(\XXX')^{\an}$.
\end{Remark}

Although we omit the precise construction of $\XXX''$,
we remind two important propositions concerning $\XXX''$.

\begin{citeProposition}[Proposition~5.7
of \cite{gubler4}] \label{gubler4Proposition5.7}
Let $\XXX''$ be the formal scheme associated to $\DDD$
as in Remark~\ref{rem:Proposition5.5gubler4}.
Let $\red_{\XXX''} : (\XXX'')^{\an} \to \tilde{\XXX''}$
be the reduction map.
Then there exists a bijective correspondence between
open faces $\tau$ of $\DDD$ and strata $R$ of $\tilde{\XXX''}$
given by
\[
R = \red_{\XXX''} ( \Val^{-1} ( \tau))
,
\quad
\tau = \Val ( \red_{\XXX''}^{-1} (Y))
,
\]
where $Y$ is any non-empty subset of $R$.
\end{citeProposition}

For a canonical simplex $\Delta_{S}$,
let $\relin \Delta_{S}$ denote the relative interior of $\Delta_{S}$.\footnote{
For a polytope $P$, let $\relin (P)$
denote the relative interior of $P$ in this article.}

\begin{citeProposition} [cf. Corollary~5.9 of \cite{gubler4}] 
\label{gubler4Corollary5.9}
Let $\XXX''$ be the formal scheme associated to $\DDD$
and let $\iota' : \XXX'' \to \XXX'$ be the morphism extending the
identity on $(\XXX')^{\an}$ (cf. Remark~\ref{rem:Proposition5.5gubler4}).
Let $u \in \DDD$ be a vertex and let $R$ be the stratum of $\tilde{\XXX''}$
corresponding to $u$ in Proposition~\ref{gubler4Proposition5.7}.
Then $S := \tilde{\iota'} (R)$ is a unique stratum of $\tilde{\XXX'}$
with $u \in \relin \Delta_{S}$.
Furthermore, if we put $r := \dim R - \dim S = \dim \Delta_S$,
then $\tilde{\iota'} |_{R} : R \to S$ has a structure
of $(\GG_{m}^{r})_{\tilde{\KK}}$-torsor.
\end{citeProposition}


Finally in this subsection, we show a lemma.

\begin{Lemma} \label{vertex-genericpoint}
Let $\DDD$ and $\XXX''$ be as above.
Let $u$ be a vertex of $\DDD$
of $S ( \XXX' )$
and let $R$ be the stratum of $\tilde{\XXX''}$ corresponding to $u$.
Let $\xi_{R} \in (\XXX'')^{\an}$ be the point corresponding to
the generic point of $R$
(cf. \S~\ref{bersp}).
Then $u = \xi_{R}$.
\end{Lemma}

\begin{Pf}
It follows from \cite[Proposition 5.7 and Corollary 5.9 (a)]{gubler4}
that $u = \Val (\xi_{R})$.
Since $\Val$ restricts to the identity on $S ( \XXX')$,
it only remains to show $\xi_{R} \in S ( \XXX'')$,
but it is done in the proof of
\cite[Corollary 5.9 (g)]{gubler4}.
\end{Pf}

\subsection{Construction of  the action
and the compatibility lemma}
\label{construction-action}

Let $S$ be a
stratum of $\tilde{\XXX'}$ of codimension $r$.
We
take a formal affine open subset $\UUU'$
of $\XXX'$ such that $S$ is the distinguished stratum of $\UUU'$
in Proposition~\ref{gubler4Proposition5.2},
and an \'etale morphism as in (\ref{ssetale}).
Let $\tilde{p} \in S \cap \tilde{\UUU'}$ be a closed point.
In this subsection, we define an affinoid-torus-action of the formal
fiber $X'_{+} (\tilde{p})$ over $\tilde{p}$.
Since 
our interest is local at 
$\tilde{p}$,
we may and do assume
$\XXX' = \UUU'
$.\footnote{We have $X'_{+} (\tilde{p}) = 
(\UUU')^{\an}_{+} (\tilde{p})$ in fact.}

Recall from the previous subsection that $\SSS = \SSS_1 \times \SSS_2$.
Since
$(\SSS_1)^{\an} \subset (\GG_{m}^{r})^{\an}$ and
$(\SSS_1)^{\an}$ is stable under the $(\GG_{m}^{r})^{\an}_{1}$-action,
we have a $(\GG_{m}^{r})^{\an}_{1}$-action on $(\SSS_1)^{\an}$.
We put the trivial $(\GG_{m}^{r})^{\an}_{1}$-action on $( \SSS_2 )^{\an}$.
Then, via the isomorphism
$\SSS^{\an} = (\SSS_1)^{\an} 
\times ( \SSS_2 )^{\an}$, we have a 
$(\GG_{m}^{r})^{\an}_{1}$-action on $\SSS^{\an}$.
Let $(\SSS^{\an})_{+} ( \tilde{\psi}(\tilde{p})) :=
\red_{\SSS}^{-1} ( \tilde{\psi}(\tilde{p}) )$
the formal fiber of $\SSS^{\an}$ over $\tilde{\psi}(\tilde{p})$.
Since 
$\tilde{\psi} ( \tilde{p} ) = ( \tilde{o} , \tilde{c})$
for some $\tilde{c} \in \tilde{\SSS_2}$
via the identification $\tilde{\SSS} = \tilde{\SSS_1} 
\times \tilde{ \SSS_2 }$, 
we see that
$(\SSS^{\an})_{+} (\tilde{\psi}(\tilde{p}))$ is stable under the
$(\GG_{m}^{r})^{\an}_{1}$-action on $\SSS^{\an}$,
and thus
$(\GG_{m}^{r})^{\an}_{1}$ acts on
$(\SSS^{\an})_{+} (\tilde{\psi}(\tilde{p}))$.

We
define the $(\GG_{m}^{r})^{\an}_{1}$-action on $X'_{+}  (\tilde{p})$:
By
\cite[Lemma 4.4]{berkovich4} or \cite[Proposition 2.9]{gubler4},
the restricted morphism
\addtocounter{Claim}{1}
\begin{align} \label{morphisms-of-formalfibers}
\psi'^{\an} |_{X'_{+}  (\tilde{p})}
: X'_{+}  (\tilde{p}) \to (\SSS^{\an})_{+} (\tilde{\psi}(\tilde{p}))
\end{align}
is an isomorphism,
and then we
define
the $(\GG_{m}^{r})^{\an}_{1}$-action on $X'_{+} (\tilde{p})$
by pulling back
the $(\GG_{m}^{r})^{\an}_{1}$-action on
$(\SSS^{\an})_{+} (\tilde{\psi}(\tilde{p}))$ via this isomorphism.

Let $\Delta_{S}$ be the canonical simplex corresponding to $S$.
Let $\DDD$ be a $\Gamma$-rational polytopal subdivision of 
the skeleton $S ( \XXX' )$,
$\XXX''$ the formal model of $(\XXX')^{\an}$ associated to $\DDD$,
$\iota' : \XXX'' \to \XXX'$ the morphism extending
the identity on $X'$
(cf. Remark~\ref{rem:Proposition5.5gubler4}),
and let 
$\red_{\XXX''} : X' \to \tilde{\XXX''}$
be the reduction map.
Suppose that $u \in \DDD$ is a vertex with $u \in \relin \Delta_{S}$
and
let $R$ be the stratum of $\tilde{\XXX''}$ corresponding to $u$
in Proposition~\ref{gubler4Proposition5.7}.
Recall that we have 
a $(\GG_{m}^{r})_{\tilde{\KK}}$-torsor
$\tilde{\iota'} : R \to S$ (cf. Proposition~\ref{gubler4Corollary5.9}).
If $\tilde{p} \in S$ is a closed point,
then $(\GG_{m}^{r})_{\tilde{\KK}}$ acts on the fiber $\{ \tilde{p} \} \times_{S} R$.

\begin{Remark} \label{newRemark1}
We have $\red_{\XXX'} = \tilde{\iota'} \circ \red_{\XXX''}$
and $\tilde{\iota'} ( \{ \tilde{p} \} \times_{S} R ) = \{ \tilde{p} \}$.
It follows that
$
(\red_{\XXX'})^{-1} ( \tilde{p} )
\supset
(\red_{\XXX''})^{-1} 
( \{ \tilde{p} \} \times_{S} R )
$,
and thus
\[
\red_{\XXX''} ( X'_{+}  (\tilde{p})  )
=
\red_{\XXX''} ( (\red_{\XXX'})^{-1} (\tilde{p})  )
\supset
\{ \tilde{p} \} \times_{S} R
.
\]
\end{Remark}

Next we show that the
$(\GG_{m}^{r})^{\an}_{1}$-action on $X'_{+} (\tilde{p})$
defined above
is compatible with
the
action of $(\GG_{m}^{r})_{\tilde{\KK}}$ on $\{ \tilde{p} \} \times_{S} R$
with respect to the reduction map $\red_{\XXX''}$.
To be precise, 
we show the following lemma.

\begin{Lemma} \label{torus-action-ff}
With the notation above,
we have a commutative diagram
\begin{align*}
\begin{CD}
(\GG_{m}^{r})^{\an}_{1} 
\times 
(\red_{\XXX''} )^{-1}(\{ \tilde{p} \} \times_{S} R)
@>{}>>
(\red_{\XXX''} )^{-1}
(\{ \tilde{p} \} \times_{S} R)
\\
@V{\red_{(\GG_m^r)_{1}^{\mathrm{f-sch}} \times \XXX''}}VV @VV{\red_{\XXX''}}V \\
(\GG_{m}^{r})_{\tilde{\KK}} \times (\{ \tilde{p} \} \times_{S} R)
@>>>
\{ \tilde{p} \} \times_{S} R
,
\end{CD}
\end{align*}
where 
the first row is the restriction of the
$(\GG_{m}^{r})^{\an}_{1}$-action on $X'_{+} (\tilde{p})$,
and the second
row is the $(\GG_{m}^{r})_{\tilde{\KK}}$-action
on $\{ \tilde{p} \} \times_{S} R$
(cf. Remark~\ref{newRemark1}).
\end{Lemma}

\begin{Pf}
Recall that we have an isomorphism
$ p_{1}^{\an} \circ \psi^{\an} : S (\XXX') \cong S ( \SSS_{1} )$
(cf. (\ref{skelton-simplex})),
where 
$p_{1} : \SSS \to \SSS_{1}$ is the canonical projection.
Let
$\DDD_{1}$ be the subdivision of $S ( \SSS_{1} )$
induced from
$\DDD$ via the isomorphism
$ p_{1}^{\an} \circ \psi^{\an}$ between the skeletons,
and let
$\SSS_{1}'$ be the formal scheme corresponding to the subdivision
$\DDD_{1}$ of $S ( \SSS_{1} )$.
As in \cite[Remark 5.6]{gubler4},
we have a cartesian diagram 
\begin{align*} 
\begin{CD}
\XXX'' @>{\psi'}>> \SSS' 
\\
@V{\iota'}VV @VV{\iota}V 
\\
\XXX' @>{\psi}>> \SSS 
,
\end{CD}
\end{align*}
where $\iota$ is the morphism 
$\mathscr{S}_1' \times \mathscr{S}_2 \to \mathscr{S}_1 \times \mathscr{S}_2$
given by the base-change of the
morphism $\mathscr{S}_1' \to \mathscr{S}_1$
arising from the subdivision $\mathscr{D}_1$.

We set $u_{1} := p_{1}^{\an} ( \psi^{\an} (u))$,
which is a vertex of $\DDD_{1}$ 
and is in the interior of $\Delta' = S(\SSS_{1})$.
Let
$T_{1}'$ be the stratum of $\tilde{\SSS'_{1}}$ corresponding to $u_{1}$.
Recall that $\tilde{\psi} \left( \tilde{p} \right) =
\left( \tilde{o} , \tilde{c} \right)$.
From Remark~\ref{newRemark1}
and the proof of \cite[Proposition~5.7]{gubler4},
we obtain a commutative diagram
\addtocounter{Claim}{1}
\begin{align} \label{eq:CDnewbyreferee}
\begin{CD}
X'_{+} \left( \tilde{p} \right) @<{\supset}<<
\left( \red_{\mathscr{X}''} \right)^{-1}
\left( \{ \tilde{p} \} \times_{S} R \right) 
@>{\red_{\mathscr{X}''}}>>
\{ \tilde{p} \} \times_{S} R \\
@V{\psi^{\an} |_{X'_{+} \left( \tilde{p} \right)}}VV 
@VVV @VV{\tilde{\psi'}|_{\{ \tilde{p} \} \times_{S} R}}V
\\
\mathscr{S}'_{+} \left( \tilde{\psi} \left( \tilde{p} \right) \right)
@<{\supset}<<
\left( \red_{\mathscr{S}'} \right)^{-1}
\left( T_1' \times \left\{ \tilde{c} \right\} \right) 
@>{\red_{\mathscr{S}'}}>>
T_1' \times \left\{ \tilde{c} \right\}
,
\end{CD}
\end{align}
where the first column is the isomorphism (\ref{morphisms-of-formalfibers}),
and the left rows are inclusions.

By
(17) in the proof of \cite[Proposition~5.7]{gubler4},
we have
$\tilde{\psi '}^{-1} \left( T_1 ' \times \tilde{\mathscr{S}_2} \right)
=
R
$.
It follows that
the last column $\{ \tilde{p} \} \times_{S} R \to
T_1' \times \left\{ \tilde{c} \right\}$ 
in (\ref{eq:CDnewbyreferee})
is an isomorphism.
Furthermore,
it follows that the middle column is also an isomorphism.

Noting (14) in the proof of \cite[Proposition~5.7]{gubler4},
we find that $\left( \red_{\mathscr{S}'} \right)^{-1}
\left( T_1' \times \left\{ \tilde{c} \right\} \right) $ is 
$(\GG_m^{r})^{\an}_1$-invariant. It follows that 
$\left( \red_{\mathscr{X}''} \right)^{-1}
\left( \{ \tilde{p} \} \times_{S} R \right) $ is $(\GG_m^{r})^{\an}_1$-invariant.
Since the isomorphism in the middle (resp. right) column
is $(\GG_m^{r})^{\an}_1$-equivariant 
(resp. $(\GG_m^{r})_{\tilde{\KK}}$-equivariant)
and since $\red_{\mathscr{S}'}$ is 
equivariant with respect to the reduction map 
$(\GG_m^{r})^{\an}_1 \to (\GG_m^{r})_{\tilde{\KK}}$,
the right upper row is also equivariant with respect to the reduction map 
$(\GG_m^{r})^{\an}_1 \to (\GG_m^{r})_{\tilde{\KK}}$.
This concludes that
Lemma~\ref{torus-action-ff} holds.
\end{Pf}


\begin{Remark} \label{nd4}
The $(\GG_{m}^{r})^{\an}_{1}$-action on $X'_{+} (\tilde{p})$ preserves
$\left( \red_{\XXX''} |_{X'_{+} (\tilde{p})} \right)^{-1} \left(
\{ \tilde{p} \} \times_S R \right)$,
as is shown above.
\end{Remark}

\section{Raynaud extension and Mumford models} \label{RaynaudextMumford}

In this section, we
put together some properties
of the Raynaud extensions of abelian
varieties and their Mumford models, which will be used later.
We
recall,
in \S~\ref{RtM}
and \S~\ref{RtM2}, basic facts
discussed in
\cite[\S 4]{gubler4}.
We show some properties
of the torus rank and the abelian rank of an abelian variety
in \S~\ref{homo-Raynaud}.

\subsection{Raynaud extension and the valuation map} 
\label{RtM}
We recall some notions on Raynaud extensions.
See \cite[\S 1]{BL} and \cite[\S 4]{gubler4} for details.

Let $A$ be an abelian variety over $\KK$.
By \cite[Theorem 1.1]{BL},
there exists 
a unique analytic subgroup
$A^{\circ} \subset A^{\an}$
with a formal model $\mathscr{A}^{\circ}$ 
having the following properties:
\begin{itemize}
\item
$\mathscr{A}^{\circ}$ is a formal group scheme and
$(\mathscr{A}^{\circ})^{\an} \cong A^{\circ}$ as analytic groups;
\item
There is 
a short exact sequence
\addtocounter{Claim}{1}
\begin{align} \label{qcptformalR}
\begin{CD}
1 @>>> \mathscr{T}^{\circ} @>>> \mathscr{A}^{\circ} @>>>
\mathscr{B} @>>> 0 ,
\end{CD}
\end{align}
where $\mathscr{T}^{\circ} \cong (\GG_{m}^{n})^{\mathrm{f-sch}}_{1}$ 
for some $n \geq 0$, and
$\mathscr{B}$ is a formal abelian variety over $\KK^{\circ}$.
\end{itemize}
By 
\cite[Satz 1.1]{bosch},
such an
$\AAA^{\circ}$ is unique,
and $\TTT^{\circ}$
and $\BBB$ are also uniquely determined.
Taking the generic fiber of (\ref{qcptformalR}),
we obtain an exact sequence
\addtocounter{Claim}{1}
\begin{align} \label{Raynaudextcirc}
\begin{CD}
1 @>>> T^{\circ} @>>> A^{\circ} @>{(q^{\circ})^{\an}}>> B @>>> 0
\end{CD}
\end{align}
of analytic groups.
We call $\TTT^{\circ}$, $\AAA^{\circ}$ and $\BBB$ 
the \emph{canonical formal models}
of $T^{\circ}$, $A^{\circ}$ and $B$ respectively.

Since $T^{\circ} \cong (\GG_{m}^{n})^{\an}_{1}$,
we regard
$T^{\circ}$ as an analytic subgroup
of the analytic torus 
$T = (\GG_{m}^{n})^{\an}$.
Pushing (\ref{Raynaudextcirc}) out by $T^{\circ} \hookrightarrow T$,
we obtain an exact sequence
\addtocounter{Claim}{1}
\begin{align} \label{Raynaudext}
\begin{CD}
1 @>>> T @>>> E @>{q^{\an}}>> B @>>> 0 ,
\end{CD}
\end{align}
which we call the \emph{Raynaud extension} of $A$.
The natural morphism $A^{\circ} \to E$
is an immersion of analytic groups.
The assertion
\cite[Theorem 1.2]{BL} says that
the homomorphism $T^{\circ} \hookrightarrow A^{\an}$ extends uniquely
to a homomorphism $T \to A^{\an}$ and hence to 
a homomorphism
$p^{\an} : E \to A^{\an}$.
It is known that
$p^{\an}$ is a surjective homomorphism and 
$M := \Ker p^{\an}$ is a lattice in $E ( \KK )$.
Thus $A^{\an}$ can be described as a quotient of
$E$ by a lattice.
This $p^{\an}$
is called the \emph{uniformization}
of $A$.

The dimension of $T$ is called the \emph{torus rank} of $A$,
and the dimension of $B$ is called \emph{abelian rank} of $A$.
We denote the abelian rank of $A$ by $b (A)$.
Note that the torus rank of $A$ equals $\dim A - b (A)$.
The abelian variety $A$ is said to be \emph{degenerate} if $b(A) < \dim A$,
or equivalently if the torus rank of $A$ is positive.
Note that ``being non-degenerate'' means ``having good reduction''.

We can take 
transition functions of the
$T$-torsor (\ref{Raynaudext}) valued in $T^{\circ}$,
and thus
we can define a continuous map
\addtocounter{Claim}{1}
\begin{align} \label{valuationmap}
\val : E \to \RR^{n}
\end{align}
as in \cite{BL},
where $n$
is the torus rank of $A$.
We recall here how it is constructed.
Fix an isomorphism
$T^{\circ} \cong (\GG_{m}^{n})^{\an}_{1}$,
with the standard coordinates $x_{1}, \ldots , x_{n}$.
We can
take a
covering $\{ V \}$ of $B$ consisting of 
rational subdomains
with trivializations
\addtocounter{Claim}{1}
\begin{align} \label{trivialization}
((q^{\circ})^{\an})^{-1} (V)
\cong V \times (\GG_{m}^{n})^{\an}_{1} 
\end{align}
as $(\GG_{m}^{n})^{\an}_{1}$-torsors
for all $V$.
Since the Raynaud extension is the push-out of
(\ref{Raynaudextcirc}) by the canonical inclusion 
$(\GG_{m}^{n})^{\an}_{1} \hookrightarrow (\GG_{m}^{n})^{\an}$,
the isomorphisms (\ref{trivialization}) extend to
isomorphisms
\addtocounter{Claim}{1}
\begin{align} \label{triv-1}
(q^{\an})^{-1} (V)
\cong V \times (\GG_{m}^{n})^{\an}
\end{align}
of $(\GG_{m}^{n})^{\an}$-torsors.
Thus
we obtain morphisms
\addtocounter{Claim}{1}
\begin{align} \label{localval}
r_{V} : (q^{\an})^{-1} (V) \cong V \times (\GG_{m}^{n})^{\an}
\to (\GG_{m}^{n})^{\an}
\end{align}
for all $V$
by composing with the second projection.
A different choice of (\ref{trivialization}) gives a different
isomorphism in
(\ref{triv-1}) and hence
a different morphism in (\ref{localval})
for each $V$,
but the difference is only the multiplication of an 
element of $(\GG_{m}^{n})^{\an}_{1}$.
Therefore, 
the morphisms $(q^{\an})^{-1} (V) \to \RR^{n}$ given by
\[
e \mapsto
(
-\log r_{V} (e) (x_{1})
,
\ldots
,
- \log r_{V} (e) (x_{n})
)
\]
patch together to be a morphism from $E$ to $\RR^{n}$.
This is our valuation map $\val : E \to \RR^{n}$.

We set $\Lambda := \val (M) \subset \RR^{n}$,
where we recall that $M = \Ker p^{\an}$ is a lattice of $E$.
Then $\Lambda$
is a 
complete lattice in 
$\RR^{n}$
and
is
contained in $\Gamma^{n}$.
There is a diagram
\begin{align*}
\begin{CD}
E @>{\val}>> \RR^{n} \\
@VVV @VVV \\
A^{\an} @>{\overline{\val}}>> \RR^{n} / \Lambda
\end{CD}
\end{align*}
that commutes.
The homomorphism $\overline{\val}$
is also called the \emph{valuation map}.
It follows from
the constructions 
of $\val$ and $\overline{\val}$
that
$A^{\circ} = \overline{\val}^{-1} ( \overline{\mathbf{0}} ) \cong
\val^{-1} (\mathbf{0})$.

\subsection{Homomorphism, products and abelian ranks} \label{homo-Raynaud}

Let $A_{1}$ and $A_{2}$ be abelian varieties over $\KK$
and let 
\[
\begin{CD}
1 @>>> T_{i} @>>> E_{i} @>{q_{i}^{\an}}>> B_{i} @>>> 0
\end{CD}
\]
be the Raynaud extension of $A_{i}$ for $i = 1,2$.
We consider a homomorphism
$\phi : A_{1} \to A_{2}$.
Then
\cite[Proposition 2.2]{BL-neron}
tells us that $\phi$ induces a homomorphism 
$A_{1}^{\circ} \to A_{2}^{\circ}$
and hence a homomorphism $T_{1}^{\circ} \to T_{2}^{\circ}$.
Thus we have an induced homomorphism
$\phi_{ab} : B_{1} \to B_{2}$ between the abelian parts.
Furthermore, it follows from \cite[Theorem~1.2]{BL}
that $\phi$ gives rise to a homomorphism
$\Phi : E_{1} \to E_{2}$.

It follows from the construction of the
valuation map that 
$\Phi$ descends to
a linear map
$\phi_{\aff} : \RR^{n_{1}} \to \RR^{n_{2}}$ via the valuation maps,
where 
$n_{1}$ and $n_{2}$ are the torus ranks of $A_{1}$
and $A_{2}$ respectively.
Further,
$\phi_{\aff} ( \Lambda_{1} ) \subset \Lambda_{2}$
holds.
Thus we obtain a homomorphism $\overline{\phi}_{\aff}
: \RR^{n_{1}} / \Lambda_{1} \to  \RR^{n_{2}} / \Lambda_{2}$ 
which makes the  diagram
\addtocounter{Claim}{1}
\begin{align} \label{tamanidetekuru}
\begin{CD}
A_{1}^{\an} @>{\phi}>> A^{\an}_{2} \\
@VVV @VVV \\
\RR^{n_{1}} / \Lambda_{1} @>{\overline{\phi}_{\aff}}>> \RR^{n_{2}} / \Lambda_{2}
\end{CD}
\end{align}
commutative.

Next we consider the direct product.
We set $A := A_{1} \times A_{2}$.
The Raynaud extensions of $A_{1}$ and $A_{2}$ gives rise to
an exact sequence
\[
1 \to T_{1}^{\circ} \times T_{2}^{\circ} 
\to
A_{1}^{\circ} \times A_{2}^{\circ}
\to
B_{1} \times B_{2}
\to
0
.
\]
It follows
from their definitions
that
$A^{\circ} = A_{1}^{\circ} \times A_{2}^{\circ}$
and
$T^{\circ} = T_{1}^{\circ} \times T_{2}^{\circ}$
hold,
and
that
\addtocounter{Claim}{1}
\begin{align} \label{product-Raynaud}
1 \to T_{1} \times T_{2} \to E_{1} \times E_{2}
\to B_{1} \times B_{2} \to 0
\end{align}
is the Raynaud extension of $A$.
Further,
\[
\overline{\val} : 
(A_{1} \times A_{2} )^{\an} \to \RR^{n_{1} + n_{2}} / 
(\Lambda_{1} \oplus \Lambda_{2})
\]
coincides with the map 
\[
 A_{1}^{\an} \times  A_{2}^{\an} \to 
\RR^{n_{1}} / \Lambda_{1} \times \RR^{n_{2}} / \Lambda_{2}
\]
given by the product of
$\overline{\val}_{1} : A_{1}^{\an} \to \RR^{n_{1}} / \Lambda_{1}$
and
$\overline{\val}_{2} : A_{2}^{\an} \to \RR^{n_{2}} / \Lambda_{2}$.

We show some properties of abelian ranks
which will be used later.

\begin{Lemma} \label{isogeny-abelianrank}
Let $A_{1}$ and $A_{2}$ be abelian varieties over $\KK$.
\begin{enumerate}
\item
We have $b (A_{1} \times A_{2}) = b (A_{1}) + b ( A_{2} )$.
\item
Suppose that  $\phi : A_{1} \to A_{2}$ is an isogeny.
Then $b (A_{1} ) = b (A_{2})$.
\end{enumerate}
\end{Lemma}

\begin{Pf}
The equality in
(1) follows from the fact that
(\ref{product-Raynaud}) is the Raynaud extension of $A_{1} \times A_{2}$.

To show (2),
note that
there is an isogeny
$A_{2} \to A_{1}$
as well as an isogeny $A_{1} \to A_{2}$
(cf. \cite[p.157, Remark]{mumford}).
Thus it is enough to show $b( A_{1}) \geq b(A_{2})$,
so that it suffices to
show that
the homomorphism $\phi_{ab} : B_{1} \to B_{2}$ 
induced from $\phi$ is surjective.

Let $y \in B_{2}$ be a point.
Since $q_{2}^{\an} |_{A_{2}^{\circ}} : A_{2}^{\circ} \to B_{2}$,
$\phi^{\an} : A_{1}^{\an} \to A_{2}^{\an}$
and $p_{1} : E_{1} \to A_{1}^{\an}$ are surjective,
there exists a point $x \in E_{1}$ with 
$q_{2}^{\an} (\phi^{\an} ( p_{1} (x) )) = y$.
Then  we have $\phi_{ab} (q_{1}^{\an} (x)) = y$,
which implies $\phi_{ab}$ is surjective.
\end{Pf}

\begin{Remark} \label{valofisogeny}
If $\phi : A_{1} \to A_{2}$ is an isogeny,
then
$\phi_{\aff}$ is an isomorphism of vector spaces,
and thus
$\overline{\phi}_{\aff}$ is a finite surjective homomorphism.
\end{Remark}

\begin{Proposition} \label{exactsequence-abelianrank}
Let
\[
0 \to A_{1} \to A_{2} \to A_{3} \to 0
\]
be an exact sequence of abelian varieties over $\KK$.
Then we have $b (A_{2}) = b (A_{1}) + b (A_{3})$
and $n_{2}  = n_{1} + n_{3}$,
where $n_{i}$ denote the torus rank of $A_{i}$.
\end{Proposition}

\begin{Pf}
Recall that Poincar\'e's complete reducibility theorem 
(cf. \cite[Proposition~12.1]{milne} or \cite[\S~19, Theorem~1]{mumford}) gives us an
abelian subvariety $A'$ of $A_2$ such that 
the natural homomorphism
$A_{1} \times A' \to A_{2}$ given by $(a_{1} , a') \mapsto a_{1} + a'$ 
is an isogeny.
Then $b (A_{2}) = b (A_{1} \times A') = b (A_{1}) + b (A') $ by
Lemma~\ref{isogeny-abelianrank}.
Since
the composite $A' \hookrightarrow A_{2} \to A_{3}$
is an isogeny,
we have $b (A') = b(A_{3})$
by Lemma~\ref{isogeny-abelianrank}~(2).
It follows that
$b (A_{2}) = b (A_{1}) + b(A_{3})$,
which shows one equality of this proposition.
The equality for torus ranks follows from
that
for abelian ranks.
\end{Pf}

\subsection{Mumford models, Torus-torsors,
and initial degenerations} \label{RtM2}

In this subsection, we recall some properties of Mumford models
associated to a $\Lambda$-periodic $\Gamma$-rational polytopal decomposition.
The basic reference is \cite[\S~4]{gubler4}.
We also define the initial degeneration of a closed subvariety
of an abelian variety.

Let $A$ be an abelian variety over $\KK$,
and let
\[
\begin{CD}
1 @>>> T @>>> E @>{q^{\an}}>> B @>>> 0
\end{CD}
\]
be the Raynaud extension of $A$.
We use the notation in \S~\ref{RtM}:
The affinoid torus
$T^{\circ}$ is the maximal affinoid subtorus of $T$;
The morphism
$p^{\an} : E \to A^{\an}$ is the uniformization,
the map $\val : E \to \RR^{n}$
is the valuation map,
where $n$ is the torus rank of $A$;
Further,
$\Lambda := \val (\Ker p^{\an})$ is the lattice in $\RR^{n}$,
and $\overline{\val} : A \to \RR^{n} / \Lambda$
is the valuation map induced from $\val$ by quotient.

Let 
$\CCC$ be a $\Lambda$-periodic 
$\Gamma$-rational
polytopal decomposition of 
$\RR^{n}$ (cf. \S~\ref{conv-term}).
For a polytope 
${\Delta} \in {\CCC}$,
the subset ${\val}^{-1} ( \Delta) \subset E
$
is an analytic subdomain,
and there exists a natural surjective morphism
\[
q^{\an}_{\Delta} := q^{\an}|_{{\val}^{-1} ( \Delta)} : 
{\val}^{-1} ( \Delta) \to B
.
\]
Since $\val$
is invariant under the action
of $T^{\circ}$,
we have a natural $T^{\circ}$-action on ${\val}^{-1} ( \Delta)$,
which is an action
over $B$
with respect to $q_{\Delta}^{\an}$.

Let 
$\overline{\CCC}$ denote the polytopal decomposition of 
$\RR^{n} / \Lambda$ induced from $\CCC$ by quotient
and let
$\overline{\Delta} \in \overline{\CCC}$ be a polytope.
Then $\overline{\val}^{-1} ( \overline{\Delta})$ 
is an analytic subdomain of $A^{\an}$
with a $T^{\circ}$-action. 
Fix a representative
$\Delta \in {\CCC}$ of $\overline{\Delta}$.
The
uniformization map $p^{\an} : E \to A^{\an}$ restricts to an isomorphism
$
p^{\an} |_{{\val}^{-1} (\Delta)} :
{\val}^{-1} (\Delta) \to
\overline{\val}^{-1} ( \overline{\Delta})
$,
via which
we define
\begin{align*} 
\overline{q^{\an}_{\Delta}} := 
q^{\an}_{\Delta} \circ ( p^{\an} |_{{\val}^{-1} (\Delta)} )^{-1}
:
\overline{\val}^{-1} ( \overline{\Delta}) \to
B
.
\end{align*}
The $T^{\circ}$-action on
$\overline{\val}^{-1} ( \overline{\Delta})$
is a $T^{\circ}$-action over $B$.

There are open subsets 
${\val}^{-1} ( \relin ({\Delta})) \subset
{\val}^{-1} ( {\Delta} )$ 
and
$\overline{\val}^{-1} ( \relin (\overline{\Delta})) \subset
\overline{\val}^{-1} ( \overline{\Delta} )$ 
with $T^{\circ}$-actions,
where $\relin ({\Delta})$ and
$\relin (\overline{\Delta})$ are
the relative interiors of $\Delta$ and
$\overline{\Delta}$, respectively.
We have morphisms
\begin{align} \label{tsuika}
{q^{\an}_{\Delta}} 
|_{{\val}^{-1} ( \relin ({\Delta}))}
:
{\val}^{-1} ( \relin ({\Delta})) \to
B
,
\quad
\overline{q^{\an}_{\Delta}} 
|_{\overline{\val}^{-1} ( \relin (\overline{\Delta}))}
:
\overline{\val}^{-1} ( \relin (\overline{\Delta})) \to
B
\end{align}
with $T^{\circ}$-actions over $B$.

The \emph{Mumford model} $p  : \EEE \to \AAA$ associated to
$\CCC$ is constructed in \cite{gubler4}.
It is a formal model of $p^{\an} : E \to A^{\an}$,
and
there exists 
a unique morphism $q : \EEE \to \BBB$ extending the morphism
$q^{\an} : E \to B$ in the Raynaud extension.
Although we do not repeat the precise definition of it here,
we recall some properties which will be used later.
See \S~4 (especially 4.7) of \cite{gubler4} for more details.

For a formal affine open subset $\mathscr{V}$ of $\BBB$,
its generic fiber $V = \mathscr{V}^{\an}$ is an analytic subdomain of $B$.
Such an analytic subdomain is called a 
\emph{formal affinoid subdomain} of $B$.
Using this notion,
we can write
${\val}^{-1} ( \Delta) = \bigcup_{V} 
\left(
{\val}^{-1} ( \Delta) \cap (q^{\an})^{-1} (V)
\right)$,
where $V$ runs through the formal affinoid subdomains of $B$.
Further, we have
\[
{\val}^{-1} ( \Delta) \cap (q^{\an})^{-1} (V)
\cong U_{\Delta} \times V
,
\]
where
$U_{\Delta}$ is the rational subdomain of $(\GG_{m}^{n})^{\an}$
as in \S~\ref{tori}.
We set $U_{\Delta, V} := U_{\Delta} \times V$.
Then
$\EEE$ has 
a suitable formal affine open covering by the sets $\mathscr{U}_{V , \Delta}$,
where
$\mathscr{U}_{V , \Delta}$ 
satisfies
$(\mathscr{U}_{V , \Delta})^{\an} = U_{V , \Delta}$.
We set $\EEE_{\Delta} := \bigcup_{V} \mathscr{U}_{V , \Delta}$,
where $V$ runs through the formal affinoid subdomains of $B$.
We have a natural morphism $q |_{\EEE_{\Delta}} : \EEE_{\Delta} \to \BBB$
by restriction.
The restriction of
$p : \EEE \to \AAA$ to $\EEE_{\Delta}$
is an isomorphism onto its image $\AAA_{\overline{\Delta}} :=
p (\EEE_{\Delta})$.
Using the isomorphism
$p |_{\EEE_{\Delta}}$,
we define
\begin{align*}
\overline{q_{\Delta}} := q|_{\EEE_{\Delta}} \circ (p |_{\EEE_{\Delta}})^{-1}
:
\AAA_{\overline{\Delta}} \to \BBB
.
\end{align*}
We notice that $\overline{q_{\Delta}}$
depends not only on $\overline{\Delta}$
but also on the choice of a representative $\Delta$ of $\overline{\Delta}$.
over $\BBB$.

The $T^{\circ}$-action on $E$ over $B$ extends to the $\TTT^{\circ}$-action
on $\EEE$ over $\BBB$,
where $\TTT^{\circ}$ is the canonical model of $T^{\circ}$
(cf. \S~\ref{RtM}).
The $\TTT^{\circ}$-action on $\EEE$
restricts to a 
$\TTT^{\circ}$-action on
$\EEE_{\Delta}$ over $\BBB$.
Via the isomorphism
$
p_{|_{\EEE_{\Delta}}} : \EEE_{\Delta} \to \AAA_{\overline{\Delta}}
$,
we have a
$\TTT^{\circ}$-action on $\AAA_{\overline{\Delta}}$
over $\BBB$.
Note that $\overline{\val}^{-1} ( \overline{\Delta} ) =
(\AAA_{\overline{\Delta}})^{\an}$
as an analytic subspace of $A^{\an}$
and 
that 
this $\TTT^{\circ}$-action on $\AAA_{\overline{\Delta}}$
induces
the canonical 
$T^{\circ}$-action on $\overline{\val}^{-1} ( \overline{\Delta} )$
on the Raynaud generic fiber.

We recall that 
\cite[Proposition 4.8]{gubler4}
gives us a bijective correspondence between
the strata of $\tilde{\EEE}$ and the set of relative interiors
of the polytopes in $\CCC$,
and a bijective correspondence between
the strata of $\tilde{\AAA}$ and the set of relative interiors
of the polytopes in $\overline{\CCC}$.
In fact,
if
we set
$Z_{\relin ( {\Delta})} : = 
\red_{\EEE} ({\val}^{-1}  ( \relin ({\Delta} )))$
and
$Z_{\relin ( \overline{\Delta})} : = 
\red_{\AAA} (\overline{\val}^{-1}  ( \relin (\overline{\Delta} )))$,
then
$Z_{\relin ( {\Delta})}$ and
$Z_{\relin ( \overline{\Delta})}$
are the strata of $\tilde{\EEE}$ and $\tilde{\AAA}$
corresponding to 
$\relin ( {\Delta})$
and
$\relin ( \overline{\Delta})$ respectively via 
these
bijective correspondences.

Taking the reductions of the morphisms in (\ref{tsuika}),
we obtain surjective morphisms
\addtocounter{Claim}{1}
\begin{align} \label{overlineqtilde}
\tilde{q_{\Delta}}
|_{Z_{\relin ({\Delta})}}
:
Z_{\relin ({\Delta})} \to \tilde{\BBB}
,
\quad
\tilde{
\overline{q_{\Delta}}}
|_{Z_{\relin (\overline{\Delta})}}
:
Z_{\relin (\overline{\Delta})} \to \tilde{\BBB}
.
\end{align}
The torus
$\tilde{\TTT^{\circ}}
$ over $\tilde{\KK}$
acts on 
$Z_{\relin ({\Delta})}$
and
$Z_{\relin (\overline{\Delta})}$,
and these actions are over $\tilde{\BBB}$
with respect to 
$\tilde{q_{\Delta}}$
and
$\tilde{
\overline{q_{\Delta}}}
|_{Z_{\relin (\overline{\Delta})}}$
respectively.
If ${\Delta}$ consists of a single point ${w}$,
then
we write 
$Z_{w}$
and
$Z_{\overline{w}}$, etc.
instead of 
$Z_{\{ {w} \}}$
and
$Z_{\{ \overline{w} \}}$, etc.
for simplicity.

\begin{Remark} \label{rem:torsor}
It follows from
\cite[Remark~4.9]{gubler4}
that
$\tilde{{q_{w}}} : Z_{{w}} \to \tilde{\BBB}$
and
$\tilde{\overline{q_{w}}} : Z_{\overline{w}} \to \tilde{\BBB}$
are $\tilde{\TTT^{\circ}}$-torsors.
\end{Remark}


We have seen that
the affinoid torus $T^{\circ}$ acts on 
$\overline{\val}^{-1} (\relin (\overline{\Delta}))$
and that the algebraic torus $\tilde{\TTT^{\circ}}$ acts on the reduction
$Z_{\relin (\overline{\Delta})}$.
These actions are compatible with respect to the
reduction map:
to be precise, we have the following lemma.

\begin{Lemma} \label{lem:new-compatibility:mumfordmodel}
With the notation above,
the diagram
\begin{align*}
\begin{CD}
T^{\circ} \times \overline{\val}^{-1} ( \relin (\overline{\Delta}))
@>>>
\overline{\val}^{-1} ( \relin ( \overline{\Delta}))
\\
@V{\red_{\TTT^{\circ} \times \AAA}}VV @VV{\red_{\AAA}}V
\\
\tilde{\TTT^{\circ}} \times Z_{\relin (\overline{\Delta})}
@>>>
Z_{\relin (\overline{\Delta})}
,
\end{CD}
\end{align*}
where the first row is the $T^{\circ}$-action
and the second row is the $\tilde{\TTT^{\circ}}$-action,
is commutative.
\end{Lemma}

\begin{Pf}
Let
$
\mu : \TTT^{\circ} \times \AAA_{\overline{\Delta}} \to \AAA_{\overline{\Delta}}
$
denote the $\TTT^{\circ}$-action 
recalled above.
Then the $T^{\circ}$-action on $\overline{\val}^{-1} ( \relin (\overline{\Delta}))$
is induced
from $\mu$ by taking the Raynaud generic fiber,
and the $\tilde{\TTT^{\circ}}$-action on $Z_{\relin (\overline{\Delta})}$
is induced from $\mu$ by taking the special fiber.
It follows from the definition of reduction maps
that the diagram is commutative.
\end{Pf}

Let $A_1$ and $A_2$ be abelian varieties over $\KK$
and let $\phi : A_1 \to A_2$ be a homomorphism.
Let
$n_i$ be the torus rank of $A_i$
and let $\overline{\val_i} : A_i^{\an} \to \RR^{n_i} / \Lambda_i$ be
the valuation map for $i = 1, 2$
and let
$\overline{\phi}_{\aff} : 
\RR^{n_1} / \Lambda_1 \to \RR^{n_2} / \Lambda_2$ be the
homomorphism as in (\ref{tamanidetekuru}).
Let
$\overline{\mathscr{C}_i}$
be a
$\Gamma$-rational polytopal decomposition of $\RR^{n_i} / \Lambda_i$
and let $\mathscr{A}_i$ 
be the Mumford model
associated to $\overline{\mathscr{C}_i}$.
Then
$\phi$ extends to a morphism $\mathscr{A}_1 \to
\mathscr{A}_2$ if, for any $\overline{\Delta_1} \in
\overline{\mathscr{C}_1}$,
there exists a polytope $\overline{\Delta_2} \in \overline{\mathscr{C}_2}$
such that $\overline{\phi}_{\aff} ( \overline{\Delta_1} )
\subset \overline{\Delta_2}$.
In fact, for any 
$\overline{\Delta_1} \in
\overline{\mathscr{C}_1}$ and 
$\overline{\Delta_2} \in \overline{\mathscr{C}_2}$
with
 $\overline{\phi}_{\aff} ( \overline{\Delta_1} )
\subset \overline{\Delta_2}$,
we can construct a unique
morphism $\left( \mathscr{A}_1\right)_{\overline{\Delta_1}}
\to \left( \mathscr{A}_2 \right)_{\overline{\Delta_2}}$
which extends the morphism $\phi |_{\overline{\val_1}^{-1} (\overline{\Delta_1})}
:
\overline{\val_1}^{-1} ( \overline{\Delta_1}) \to 
\overline{\val_2}^{-1} ( \overline{\Delta_2} )$,
and these morphisms patch together to be a morphism $\mathscr{A}_1 \to
\mathscr{A}_2$ extending $\phi$.


Finally in this subsection,
we define the notion of initial degenerations
of a closed subvariety $X$ of $A$.
For any $\Gamma$-rational polytope $\overline{\Delta}$ of $\RR^{n} / \Lambda$,
let $\XXX_{\overline{\Delta}}$ be the closure of 
$X^{\an} \cap \overline{\val}^{-1} ( \overline{\Delta} )$
in $\AAA_{\overline{\Delta}}$.
We define the \emph{initial degeneration 
$\initial_{\overline{\Delta}} (X)$
of $X$ over $\overline{\Delta}$}
by
\addtocounter{Claim}{1}
\begin{align} \label{initialdeg}
\initial_{\overline{\Delta}} (X) := \tilde{\XXX}_{\overline{\Delta}}
,
\end{align}
the special fiber of ${\XXX}_{\Delta}$.
It is determined from $X$ and $\overline{\Delta}$
and is independent of the polytopal decomposition $\mathscr{C}$
that was used before.
If $\overline{\Delta} = \{ \overline{w} \}$,
we write $\initial_{\overline{w}} (X)$
for
$\initial_{\overline{\Delta}} (X)$.
Note that $\initial_{\overline{w}} (X)$ is a closed subscheme of
$Z_{\overline{w}}$.
Moreover, since
$\mathscr{X}_{\overline{\Delta}} = 
\mathscr{X} \cap \mathscr{A}_{\overline{\Delta}}$,
we have
$\initial_{\overline{\Delta}} (X) =
\tilde{\mathscr{X}} \cap \tilde{\mathscr{A}_{\overline{\Delta}}}$
and $\initial_{\overline{w}} (X) = \tilde{\mathscr{X}} \cap Z_{\overline{w}}$.

\begin{Remark} \label{initial-dim}
Set $d: = \dim X$.
Let $W$ be an irreducible component of $\initial_{\overline{\Delta}} (X)$.
Since a formal model of $X$ is flat over $\KK^{\circ}$,
it follows from the definition of $\initial_{\overline{\Delta}} (X)$
that $\dim W = d$.
\end{Remark}

\begin{Remark} \label{rem:initialdegeneration}
Let $W$ be an irreducible component of $\initial_{\overline{\Delta}} (X)$.
By the definition of $\initial_{\overline{\Delta}} (X)$,
we have the reduction map
$X^{\an} \cap \overline{\val}^{-1} ( \overline{\Delta} ) 
\to \initial_{\overline{\Delta}} (X)$,
which is surjective.
Thus
there is a point $\xi_{W} \in X^{\an} \cap \overline{\val}^{-1} ( \overline{\Delta} ) $ which maps to
the generic point of $W$ by this reduction map.
Then we say that
$\xi_{W}$
\emph{reduces to the generic point of $W$}.
\end{Remark}


\section{Torus-equivalence between strata and canonical simplices} \label{non-degstrata2}

We begin by fixing the notation throughout this section.
Let $A$ be an abelian variety over $\KK$
with Raynaud extension
\[
\begin{CD}
1 @>>>
(\GG_{m}^{n})^{\an}
@>>>
E
@>{q^{\an}}>>
B
@>>>
0
.
\end{CD}
\]
Let $p^{\an} : E \to A^{\an}$ be the uniformization.
Let $\val : E \to \RR^{n}$ 
and
$\overline{\val} : A^{\an} \to \RR^{n} / \Lambda$
be the valuation maps 
with respect to the standard coordinates $x_{1} , \ldots , x_{n}$
(cf. (\ref{valuationmap})),
where $\Lambda = \val (\Ker p^{\an})$.
Fix a $\Lambda$-periodic
$\Gamma$-rational polytopal decomposition 
$\CCC_{0}$ of $\RR^{n}$,
and let
$\overline{\CCC_{0}}$
denote the polytopal decomposition
of $\RR^{n} / \Lambda$
induced from $\CCC_{0}$ by quotient.
Let $p_{0} : \EEE_{0} \to 
\AAA_{0}$ be the Mumford model associated to $\CCC_{0}$
(cf. \S~\ref{RtM2}).
Let 
$\XXX'$ be a 
connected 
strictly semistable formal scheme over $\KK^{\circ}$
(cf. \S~\ref{ssfss})
and let
$\varphi_{0} : \XXX' \to \AAA_{0}$ be a 
generically finite
morphism.
We put $X' := (\XXX')^{\an}$ and $d := \dim X'$.
Let
$f : X' \to A^{\an}$ be the restriction
of $\varphi_{0}$
to the Raynaud generic fiber.

In this section, we show
that $f$ induces a torus-equivariant morphism
over the formal fiber over a point in a stratum of $\tilde{\mathscr{X}'}$
(cf. Lemma~\ref{restrictionlemma}).
Using this fact together with Lemma~\ref{torus-action-ff}, we show that the 
morphism between some models of $X'$ and $A^{\an}$ induces 
a torus-equivariant morphism 
between their strata (cf. Proposition~\ref{equivariantdiagram}).
Then, we give 
in Lemma~\ref{mainproposition}
a sufficient condition for a canonical simplex
$\Delta_S$ of $S ( \mathscr{X}' )$ not to be collapsed by $f$.
This lemma will play a key role in the sequel.

\subsection{Associated affine map} \label{associatedaffine}

Recall that $\Val : X' \to S(\XXX')$ is
the retraction map to the skeleton
(cf. \S~\ref{ssfss}).
The following assertion is due to Gubler:

\begin{citeProposition} [Proposition~5.11 of \cite{gubler4}]
\label{Val}
Under the setting above,
there is a unique map $\overline{f}_{\aff}: S ( \XXX' )
\to \RR^{n} / \Lambda$ with $\overline{f}_{\aff} 
\circ \Val = \overline{\val} \circ f$
on $X'$.
The map $\overline{f}_{\aff}$
is continuous.
For any $S \in \str ( \tilde{\XXX'})$,
the restriction of $\overline{f}_{\aff}$
to the canonical simplex $\Delta_{S}$
(cf. \S~\ref{ssfss})
is an affine map and there exists a unique $\overline{\Delta}
\in \CCC_{0}$
with $\overline{f}_{\aff} ( 
\relin
(\Delta_{S}) ) \subset \relin (\overline{\Delta})$.
\end{citeProposition}

We
recall how $\overline{f}_{\aff}$
is described over $\Delta_{S}$.
Set 
$r := d - \dim S$,
the codimension of $S$ in $\tilde{\XXX'}$.
%
%
We take an affine open subset $\UUU' \subset \XXX'$
and an \'etale morphism $\psi : \UUU' \to \SSS = \SSS_{1} \times \SSS_{2}$,
where
$\SSS_{1} = \KK^{\circ} 
\langle
x_{0}' , \ldots , x_{r}'
\rangle
/
(x_{0}'  \cdots  x_{r}' - \pi)$
with $\pi \in \KK^{\circ \circ} \setminus \{ 0 \}$,
such that $S$ is a distinguished stratum associated to $\UUU'$
as in
Proposition~\ref{gubler4Proposition5.2}.
Since $p^{\an} : E \to A^{\an}$ is a local isomorphism,
replacing $\UUU'$ by a non-empty open subset of it,
we can take
a local lift
$F : (\UUU')^{\an} \to q^{-1} (V) 
$ of $f : X' \to A^{\an}$,
where $V$ is a formal affinoid subdomain of $B$.
We fix an identification $q^{-1} (V) \cong (\GG_{m}^{n})^{\an} \times V$.
By \cite[Proposition 2.11]{gubler1},
we have an expression
\addtocounter{Claim}{1}
\begin{align} \label{localF}
F^{\ast} (x_{i}) = \lambda_{i} v_{i} \psi^{\ast} (x_{1}')^{m_{i1}}
\cdots \psi^{\ast} (x_{r}')^{m_{ir}}
\end{align}
with some $\lambda_{i} \in \KK^{\times}$, $v_{i} \in \OO (\UUU')^{\times}$
and $\mathbf{m}_{i} = (m_{i1} ,\ldots , m_{ir}) \in \ZZ^{r}$.
We define $f_{\aff} : S ( \UUU' ) \to \RR^{n}$,
via the identification 
\[
S ( \UUU' )=
\Delta_{S} =
\{
(u'_{1} , \ldots , u'_{r})
\in \RR_{\geq 0}^{r}
\mid
u'_{1} + \cdots + u'_{r} \leq v (\pi)
\}
\] 
in (\ref{canonicalsimplexidentification}),
by
\addtocounter{Claim}{1} 
\begin{align} \label{faff}
f_{\aff} (\mathbf{u}') =
( \mathbf{m}_{1} \cdot \mathbf{u}' + v (\lambda_{1})
, \ldots,
\mathbf{m}_{n} \cdot \mathbf{u}' + v (\lambda_{n}))
,
\quad
\mathbf{u}' \in \Delta_{S} = S ( \UUU')
.
\end{align}
Then
$\overline{f}_{\aff}$
is 
the composite of 
$f_{\aff} : S ( \UUU' ) \to \RR^{n}$
with the quotient $\RR^{n} \to \RR^{n} / \Lambda$.

Let $\rank {\overline{f}_{\aff}}$ denote the 
rank of the linear part of $f_{\aff}$, i.e.,
the rank of the matrix $(m_{ij})$.
Note that $\rank {\overline{f}_{\aff}}$
does not depend on
the choice of a lift $F$
or the identification $q^{-1} (V) \cong (\GG_{m}^{n})^{\an} \times V$,
and thus
well-defined from $f$ and $\Delta_{S}$.

\subsection{Torus-equivariance between strata} \label{tebos}

Let 
$S$ be a stratum of $\tilde{\XXX'}$
of codimension $r$
and let $\Delta_S$ be the canonical simplex corresponding to $S$.
Let $m_{ij}$, for $1 \leq i \leq n$ and $1 \leq j \leq r$,
be integers which give the linear part of
$f_{\aff}$ (cf. (\ref{faff})).
We define a homomorphism
$h_{f,\Delta_{S}} : (\GG_{m}^{r})^{\an}_{1} \to (\GG_{m}^{n})^{\an}_{1}$
by
\[
h_{f,\Delta_{S}}^{\ast} (x_{i}) = (x_{1}')^{{m}_{i1}} \cdots (x_{r}')^{{m}_{ir}},
\quad
i = 1 , \ldots , n.
\]

Since
$(\GG_m^{n})^{\an}$ is the torus part of the Raynaud extension of $A$,
the affinoid torus $(\GG_{m}^{n})^{\an}_{1}$ 
acts on $A^{\an}$.
Then
we make $(\GG_{m}^{r})^{\an}_{1}$
act on $A^{\an}$
via $h_{f,\Delta_{S}}$.

\begin{Remark} \label{rem:torus-paralell}
Let $\tilde{h}_{f,\Delta_{S}} : 
(\GG_{m}^{r})_{\tilde{\KK}} \to (\GG_{m}^{n})_{\tilde{\KK}}$
be the homomorphism given by the reduction of $h_{f,\Delta_{S}}$.
Then the image $
\tilde{h}_{f,\Delta_{S}}
\left(
(\GG_{m}^{r})_{\tilde{\KK}}
\right)
$
is a subtorus of $(\GG_{m}^{n})_{\tilde{\KK}}$
of
$\dim \overline{f}_{\aff} ( \Delta_{S})$,
which is the rank of (the linear part of) $f_{\aff}$.
\end{Remark}

\begin{Lemma} \label{restrictionlemma}
With the notation above,
let $\tilde{p} \in S 
$ be a closed point
and let $X'_{+} (\tilde{p}) := \red_{\XXX'}^{-1} (\tilde{p})$
be the formal fiber over $\tilde{p}$.
Then 
the morphism
\begin{align*} 
f|_{X'_{+} (\tilde{p})} : X'_{+} (\tilde{p}) 
\to
A^{\an}
\end{align*}
is 
$(\GG_{m}^{r})^{\an}_{1}$-equivariant
with respect to the
$(\GG_{m}^{r})^{\an}_{1}$-action on $X'_{+} (\tilde{p})$
given in \S~\ref{construction-action}
and that on $A^{\an}$
induced by $h_{f,\Delta_{S}}$.
\end{Lemma}

\begin{Pf}
We take a formal affine open subscheme $\UUU' \subset \XXX'$
such that $\tilde{p} \in \tilde{\UUU'}$. 
We then note
$X'_{+} (\tilde{p}) \subset (\UUU')^{\an}$.
Replacing $\UUU'$ by a
formal open subscheme containing $\tilde{p}$,
we have a local lift $F$
of $f$ as in \S~\ref{associatedaffine},
so that
we have a commutative diagram
\[
\begin{CD}
(\UUU')^{\an} @>{F}>> (\GG_{m}^{n})^{\an} \times V \\
@AA{\cup}A @VV{p^{\an}}V \\
X'_{+}( \tilde{p} ) @>{f}>> A^{\an}
,
\end{CD}
\]
where $V$ is a formal affinoid subdomain of $B$,
and
$(\GG_{m}^{n})^{\an} \times V$ is regarded as a subdomain of $E$.

Since the right column in the above diagram 
is $(\GG_{m}^{r})^{\an}_{1}$-equivariant,
it suffices to show that $F$ is a
$(\GG_{m}^{r})^{\an}_{1}$-equivariant morphism.
Recall that we have an isomorphism
$\psi'^{\an} |_{X'_{+}  (\tilde{p})}: X'_{+}  (\tilde{p}) \to (\SSS^{\an})_{+} (\tilde{\psi}(\tilde{p}))$
in (\ref{morphisms-of-formalfibers}).
Let $G$ be the composite
\[
\begin{CD}
(\SSS^{\an})_{+} ( \tilde{\psi} (\tilde{p})) @>{(\psi'^{\an}|_{X'_{+}  (\tilde{p})})^{-1}}>> 
X'_{+} (\tilde{p})
@>{F}>> (\GG_{m}^{n})^{\an} \times V
@>{pr_{1}}>> (\GG_{m}^{n})^{\an}
,
\end{CD}
\]
where 
$pr_{1}$ is the first projection.
We note that, if $G$ is $(\GG_{m}^{r})^{\an}_{1}$-equivariant,
then $F$ is also $(\GG_{m}^{r})^{\an}_{1}$-equivariant.
Indeed,
suppose that $G$ is $(\GG_{m}^{r})^{\an}_{1}$-equivariant.
Since the $(\GG_{m}^{r})^{\an}_{1}$-action on $X'_{+} ( \tilde{p} )$
in Lemma~\ref{torus-action-ff}
is defined through 
$\psi'^{\an}|_{X'_{+}  (\tilde{p})}$,
this morphism $\psi'^{\an}|_{X'_{+}  (\tilde{p})}$
is $(\GG_{m}^{r})^{\an}_{1}$-equivariant by definition.
Thus
$G \circ \psi'^{\an}|_{X'_{+}  (\tilde{p})}$
is
$(\GG_{m}^{r})^{\an}_{1}$-equivariant, and hence 
$pr_{1} \circ F$, which equals $G \circ \psi'^{\an}|_{X'_{+}  (\tilde{p})}$,
is $(\GG_{m}^{r})^{\an}_{1}$-equivariant.
Since the $(\GG_{m}^{r})^{\an}_{1}$-action on $(\GG_{m}^{n})^{\an} \times V$
is the pull-back of the $(\GG_{m}^{r})^{\an}_{1}$-action on $(\GG_{m}^{n})^{\an}$
by $pr_{1}$,
it follows that $F$ is $(\GG_{m}^{r})^{\an}_{1}$-equivariant.

Thus we only have to show that $G$ is $(\GG_{m}^{r})^{\an}_{1}$-equivariant.
By the description (\ref{localF}),
we have
$
G^{\ast} (x_{i}) = 
\lambda_{i} v_{i}' (x_{1}')^{m_{i1}}
\cdots (x_{r}')^{m_{ir}}
$,
where $\lambda_{i} \in \KK^{\times}$
and $v_{i}' \in \OO (\SSS)^{\times}$.
It follows
from this description
and the definition of $h_{f,\Delta_{S}}$
that
$G$ is $(\GG_{m}^{r})^{\an}_{1}$-equivariant,
which concludes the lemma.
\end{Pf}

Let $\CCC$ be a $\Lambda$-periodic 
$\Gamma$-rational polytopal subdivision
of ${\CCC_{0}}$,
with the induced polytopal subdivision 
$\overline{\CCC}$ of $\overline{\CCC_{0}}$,
and let
$p : \EEE \to \AAA$ be the Mumford model
of the uniformization $p^{\an} : E \to A^{\an}$ 
associated to $\CCC$.
Let $\BBB$ be the formal abelian scheme with $\BBB^{\an} = B$.
Recall that there is a morphism
$q : \EEE \to \BBB$ which restricts to 
the morphism
$q^{\an} : E \to A^{\an}$
on the Raynaud extension.
Also, recall that,
for any $\overline{\Delta} \in \overline{\CCC}$,
we obtain a natural $(\GG_{m}^{n})^{\an}_{1}$-action on 
$\overline{\val}^{-1} ( \relin (\overline{\Delta}))$
(cf. \S~\ref{RtM2}).
Via the homomorphism $h_{f,\Delta_{S}}$, we have
a $(\GG_{m}^{r})^{\an}_{1}$-action on 
$\overline{\val}^{-1} ( \relin (\overline{\Delta}))$.

Let
$\DDD$ be the subdivision of $S (\XXX')$ given by
\[
\DDD :=
\left\{ 
\left.\Delta_{S'} \cap \overline{f}_{\aff}^{-1} (\overline{\Delta}) 
\right|
S' \in \str (\tilde{\XXX'} ),
\overline{\Delta} \in \overline{\CCC}
\right\}
.
\]
Let
$\XXX''$
be the formal model of $X'$ associated
to $\DDD$
and let $\iota' : \XXX'' \to \XXX'$
be
the morphism extending the identity on $X'$
as in Remark~\ref{rem:Proposition5.5gubler4}.
Then we have a morphism
$\varphi' : \XXX'' \to \AAA$
extending $f : X' \to A^{\an}$
by \cite[Proposition 5.14]{gubler4}.

Let $\tilde{h}_{f,\Delta_{S}} : 
(\GG_{m}^{r})_{\tilde{\KK}} \to (\GG_{m}^{n})_{\tilde{\KK}}$
be the homomorphism given by the reduction of $h_{f,\Delta_{S}}$.
Since
$Z_{\relin \overline{\Delta}} = \red_{\AAA} 
\left( \overline{\val}^{-1} ( \relin (\overline{\Delta}))
\right)$
has a $(\GG_{m}^{n})_{\tilde{\KK}}$-action
(cf. \S~\ref{RtM2}),
we obtain
a $(\GG_{m}^{r})_{\tilde{\KK}}$-action on
$Z_{\relin \overline{\Delta}}$ via
$\tilde{h}_{f,\Delta_{S}}$.

\begin{Proposition} \label{equivariantdiagram}
Let $u \in \DDD$ be a vertex,
$\Delta_{S}$ the canonical simplex of $S ( \XXX' )$ with
$u \in \relin \Delta_{S}$
and let
$R$ be the stratum of $\tilde{\XXX''}$ 
corresponding to $u$
(cf. Proposition~\ref{gubler4Proposition5.7}).
Let $\overline{\Delta} \in 
\overline{\CCC}$ be 
the polytope with 
$\overline{f}_{\aff} (u) \in \relin ( \overline{\Delta})$.
Take a representative $\Delta \in \CCC$ of $\overline{\Delta}$
and let
$\tilde{\overline{q_{\Delta}}} |_{Z_{\relin \overline{\Delta}} }
: Z_{\relin \overline{\Delta}} 
\to \tilde{\BBB}$
be the morphism in (\ref{overlineqtilde}).
Then,
there exists a unique morphism $\beta_{\Delta} : S \to \tilde{\BBB}$
such that the diagram
\begin{align*}
\begin{CD}
R @>{\tilde{\varphi'}}>> Z_{\relin \overline{\Delta}} \\
@V{\tilde{\iota'}}VV 
@VV{\tilde{\overline{q_{\Delta}}} |_{Z_{\relin \overline{\Delta}} }}V \\
S @>{\beta_{\Delta}}>> \tilde{\BBB}
\end{CD}
\end{align*}
commutes.
Moreover, 
this diagram is $(\GG_{m}^{r})_{\tilde{\KK}}$-equivariant
with respect to the
$(\GG_{m}^{r})_{\tilde{\KK}}$-action on $R$ in
Proposition~\ref{gubler4Corollary5.9},
that on
$Z_{\relin \overline{\Delta}}$ induced by $\tilde{h}_{f , \Delta_{S}}$
and the trivial actions on $S$ and $\tilde{\BBB}$.
\end{Proposition}

\begin{Pf}
We first define the morphism $S \to \tilde{\BBB}$.
Since there is a local section of
$\tilde{\iota'}$,
we define locally on $S$ a morphism from $S$ to $\tilde{\BBB}$
which is compatible with 
$\tilde{\overline{q_{\Delta}}} |_{Z_{\relin \overline{\Delta}} } 
\circ \tilde{\varphi'}$.
Since the fiber of $\tilde{\iota'}$ is an algebraic torus
and $\tilde{\BBB}$ is an abelian variety,
the morphism
$\tilde{\overline{q_{\Delta}}} |_{Z_{\relin \overline{\Delta}} } 
\circ \tilde{\varphi'}$
contracts
any fiber of $\tilde{\iota'}$ to a point.
This implies that
the local morphism from $S$ to $\tilde{\BBB}$ defined in this way
does not depend on the choice of local sections of 
$\tilde{\iota'}$.
It follows that the local morphisms patch together
to be a global morphism $\beta_{\Delta} : S \to \tilde{\BBB}$,
which satisfies the commutativity of the diagram.

The uniqueness of $\beta_{\Delta}$
follows from the construction.
It only remains to show that 
$\tilde{\varphi'} : R \to Z_{\relin \overline{\Delta}}$
is $(\GG_{m}^{r})_{\tilde{\KK}}$-equivariant.
Take an arbitrary closed point $\tilde{p} \in S
$.
We have relations
$
R = \red_{\XXX''} ( \Val^{-1} ( u ))$
and $
\{ u \}
=
\Val ( (\red_{\XXX''})^{-1} ( \{ \tilde{p} \} \times_{S} R ))
$ by
Proposition~\ref{gubler4Proposition5.7}.
We also have
$\overline{\val} \circ f = \overline{f}_{\aff} \circ \Val$
by
Proposition~\ref{Val}.
It follows that
\[
\overline{\val}
\left(
f ( (\red_{\XXX''})^{-1} ( \{ \tilde{p} \} \times_{S} R ))
\right)
=
\overline{f}_{\aff}
\left(
\Val ( (\red_{\XXX''})^{-1} ( \{ \tilde{p} \} \times_{S} R ))
\right)
=
\{ \overline{f}_{\aff} ( u ) \}
\subset \relin \overline{\Delta}
,
\]
and thus
$f$ restricts to a morphism
\addtocounter{Claim}{1}
\begin{align} \label{proofofequiv}
(\red_{\XXX''})^{-1} ( \{ \tilde{p} \} \times_{S} R )
\to
\overline{\val}^{-1} (\relin (\overline{\Delta})
)
.
\end{align}
Recall that
we have a $(\GG_{m}^{r})^{\an}_{1}$-action
on
$(\red_{\XXX''})^{-1} ( \{ \tilde{p} \} \times_{S} R )
\subset X'_{+} ( \tilde{p} )$ 
(cf. Remark~\ref{nd4}).
Since
$\overline{\val}^{-1} (\relin (\overline{\Delta}))$
is stable under the
$(\GG_{m}^{n})^{\an}_{1}$-action (cf. \S~\ref{RtM2}),
it is stable under the $(\GG_{m}^{r})^{\an}_{1}$-action
induced by $h_{f, \Delta_{S}}$.
It follows from Lemma~\ref{restrictionlemma}
that
the morphism (\ref{proofofequiv}) is $(\GG_{m}^{r})^{\an}_{1}$-equivariant.

We deduce 
the 
$(\GG_{m}^{r})_{\tilde{\KK}}$-equivariance
of
$\tilde{\varphi'} : R \to Z_{\relin \overline{\Delta}}$
from the $(\GG_{m}^{r})^{\an}_{1}$-equivariance of
(\ref{proofofequiv}).
Indeed,
it follows from
Lemma~\ref{lem:new-compatibility:mumfordmodel}
that 
the diagram
\begin{align*}
\begin{CD}
(\GG_{m}^{r})^{\an}_{1} \times \overline{\val}^{-1} (\relin (\overline{\Delta}))
@>>>
\overline{\val}^{-1} ( \relin (\overline{\Delta}))
\\
@V{\red_{(\GG_{m}^{r})^{\mathrm{f-sch}}_{1} \times \AAA}}VV @VV{\red_{\AAA}}V
\\
(\GG_{m}^{r})_{\tilde{\KK}} \times Z_{\relin (\overline{\Delta})}
@>>>
Z_{\relin (\overline{\Delta})}
,
\end{CD}
\end{align*}
where the rows are the torus-actions given by ${h}_{f , \Delta_{S}}$
and $\tilde{h}_{f , \Delta_{S}}$,
is commutative, 
that is,
the $(\GG_{m}^{r})^{\an}_{1}$-action on
$\overline{\val}^{-1} (\relin (\overline{\Delta}))$
and the $(\GG_{m}^{r})_{\tilde{\KK}}$-action on $Z_{\relin (\overline{\Delta})}$
are compatible with respect to reduction.
On the other hand,
by
Lemma~\ref{torus-action-ff},
the
$(\GG_{m}^{r})^{\an}_{1}$-action on 
$(\red_{\XXX''})^{-1} ( \{ \tilde{p} \} \times_{S} R )$
and the $(\GG_{m}^{r})_{\tilde{\KK}}$-action on
$\{ \tilde{p} \} \times_{S} R$
are compatible with respect to the reduction map $\red_{\XXX''}$.
Since
(\ref{proofofequiv})
is $(\GG_{m}^{r})^{\an}_{1}$-equivariant,
it
follows that
$
\{ \tilde{p} \} \times_{S} R \to Z_{\relin \overline{\Delta}}
$
is $(\GG_{m}^{r})_{\tilde{\KK}}$-equivariant.
Since 
$\tilde{p} \in S$ is any closed point,
this implies
$R \to Z_{\relin \overline{\Delta}}$ is $(\GG_{m}^{r})_{\tilde{\KK}}$-equivariant.
Thus we complete the proof.
\end{Pf}

\subsection{Key lemma}

The following
is the key lemma in this paper.
It is  crucially used in
the proof of
Proposition~\ref{nondeg-supp2}.

\begin{Lemma} \label{mainproposition}
Let $X$ be a $d$-dimensional closed subvariety of $A$
and let $\XXX_{0}$ be the closure of $X$ in $\AAA_{0}$.
Assume that the morphism $\varphi_{0} : \XXX' \to \AAA_{0}$
factors through $\XXX_{0}$ to be
a generically finite 
surjective morphism $\XXX' \to \XXX_{0}$.~\footnote{
If this morphism $\XXX' \to \XXX_{0}$ proper moreover, it is called a
semistable alteration (cf. \S~\ref{sa}).}
Let $\Delta_{S}$ be a canonical simplex of $S ( \XXX' )$.
Suppose that
there exist a
$\Gamma$-rational point
$\overline{w} \in \overline{\val} (X^{\an})$,
 an irreducible component $W$ of 
$\initial_{\overline{w}} (X)$,\footnote{Recall that 
$\initial_{\overline{w}} (X)$ denotes the initial degeneration of $X$
at
$\overline{w}$ (cf. (\ref{initialdeg})).}
and a 
$\Gamma$-rational
point $u \in \relin ( \Delta_{S} )$
such that
$f (u)$ reduces to the generic point of $W$
(cf. Remark~\ref{rem:initialdegeneration}).
Then we have
$\dim \overline{f}_{\aff} ( \Delta_{S}) = \dim \Delta_{S}$.
\end{Lemma}

\begin{Pf}
Let $\overline{\CCC}$ be a
$\Gamma$-rational subdivision of $\overline{\CCC_{0}}$ 
such that $\overline{w}$ itself is a vertex of $\overline{\CCC}$.
Let $\XXX''$ be the formal model of $X'$
associated to
the polytopal subdivision
\[
\DDD :=
\left\{
\left.
\Delta_{S'} \cap \overline{f}_{\aff}^{-1} (\overline{\Delta})
\right|
S' \in \str( \tilde{\XXX'})
,
\overline{\Delta} \in \overline{\CCC}
\right\}
\]
of 
$S (\XXX')$,
and let
$\iota' : \XXX'' \to \XXX'$ denote the natural morphism
extending the identity on the generic fiber.
Then $u$ is a vertex of $\DDD$.
Let $R \in \str (\tilde{\XXX''})$ be the stratum corresponding to $u$.
Then
$(\GG_{m}^{r})_{\tilde{\KK}}$ acts on $R$,
which makes
$\tilde{\iota'} : R \to S$
a torus bundle of
relative dimension 
$r = \dim \Delta_{S}$
(cf. Proposition~\ref{gubler4Corollary5.9}).

Let $\EEE$ and $\AAA$ be the Mumford models of $E$ and $A$ 
associated with $\overline{\CCC}$.
Let $q : \EEE \to \BBB$
be the surjective morphism 
obtained from the Raynaud extension
(cf. \S~\ref{RtM2}).
By \cite[Proposition~5.14]{gubler4},
there is
a unique morphism $\varphi' : \XXX'' \to \AAA$
extending $f$.
Since $f (u)$ reduces to the generic point of $W$,
we have
\[
\overline{f}_{\aff} (u) = \overline{\val}( f (u))
= \overline{w}
\in
\relin \{ \overline{w} \}
.
\]
Thus
we 
have a 
morphism
$\tilde{\varphi'} : R \to Z_{\overline{w}}$,
which is $\left( \GG_{m}^{r} \right)_{\tilde{\KK}}$-equivariant 
by Proposition~\ref{equivariantdiagram}.

Since $f (u)$ reduces to the generic point of $W$,
Lemma~\ref{vertex-genericpoint}
tells us
that 
the generic point of $R$ maps to that of $W$.
Since
$W$ is a closed subset of
$Z_{\overline{w}}$,
it follows that
$W$ coincides with the closure of
$\tilde{\varphi'} (R)$ in $Z_{\overline{w}}$.
Since the morphism
$\tilde{\varphi'} : R \to Z_{\overline{w}}$
is $(\GG_{m}^{r})_{\tilde{\KK}}$-equivariant,
we see that $W$ is stable under the 
$(\GG_{m}^{r})_{\tilde{\KK}}$-action
induced by $\tilde{h}_{f,\Delta_{S}}$.
It follows 
that
$W$ is stable under the action of $\TT''
:= \Image  \tilde{h}_{f,\Delta_{S}}$,
the image of $\tilde{h}_{f,\Delta_{S}}$.
Thus, we obtain a $\TT''$-action on $W$.
Since the action of $\TT''$ on $Z_{\overline{w}}$ is free 
by Remark~\ref{rem:torsor},
we see that the $\TT''$-action on $W$ is free.

We set $\Xi := W/ \TT''$.
Since $\TT''$-action on $W$ is free,
%
we have
$\dim \Xi = d - \dim \overline{f}_{\aff} ( \Delta_{S} )$
(cf. Remarks~\ref{initial-dim}
and \ref{rem:torus-paralell}).
Thus we have
\addtocounter{Claim}{1}
\begin{align} \label{S-Xi}
\dim S = d - \dim \Delta_{S} \leq d - \dim \overline{f}_{\aff} ( \Delta_{S} )
= \dim \Xi
\end{align}
on one hand.

Since 
$\Xi$ is the quotient of $W$ by the action of
$(\GG_{m}^{r})_{\tilde{\KK}}$ induced by
$\tilde{h}_{\Delta_{S} , f} : (\GG_{m}^{r})_{\tilde{\KK}} \to \TT''$
and since
$R / (\GG_{m}^{r})_{\tilde{\KK}} = S$,
the composite $R \to W \to \Xi$ factors through $R \to S$.
Since 
$\tilde{\varphi'} : R \to W$
is dominant and since $W \to \Xi$ is surjective,
the morphism
$S \to \Xi$ is dominant.
Thus we have
$\dim S \geq \dim \Xi$,
on the other hand.
This inequality together with (\ref{S-Xi})
shows
$\dim S = \dim \Xi$. 
Further,
this equality together with (\ref{S-Xi}) gives us
$\dim \Delta_{S} = 
\dim \overline{f}_{\aff} ( \Delta_{S} )$.
This completes the proof of the key lemma.
\end{Pf}

\begin{Remark}
In the proof above,
we have shown that $
\Image \tilde{h}_{f,\Delta_{S}}$ 
acts on $W$ freely.
\end{Remark}

\subsection{Non-degenerate canonical simplices}

In this subsection,
we recall the notion of 
non-degeneracy
with respect to $f$ for canonical simplices
introduced
by Gubler
in
\cite[6.3]{gubler4}.
We show some properties of it also.

We begin by recalling what the non-degenerate 
canonical simplices are.
The morphism $\varphi_{0} : \XXX' \to \AAA_{0}$
gives us a morphism $\tilde{\varphi_{0}} : S \to \tilde{\AAA_{0}}$.
Then \cite[Lemma 5.15]{gubler4} gives us
a morphism $\tilde{\Phi}_{0} : S \to \tilde{\EEE_{0}}$
with $\tilde{p_{0}} \circ \tilde{\Phi}_{0} = \tilde{\varphi_{0}}$,
called a lift
of $\tilde{\varphi_{0}}$.
Let $q_{0} : \EEE_{0} \to \BBB$ be the extension of $q^{\an} : E \to A^{\an}$
between the models.
A canonical simplex $\Delta_{S}$ is said to be
\emph{non-degenerate with respect to $f$}
if $\dim \overline{f}_{\aff} ( \Delta_{S} ) = \dim \Delta_{S}$
and $\dim \tilde{q_{0}} (\tilde{\Phi}_{0}  ( S)) = \dim S$.
This notion does not depend on the choice of $\tilde{\Phi}_{0}$.

Let $\CCC$,
$\AAA$,
$\DDD$,
and $\iota' : \XXX'' \to \XXX'$
be as in \S~\ref{tebos}.
Recall that $\varphi_{0}$ lifts to $\varphi' : \XXX'' \to \AAA$
by \cite[Proposition~5.14]{gubler4}.

\begin{Lemma} \label{lem:non-degenerate}
With the notation above,
let $u$ be a vertex of $\DDD$
and let $\Delta_{S}$ be the canonical simplex of $S (\XXX')$
with $u \in \relin \Delta_{S}$.
Let
$R$ 
be the stratum of $\tilde{\XXX''}$ corresponding to $u$.
Assume that $\overline{w} := \overline{f}_{\aff} (u)$
is a vertex of $\overline{\CCC}$.
Let $w \in \RR^{n}$ be a representative of $\overline{w}$
and let
$\beta_{w} : S \to \tilde{\mathscr{B}}$ be the morphism in Proposition~\ref{equivariantdiagram}
for $\{ w \} \in \CCC$.
Then,
$\Delta_{S}$ is non-degenerate with respect to $f$
if and only if
$\dim \overline{f}_{\aff} ( \Delta_{S} ) = \dim \Delta_{S}$
and 
$\dim \beta_{w} (S) = \dim S$.
\end{Lemma}

\begin{Pf}
It suffices to show $\beta_{w} = \tilde{q_{0}} \circ \tilde{\Phi}_{0}$
for some lift $\tilde{\Phi}_{0}$ of $\tilde{\varphi_{0}}$.
Since $\CCC$ is a $\Gamma$-rational subdivision of $\CCC_{0}$,
we have
an extension
$\overline{\iota_{0}} : \AAA \to \AAA_{0}$
of the
identity morphism of $A$ between the Mumford models
(cf. \cite[\S~4]{gubler4}).
Let 
$\Delta_{0} \in \CCC_{0}$ be the polytope with
$w \in  \relin \Delta_{0}$.
Then $\tilde{\overline{\iota_{0}}}$ restricts to $Z_{\overline{w}}
\to Z_{\relin \overline{\Delta_{0}}}$,
and there is a commutative diagram
\[
\begin{CD}
R @>{\tilde{\varphi'}}>> Z_{\overline{w}} 
@<{\tilde{p}|_{Z_{w}}}<{\cong}< Z_{w} 
@>{\tilde{q}
|_{Z_{w}}
}>> \tilde{\BBB} 
\\
@VV{\tilde{\iota'}}V @VV{\tilde{\overline{\iota_{0}}}}V @VV{\tilde{\iota_{0}}}V @VV{id}V \\
S @>{\tilde{\varphi_{0}}}>> Z_{\relin \overline{\Delta_{0}}}
@<{\tilde{p_{0}}|_{Z_{\relin \Delta_{0}}}}<{\cong}<
Z_{\relin \Delta_{0}} @>{\tilde{q_{0}}|_{Z_{\relin \Delta_{0}}}
}>> \tilde{\BBB}
.
\end{CD}
\]

We take the lift $\tilde{\Phi}_{0}$ of $\tilde{\varphi_{0}}$ such that
$\tilde{\Phi}_{0} (S) \subset Z_{\relin \Delta_{0}}$.
Then
\addtocounter{Claim}{1}
\begin{align} \label{waketa}
\tilde{q_{0}} \circ \tilde{\Phi}_{0} \circ \tilde{\iota'} 
=
\tilde{q_{0}} \circ (\tilde{p_{0}}|_{Z_{\relin \Delta_{0}}})^{-1} \circ
\tilde{\varphi_{0}} \circ \tilde{\iota'} 
= 
\tilde{q} \circ (\tilde{p}|_{Z_{w}})^{-1} \circ \tilde{\varphi'}
.
\end{align}
Since $\tilde{q} \circ (\tilde{p}|_{Z_{w}})^{-1} =
\tilde{\overline{q_{w}}}$
by definition
and
since
$
\tilde{\overline{q_{w}}} \circ \tilde{\varphi'}
= \beta_{w} \circ \tilde{\iota'}
$ by Proposition~\ref{equivariantdiagram},
it follows from
(\ref{waketa}) 
that
$\tilde{q_{0}} \circ \tilde{\Phi}_{0} \circ \tilde{\iota'} 
=
\beta_{w} \circ \tilde{\iota'}$.
Since $\tilde{\iota'}$ is surjective,
this concludes
$\tilde{q_{0}} \circ \tilde{\Phi}_{0}
= \beta_{w}$,
as required.
\end{Pf}

\begin{Lemma} \label{Snotuburekata}
Let $u$, $\Delta_{S}$ and $R$ be as in 
Lemma~\ref{lem:non-degenerate}.
Assume
that $\overline{w} := \overline{f}_{\aff} (u)$
is a vertex of $\overline{\CCC}$.
If $\Delta_{S}$ is non-degenerate with respect to $f$,
then $\dim \tilde{\varphi'} (R) = d$.
\end{Lemma}

\begin{Pf}
Recall that $(\GG_{m}^{n})_{\tilde{\KK}}$ acts on $Z_{\overline{w}}$
and that the  $(\GG_{m}^{r})_{\tilde{\KK}}$-action on $Z_{\overline{w}}$
is given by the homomorphism $\tilde{h}_{f,\Delta_{S}}$.
We set $\TT'' 
:= \Image \tilde{h}_{f,\Delta_{S}}$.
Then
there is a natural
$\TT''$-action on $Z_{\overline{w}}$.
By Proposition~\ref{equivariantdiagram},
the morphism 
$\tilde{\varphi'} : R \to Z_{\overline{w}}$
is $(\GG_{m}^{r})_{\tilde{\KK}}$-equivariant.
It follows that 
$\tilde{\varphi'} (R) \subset Z_{\overline{w}}$
is stable under the $\TT''$-action.
Thus we have a $\TT''$-action
on
$\tilde{\varphi'} (R)$.
Since the $(\GG_{m}^{n})_{\tilde{\KK}}$-action on $Z_{\overline{w}}$ is free
(cf. Remark~\ref{rem:torsor}),
the $\TT''$-action on $\tilde{\varphi'} (R)$ is free.

We take a representative $w \in \CCC$ of $\overline{w}$.
The morphism $\tilde{\overline{q_{w}}} : Z_{\overline{w}} \to \tilde{\BBB}$
in (\ref{overlineqtilde})
is the quotient by $(\GG_{m}^{n})_{\tilde{\KK}}$,
so that
any $\TT''$-orbit in $\tilde{\varphi'} (R)$
contracts to a point by $\tilde{\overline{q_{w}}}$.
It follows that
$\dim \tilde{\varphi'} (R) 
\geq \dim \tilde{\overline{q_{w}}} ( \tilde{\varphi'} (R) ) 
+ 
\dim \TT''$.
Since $\dim \TT'' = \dim \overline{f}_{\aff} ( \Delta_{S})$
by Remark~\ref{rem:torus-paralell},
we thus obtain
$\dim \tilde{\varphi'} (R) 
\geq \dim \tilde{\overline{q_{w}}} ( \tilde{\varphi'} (R) ) 
+ 
\dim \overline{f}_{\aff} ( \Delta_{S})$.

Suppose that $\Delta_{S}$ is non-degenerate with respect to $f$.
Then
$\dim \overline{f}_{\aff} ( \Delta_{S}) = \dim \Delta_{S}$
and $\dim S = \dim \beta_{w} (S)$
by Lemma~\ref{lem:non-degenerate}.
Since $\beta_{w} (S) = \tilde{\overline{q_{w}}} ( \tilde{\varphi'} (R) )$
(cf. Proposition~\ref{equivariantdiagram}),
we obtain
$\dim S = \dim \tilde{\overline{q_{w}}} ( \tilde{\varphi'} (R) )$.
Thus we have
\[
d \geq \dim \tilde{\varphi'} (R) 
\geq \dim \tilde{\overline{q_{w}}} ( \tilde{\varphi'} (R) )
 + \dim \overline{f}_{\aff} ( \Delta_{S}) = 
\dim S + \dim \Delta_{S} 
= d
,
\]
which concludes
$\dim \tilde{\varphi'} (R) = d$.
\end{Pf}

\section{Strict supports of canonical measures} \label{strictsupport}

In this section,
let $\mathbf{K}$
be
a subfield of $\KK$
such that the absolute value of $\KK$
restricts to a discrete absolute value on $\bfK$,
and assume that $\bfK$ is complete.
The value group $\Gamma_{\bfK}$ of $\bfK$ 
may be assumed to be a discrete subgroup of $\QQ$.
Let $\bfK^{\circ}$ denote the ring of integers of $\bfK$,
which is a discrete valuation ring.
For a variety $X$ over $\bfK$ and a 
flat formal scheme $\mathcal{X}$ of
finite type over $\bfK^{\circ}$,
we let
$X^{\an}$ and $\mathcal{X}^{\an}$
stand for $\left( X \times_{\Spec \bfK} \Spec \KK \right)^{\an}$
and $\left( \mathcal{X} \times_{\Spf \bfK^{\circ}} \Spf \KK^{\circ} \right)^{\an}$
respectively.
They are analytic spaces over $\KK$.
We deal only with
analytic spaces which 
arise from varieties 
defined over $\bfK$ and formal schemes defined over $\bfK^{\circ}$
in this section
because it
is enough for our latter applications.

\subsection{Mumford models over  a discrete valuation ring}

Let $A$ be an abelian variety over $\bfK$.
We recall that, 
replacing $\bfK$ by a finite extension in $\KK$ if necessary, 
we have the Raynaud extension over $\bfK$.
Indeed,
by \cite[Theorem~1.1]{BL}, 
replacing $\bfK$ by a finite extension in $\KK$ if necessary, 
we have an exact sequence
\begin{align*} 
\begin{CD}
1 @>>> (\GG_{m}^{n})^{\mathrm{f-sch}}_{\bfK^{\circ}} @>>> 
{\mathcal{A}}^{\circ} @>>> {\mathcal{B}} @>>> 0
\end{CD}
\end{align*}
of formal group schemes over $\bfK^{\circ}$,
where $n$ is the torus rank of $A$,
${(\GG_{m}^{n})^{\mathrm{f-sch}}_{\bfK^{\circ}}}$
is a split formal torus over $\bfK^{\circ}$,
${\mathcal{A}}^{\circ}$
is the formal completion of a semiabelian scheme over $\bfK^{\circ}$,
and
${\mathcal{B}}$ is the formal completion
of an abelian scheme over $\bfK^{\circ}$.
Further,
the base-change to $\KK^{\circ}$
is nothing but the exact sequence (\ref{qcptformalR}) 
for $A \times_{\Spec \bfK} {\Spec \KK}$.
We can also construct, in the same way as we did in \S~\ref{RtM},
a short exact sequence
of analytic groups over $\bfK$ whose base-change to $\KK$
coincides with (\ref{Raynaudext}).

Let $p^{\an} : E \to A^{\an}$ be the uniformization
and consider the lattice $\Lambda = \val (\Ker p^{\an})$ in $\RR^{n}$.
We see from \cite[Theorem~1.2]{BL} and its proof that
$\Lambda = \val (\Ker p^{\an}) \subset \Gamma_{\bfK}^{n}$,
and hence we have $\Lambda \subset \QQ^{n}$.
Thus
there exists 
a $\Lambda$-periodic rational
polytopal decomposition of $\RR^{n}$
in this setting.

\begin{Lemma} \label{lem:Mumford-model-DVR}
Let $A$ be an abelian variety over $\bfK$ of torus rank $n$.
Let
$\CCC$ be a $\Lambda$-periodic rational
polytopal decomposition of $\RR^{n}$
and let
$\overline{\CCC}$ be the induced rational polytopal decomposition
of $\RR^{n} / \Lambda$.
Then,
replacing $\bfK$ with a finite extension in
$\KK$ if necessary,
we have a proper flat formal scheme $\mathcal{A}$ over $\bfK^{\circ}$
such that $\mathcal{A} \times_{\Spf \bfK^{\circ}} \Spf \KK^{\circ}$
is the Mumford model of $A$ associated to $\overline{\CCC}$.
\end{Lemma}

\begin{Pf}
Replacing $\bfK$ by a finite extension in $\KK$,
we may assume that any vertex in $\CCC$
is in $\Gamma_{\bfK}^{n}$.
Then 
the construction of the Mumford models in \cite[4.7]{gubler4} 
works over $\bfK^{\circ}$ without any change,
so that
the Mumford model $\AAA$ associated to $\overline{\CCC}$
is defined over $\bfK^{\circ}$,
that is, there exists a formal model $\mathcal{A}$ of $A$
over $\bfK^{\circ}$ such that 
$\AAA = \mathcal{A} \times_{\Spf \bfK^{\circ}} \Spf \KK^{\circ}$.
\end{Pf}

\begin{Remark} \label{rem:Mm}
In the sequel,
when we say a Mumford model of $A$,
this means that it is the Mumford model
associated to a $\Lambda$-periodic rational polytopal
decomposition of $\RR^{n}$,
so that this Mumford model 
can be defined over a finite extension of $\bfK^{\circ}$
by Lemma~\ref{lem:Mumford-model-DVR}.
\end{Remark}

\subsection{Semistable alterations} \label{sa}
Let $\XXX_{0}$ be a connected admissible formal scheme over $\KK^{\circ}$.
A morphism $\XXX' \to \XXX_{0}$ 
of formal schemes  over $\KK^{\circ}$
is called a 
\emph{semistable alteration} for $\XXX_{0}$
if $\XXX'$ is a connected
strictly semistable formal scheme
and $\XXX' \to \XXX_{0}$ is a proper surjective generically finite morphism.
We say that a proper surjective generically finite morphism 
$\XX' \to \XX_{0}$
of 
flat
formal schemes
of finite type over $\bfK^{\circ}$ is a semistable alteration
if the base change of this morphism to $\KK^{\circ}$ is a semistable
alteration in the above sense.

In this subsection,
we give the existence of semistable alteration
for a model of a closed subvariety of an abelian variety
defined over $\bfK$.

\begin{Lemma} \label{lem:semistable-alteration}
Let $A$ be an abelian variety over $\bfK$ and let $X$
be a closed subvariety of $A$.
Let $\mathcal{A}_{0}$ be a proper flat formal scheme over $\bfK$
with $\mathcal{A}_{0}^{\an} = A^{\an}$,
and let ${\mathcal{X}_{0}}$ be the closure of $X$ in $\mathcal{A}_{0}$.
Then there exists a projective formal
scheme $\mathcal{X}'$ flat over $\Spf \bfK^{\circ}$
and a semistable alteration
$\varphi_{0} : {\mathcal{X}'} \to 
{\mathcal{X}_{0}}$.
\end{Lemma}

\begin{Pf}
Note that \cite[Proposition 10.5]{gubler0}
gives us
a projective scheme $\mathbf{X}_{1} \to \Spec \bfK^{\circ}$
with generic fiber $X$,
and a dominating morphism
$
\widehat{\mathbf{X}}_{1}
\to \XX_{0}
$
extending the identity on the generic fiber,
where $\widehat{\mathbf{X}}_{1}$ is the formal completion of $\mathbf{X}_{1}$
along its special fiber.
Indeed,
the statement of \cite[Proposition 10.5]{gubler0}
is given under the assumption that the base field $\bfK$
is algebraically closed, but the proof delivered there
works well also in our situation without any change.

Since $\widehat{\mathbf{X}}_{1}$ is algebraizable,
we apply
\cite[Theorem 6.5]{dejong} to $\widehat{\mathbf{X}}_{1}$
to obtain
a semistable alteration
$\XX' \to \widehat{\mathbf{X}}_{1}$ over $\bfK^{\circ}$
such that $\XX'$ is a projective formal scheme flat over $\Spf \bfK^{\circ}$.
Let $\varphi_{0} : \hat{\XX'} \to \hat{\XX_{0}}$ be the composite.
Then it is a semistable alteration, as required.
\end{Pf}

\begin{Remark} \label{assumption-sa}
Let $\AAA_{0}$ be such a Mumford model
as in Remark~\ref{rem:Mm}.
Let $X$ be a closed subvariety of $A$,
and let $\XXX_{0}$ the closure of $X$.
Since $\AAA_{0}$ is a formal scheme which
can be defined over a finite extension of $\bfK^{\circ}$,
$\XXX_{0}$ is a formal scheme which
can be defined over the finite extension as well.
By Lemma~\ref{lem:semistable-alteration},
we have a
semistable alteration $\varphi_{0} : \XXX' \to \XXX_{0}$ 
which can be defined over 
a finite extension of $\bfK^{\circ}$.
Further,
the restriction of
$\varphi_{0}$ to the Raynaud generic fibers is
regarded as a morphism
of projective varieties over a finite extension of $\bfK$.
In the sequel,
we consider such
semistable alterations only.
\end{Remark}

\subsection{Canonical measures and the canonical subset} \label{subsection:canonical}
Let $X$ be a proper variety over $\bfK$
and let $L$ be a line bundle on $X$.
In \cite{gubler4}, Gubler defined the notion of \emph{admissible metrics on $L$}
(cf. \cite[3.5]{gubler4}).
If
$\overline{L}$ is a line bundle on $X$ endowed with
an admissible metric,
then one can define a regular Borel measure
$\cherncl_{1} ( \overline{L})^{\wedge d}$
on $X^{\an}$
with suitable properties
(\cite[Proposition 3.8]{gubler4}).
It was originally introduced by Chambert-Loir
in \cite{chambert-loir}.
These measures satisfy the projection formula:
if $f : X' \to X$ is a morphism of $d$-dimensional geometrically integral
proper varieties over $\bfK$, then 
$f^{\ast} \overline{L}$ is an admissibly metrized line bundle
on $X'$, and
\begin{align*} 
f_{\ast} ( \cherncl_{1} ( f^{\ast} \overline{L})^{\wedge d} )
=
\deg (f)
\cherncl_{1} ( \overline{L})^{\wedge d}
.
\end{align*}

Let $A$ be an abelian variety over $\bfK$
and let $L$ be a line bundle on $A$.
We say that $L$ is \emph{even}
if
$[-1] L^{\ast} \cong L$, 
where, for any $m \in \ZZ$,  $[m] : A \to A$ denotes the homomorphism
given by $[m] (a) = ma$.
Suppose that
$L$ is even and ample.
As mentioned in \cite[Example 3.7]{gubler4},
there is
an important metric, called {canonical metric},
on $L$.
This metric is described as follows.
Let $m$ be a positive integer.
Since $L$ is even, we have an identification $[m]^{\ast} L = L^{m^{2}}$.
Then a metric $|| \cdot ||$ on $L$ is called a \emph{canonical metric}
if $[m]^{\ast} || \cdot || = || \cdot ||^{m^{2}}$ via that identification.
Note that a canonical metric depends on the choice of the
identification $[m]^{\ast} L = L^{m^{2}}$,
but it is unique up to positive rational multiples.

In the sequel, let $\overline{L}$ always denote a
line bundle endowed with a canonical metric
for a line bundle $L$ on an abelian variety.
For a closed subvariety $X$ of $A$ of dimension $d$,
the restriction $\overline{L}|_{X}$
is a line bundle on $X$ with an admissible metric
(cf. \cite[Proposition 3.6]{gubler4}).
We define a 
\emph{canonical measure}
on $X^{\an}$
to be
\[
\mu_{X^{\an}, L}
:=
\frac{1}{\deg_{L} X}
\cherncl_{1} ( \overline{L} |_{X} )^{\wedge d}
,
\]
which is a probability measure.
Although the canonical metric is unique only up to 
positive rational multiple,
the canonical measure is uniquely determined from $L$.

We fix the notation which is used in the remaining of this section.
Let $p^{\an} : E \to A^{\an}$ be the uniformization
and let $\val : E \to \RR^{n}$ be the valuation map,
where $n$ is the torus rank of $A$.
We set $\Lambda := \val ( \Ker p^{\an})$
and let $\overline{\val} : A^{\an}  \to \RR^{n} / \Lambda$ be
the valuation map.
Let $\CCC_{0}$ be a $\Lambda$-periodic rational polytopal decomposition
of $\RR^{n}$
and let $\overline{\CCC_{0}}$ be the rational polytopal decomposition
of $\RR^{n} / \Lambda$ 
induced
by quotient.
We take the Mumford model $\AAA_{0}$ of $A$ associated to $\overline{\CCC_{0}}$,
which can be defined over a finite extension of $\bfK^{\circ}$
(cf. Remark~\ref{rem:Mm}).
Let $\XXX_{0}$ be the closure of $X \subset A$ in
$\AAA_{0}$.
We
take a semistable alteration
$\varphi_{0} : \XXX' \to \XXX_{0}$ as in
Remark~\ref{assumption-sa}.
Put $X' := (\XXX')^{\an}$ and let $f : X' \to A^{\an}$ be the composite 
$(\XXX')^{\an}
 \to X^{\an} \to A^{\an}$.
Then
we have
a measure
\[
\mu_{X',f^{\ast} \overline{L}}
:=
\frac{1}{
\deg_{f^{\ast} L} X'}
\cherncl_{1} ( f^{\ast} \overline{L})^{\wedge d}
\]
on $X'$.

Since 
the morphism $X' \to X^{\an}$ is the analytification
of a morphism between projective varieties over a finite
extension of $\bfK$,
we have 
\addtocounter{Claim}{1}
\begin{align} \label{projectionformula2}
f_{\ast} \mu_{X',f^{\ast} \overline{L}} = \mu_{X^{\an} , L}
\end{align}
by the projection formula noted above\footnote{
We often identify $\mu_{X^{\an} , L}$
with its push-forward by the 
canonical closed immersion $X^{\an} \hookrightarrow A^{\an}$}.
Recall that
we have an expression
\addtocounter{Claim}{1}
\begin{align} \label{pullbackcanonicalmeasure}
\mu_{X',f^{\ast} \overline{L}}
=
\sum_{\Delta_{S}} r_{S} \delta_{\Delta_{S}}
\quad
(r_{S} > 0)
\end{align}
by \cite[Corollary 6.9]{gubler4},
where $\Delta_{S}$ runs through the set of non-degenerate
canonical simplices with respect to $f$,
and
$\delta_{\Delta_{S}}$ is the push-out of the
Lebesgue measure on $\Delta_{S}$.
Note that
all the coefficients $r_{S}$ are positive.
Let $S ( \XXX' )_{nd-f}$
be
the union of the non-degenerate canonical simplices of $S ( \XXX' )$
with respect to $f$.
Then (\ref{pullbackcanonicalmeasure})
shows that
$\mu_{X',f^{\ast} \overline{L}}$ is supported by
$S ( \XXX' )_{nd-f}$.
Since the notion of non-degeneracy of $\Delta_{S}$
is independent of
$L$, the support of $\mu_{X',f^{\ast} \overline{L}}$ is independent of
$L$.
Thus the support $S_{X^{\an}}$ of $\mu_{X^{\an} , L}$
is exactly the image of $S ( \XXX' )_{nd-f}$,
and in particular,
it 
does not depend of $L$.
Thus 
$S_{X^{\an}}$
depends only on $X$.
It
is called
the \emph{canonical subset} of $X^{\an}$ in \cite[Remark 6.11]{gubler4}.

By \cite[Theorem 6.12]{gubler4},
$S_{X^{\an}}$ has a canonical rational piecewise linear structure.
This
is characterized
by the property that,
for any model $\XXX$ of $X$
in a Mumford model
and
for any semistable alteration 
$\psi : \ZZZ \to \XXX$,
the induced map $f^{\an} |_{S ( \ZZZ )_{nd-f}} : 
S ( \ZZZ )_{nd-f} \to S_{X^{\an}}$
is a finite rational piecewise linear map
(cf. \cite[Theorem 6.12]{gubler4}).

If $\Delta_{S}$ is non-degenerate with respect to $f$,
then $f|_{\Delta_{S}} : \Delta_{S} \to f (\Delta_{S})$ is bijective.
Indeed, 
suppose that $\Delta_{S}$ is non-degenerate.
Since 
$\overline{f}_{\aff} |_{\Delta_{S}} : \Delta_{S} 
\to \RR^{n} / \Lambda$
is an affine map by Proposition~\ref{Val}
and since $\dim \overline{f}_{\aff} (\Delta_{S}) = \dim \Delta_{S}$,
the map
$\overline{f}_{\aff} |_{\Delta_{S}}$
is a finite affine map,
which means that
$\overline{f}_{\aff} |_{\Delta_{S}}$
is an injective affine map.
Since $\overline{\val} \circ f|_{\Delta_{S}} = \overline{f}_{\aff} |_{\Delta_{S}}$,
it follows that $f|_{\Delta_{S}} : \Delta_{S} \to f (\Delta_{S})$ is 
injective, and hence bijective.

\begin{Remark} \label{rem:ftofaff}
The valuation map $\overline{\val}$ restricts
to a finite piecewise linear map $\overline{\val} : S_{X^{\an}} \to 
\RR^{n} / \Lambda$.
Further,
for any non-degenerate stratum $\Delta_{S}$ with respect to $f$,
we have 
\[
\dim \Delta_{S} = 
\dim \overline{f}_{\aff} (\Delta ) = \dim f ( \Delta_{S} ).
\]
\end{Remark}

Since $S_{X^\an}$ has a canonical rational piecewise linear structure,
we can consider 
a {rational polytopal decomposition}
of $S_{X^\an}$.
In describing the canonical measure,
the following notion will be convenient.

\begin{Definition} \label{def:subdivisional}
A rational polytopal decomposition $\Sigma$ of $S_{X^{\an}}$ is said to be
\emph{$\varphi_{0}$-subdivisional}
if,
for any non-degenerate simplex $\Delta_{S}$ of $S (\XXX')$,
the image $f ( \Delta_{S} )$ is a finite union of polytopes in $\Sigma$.
\end{Definition}

\begin{Remark} \label{rem:subdivisional}
Let $\Sigma$ be a $\varphi_{0}$-subdivisional
rational polytopal decomposition of $S_{X^{\an}}$
and let
$\sigma \in \Sigma$ be a polytope.
Then there exists a non-degenerate canonical simplex $\Delta_{S}$
such that $f (\Delta_{S}) \supset \sigma$.
Since $\overline{f}_{\aff}|_{\Delta_{S}} 
= \overline{\val} \circ f |_{\Delta_{S}}$ is injective,
it follows that
$\overline{\val} |_{\sigma}$ is injective.
\end{Remark}

The following lemma gives us 
a
$\varphi_{0}$-subdivisional 
rational polytopal decomposition of $S_{X^{\an}}$.

\begin{Lemma} \label{lem:existence-subdivisional}
For any polytopal decomposition $\Sigma_{0}$ of $S_{X^{\an}}$,
there exists a 
rational polytopal
subdivision of $\Sigma_{0}$ which is a $\varphi_{0}$-subdivisional 
rational polytopal decomposition of $S_{X^{\an}}$.
\end{Lemma}

\begin{Pf}
Let $\Delta_S$ be a stratum non-degenerate with respect to$f$.
Since 
$f|_{\Delta_{S}} : \Delta_{S} \to f (\Delta_{S})$
is a bijective piecewise linear map,
there exists a rational subdivision $\Sigma_{S}'$ of $\Delta_{S}$
such that, for each $\sigma' \in \Sigma_{S}'$,
the restriction $f|_{\sigma'}$ is an affine map.
Then $f (\Sigma_{S}')$ is a rational polytopal decomposition
of $f (\Delta_{S})$.
Since $f : S(\XXX')_{nd-f} \to S_{X^{\an}}$ is 
a finite rational piecewise linear map,
we can take a rational polytopal subdivision $\Sigma$ of $\Sigma_{0}$
such that, for any non-degenerate $\Delta_{S}$,
$f (\Delta_{S} ) \cap \Sigma :=
\{ \sigma \in \Sigma \mid \sigma \subset f (\Delta_{S}) \}$
is a subdivision of $f (\Sigma_{S}')$.
Then we see that $\Sigma$ 
is a $\varphi_{0}$-subdivisional 
rational polytopal decomposition.
\end{Pf}

Let $\Sigma$ be a $\varphi_{0}$-subdivisional rational polytopal
decomposition of $S_{X^{\an}}$.
Then, we have an expression
\addtocounter{Claim}{1}
\begin{align} \label{canonicalmeasure-no-katachi}
\mu_{X^{\an} , L}
=
\sum_{\sigma \in \Sigma}
r'_{\sigma} \delta_{\sigma}
,
\end{align}
where $r'_{\sigma} \geq 0$ 
and $\delta_{\sigma}$ is the push-out of the Lebesgue measure on $\sigma$.
Indeed,
since $\Sigma$ is $\varphi_{0}$-subdivisional,
we have
$\delta_{\Delta_{S}} = \sum_{\sigma}
\delta_{(f|_{\Delta_{S}})^{-1}(\sigma)}$,
where 
$\sigma$ runs through the polytopes in $\Sigma$
such that $\sigma \subset f (\Delta_{S})$ and $\dim \sigma = \dim \Delta_{S}$.
Thus we write
$
f_{\ast} \delta_{\Delta_{S}}
=
\sum_{\sigma}
\alpha_{\sigma} \delta_{\sigma}
$
for some $\alpha_{\sigma} \geq 0$,
where $\sigma$ runs through polytopes in $\Sigma$
such that $\sigma \subset f (\Delta_{S})$.
It follows from (\ref{pullbackcanonicalmeasure}) 
and (\ref{projectionformula2}) that
we can write
$
\mu_{X^{\an} , L}=
f_{\ast} \mu_{X',f^{\ast} \overline{L}} = 
\sum_{\sigma \in \Sigma} r'_{\sigma} \delta_{\sigma}
$
for some $r'_{\sigma} \geq 0$.

\subsection{Strict supports of canonical measures} \label{ssca}

In this subsection,
we define the notion of strict supports and
investigate the strict supports of a canonical measure
on the canonical subset.
We follow the notation in \S~\ref{subsection:canonical}
for $A$, $L$, $X$, $\varphi_{0} : \XXX' \to \AAA_{0}$, $f :X' \to A^{\an}$,
$\mu_{X^{\an} . L}$ and so on.

\begin{Definition} \label{defSS}
Let $\PPP$ be a polytopal set with a finite 
polytopal decomposition $\Sigma$.
Let
$\mu$ be a semipositive Borel measure
on $\PPP$.
We say that $\sigma \in \Sigma$
is a \emph{strict support}
of
$\mu$
if there exists an $\epsilon > 0$ such that
$\mu - \epsilon \delta_{\sigma}$ is semipositive,
where $\delta_{\sigma}$ is the push-out of the Lebesgue measure on $\sigma$.
\end{Definition}

For example, let
$\Sigma$ be a $\varphi_{0}$-subdivisional rational polytopal
decomposition of $S_{X^{\an}}$.
Then
$\sigma \in \Sigma$
is a strict support 
of the canonical measure
$\mu_{X^{\an} , L}$
if and only if $r'_{\sigma} > 0$
in the expression (\ref{canonicalmeasure-no-katachi}).

\begin{Lemma} \label{nondeg-supp}
Let $\Sigma$ be a 
$\varphi_{0}$-subdivisional
rational polytopal decomposition
of $S_{X^{\an}}$
and take any
$\sigma \in \Sigma$.
Then, 
$\sigma$ is a strict support of $\mu_{X^{\an} , L}$
if and only if
there exists a canonical simplex $\Delta_{S}$ 
non-degenerate with respect to $f$ such that
$\dim \Delta_{S} = \dim \sigma$ and $\sigma \subset
f ( \Delta_{S} )$.
\end{Lemma}

\begin{Pf}
Recall that a canonical simplex $\Delta_{S}$ of $S(\XXX')$
is a strict support of $\mu_{X',f^{\ast} \overline{L}}$
if and only if $\Delta_{S}$ is non-degenerate with respect to $f$
(cf. (\ref{pullbackcanonicalmeasure})).
Since
$\Sigma$ is 
$\varphi_{0}$-subdivisional
and $\mu_{X^{\an} , L} = f_{\ast} \mu_{X', f^{\ast} \overline{L}}$,
it follows 
that 
$\sigma$ is a strict support of $\mu_{X^{\an} , L}$
if and only if
there exists a non-degenerate $\Delta_{S}$ such that
$\dim f (\Delta_{S}) = \dim \sigma$ and that $\sigma \subset
f ( \Delta_{S} )$.
Since 
$\dim f (\Delta_{S}) = \dim \Delta_{S}$
for non-degenerate $\Delta_{S}$
(cf. Remark~\ref{rem:ftofaff}),
that shows our lemma.
\end{Pf}

Recall that,
for
a rational point
$\overline{w} \in \overline{\val}(X^{\an})$,
$\initial_{\overline{w}} (X)$ denotes the initial degeneration of $X$ over 
$\overline{w}$.

\begin{Lemma} \label{uenolemma}
Let  $\Sigma$ be a 
$\varphi_{0}$-subdivisional
rational
polytopal decomposition of $S_{X^{\an}}$.
Suppose that $\sigma \in \Sigma$
is a strict support of $\mu_{X^{\an},L}$.
Then for any rational point $t \in \relin \sigma$,
there exists an
irreducible component $W$ of
$\initial_{\overline{\val}(t)} (X)$
such that
$t$ reduces to the generic point of $W$
(cf. Remark~\ref{rem:initialdegeneration}).
\end{Lemma}

\begin{Pf}
%
Since $\Sigma$ is $\varphi_{0}$-subdivisional
and $\sigma$ is a strict support of $\mu_{X^{\an},L}$,
Lemma~\ref{nondeg-supp} gives us 
a non-degenerate stratum $\Delta_{S}$ of $S ( \XXX')$
with respect to $f$
such that
$\dim \Delta_{S} = \dim \sigma$
and that $f ( \Delta_{S}) \supset \sigma$.
Note that $\relin \sigma \subset f (\relin \Delta_{S} )$,
and we
take a point
$u \in \relin \Delta_{S}$ with $f(u) = t$.
We
put $\overline{w} := \overline{\val} (t) = \overline{f}_{\aff} (u)$.

Let $\overline{\CCC}$ be a rational subdivision of $\overline{\CCC_{0}}$
such that 
$\overline{w}$ itself is a vertex of 
$\overline{\CCC}$.
Let $\AAA$ be the Mumford model of $A$ associated to $\overline{\CCC}$.
Let $\DDD$ be the subdivision of $S ( \XXX')$ given by
\[
\DDD = 
\left\{ 
\Delta_{S'} \cap
(\overline{f}_{\aff})^{-1} ( \overline{\Delta} )
\left|
S' \in \str ( \tilde{\XXX'}),
\overline{\Delta} \in \overline{\CCC}
\right.
\right\}
.
\]
Since $\Delta_{S}$ is non-degenerate, 
the map $\overline{f}_{\aff} |_{\Delta_{S}}$ is injective,
which show 
that $u$ is a vertex of $\DDD$.

Let 
$\XXX''$ be the formal model associated to  the subdivision $\DDD$
(cf. Remark~\ref{rem:Proposition5.5gubler4}).
Then we have an extension $\varphi' : \XXX'' \to \AAA$
of $f : X' \to A^{\an}$ by \cite[Proposition~5.13]{gubler4}.
Let
$R$ be the stratum of 
$\tilde{\XXX''}$
corresponding to $u$
in Proposition~\ref{gubler4Proposition5.7}.
Since $\{ \overline{w} \}$ is an open face of $\overline{\CCC}$
with
$\overline{f}_{\aff} (u) = \overline{w}$,
we have $\tilde{\varphi'} (R) \subset Z_{\overline{w}}$
by \cite[Proposition~5.14]{gubler4}.

Let $\XXX$ be the closure of $X$ in $\AAA$.
Since $\tilde{\varphi'} (R)$ is also contained in 
$\tilde{\XXX}$
and since $\initial_{\overline{w}}(X) = \tilde{\XXX}
\cap Z_{\overline{w}}$,
we have $\tilde{\varphi'} (R) \subset \initial_{\overline{w}} (X)$.
Let
$W$ be the closure of $\tilde{\varphi'} (R)$ in 
$\initial_{\overline{w}} (X)$.
It is an irreducible closed subset of $\initial_{\overline{w}} (X)$.
Since $\Delta_{S}$ is non-degenerate,
it follows from
Lemma~\ref{Snotuburekata} 
that
$\dim W = d : = \dim X$.
Since any irreducible component of $\initial_{\overline{w}} (X)$
has dimension $d$ (cf. Remark~\ref{initial-dim}),
$W$ is an irreducible component of 
$\initial_{\overline{w}} (X)$.

It only remains to show 
that $t$ reduces to the generic point of $W$.
By Lemma~\ref{vertex-genericpoint},
the point
$u \in (\XXX')^{\an}$ reduces
to the
generic point of
$R$.
Since
$\tilde{\varphi'}$ maps the generic point of $R$
to that of $W$,
it follows that
$t = f (u)$ reduces to the generic point of $W$.
\end{Pf}

Now we can show the following statement.

\begin{Proposition} \label{nondeg-supp2}
Let $\Sigma$ be a 
$\varphi_{0}$-subdivisional 
rational polytopal decomposition
of $S_{X^{\an}}$
and let
$\sigma \in \Sigma$ be a polytope.
Let $\Delta_S$ be a canonical simplex of $S (\XXX' )$.
Assume that there exists a
rational point $u \in \relin ( \Delta_{S} )$ with
$f ( u ) \in \relin (\sigma)$.
Then,
if 
$\sigma$ is a strict support of $\mu_{X^{\an} , L}$,
then
$\dim \overline{f}_{\aff} ( \Delta_{S}) = \dim \Delta_{S}$ holds.
\end{Proposition}

\begin{Pf}
Let $\sigma \in \Sigma$ be a strict support of $\mu_{X^{\an} , L}$.
We set $t := f (u)$ and $\overline{w} := \overline{f}_{\aff} (u)
=
\overline{\val}(f (u))$.
Note that they are rational points in $\relin \sigma$
and $\overline{\val}(X^{\an})$ respectively.
Since $\sigma$ is a strict support of $\mu_{X^{\an} , L}$,
Lemma~\ref{uenolemma}
gives us
an 
irreducible component $W$ of $\initial_{\overline{w}} (X)$
such that
$t$ reduces to the generic point of $W$.
Then
Lemma~\ref{mainproposition}
concludes that
$\dim \overline{f}_{\aff} ( \Delta_{S}) = 
\dim \Delta_{S}$.
\end{Pf}

The following lemma
gives us a sufficient condition 
for the assumption in Proposition~\ref{nondeg-supp2}.

\begin{Lemma} \label{rem:assumptionofprop}
Let $\Sigma$ be a $\varphi_{0}$-subdivisional
rational polytopal decomposition of $S_{\an}$
and let $\sigma \in \Sigma$ be a polytope.
Let $\Delta_{S}$ be a canonical simplex in $S( \XXX' )$ such that
$f ( \Delta_{S}) \supset \sigma$
and $\dim \overline{f}_{\aff} (\Delta_{S}) = \dim \overline{\val} (\sigma)$.
Then there exists a rational point $u \in \relin \Delta_{S}$
such that $f (u) \in \relin ( \sigma )$.
\end{Lemma}

\begin{Pf}
Since $\overline{\val}|_{\sigma} : \sigma \to \RR^{n} / \Lambda$ is
a piecewise linear map, there exists a polytope $\sigma' \subset \sigma$
with $\dim \sigma' = \dim \sigma$ such that $\overline{\val}$ is
affine over $\sigma'$.
Then $\overline{\val} ( \sigma' )$ is a polytope in $\RR^{n} / \Lambda$
such that $\relin \left( \overline{\val}(\sigma') \right) 
\subset \relin \left( \overline{f}_{\aff} (\Delta_S) \right)$,
so that there exists
a rational point $u \in \relin (\Delta_{S})$
such that $\overline{f}_{\aff} (u) \in \relin
\left( \overline{\val} ( \sigma' ) \right)$.
Since $\overline{\val}|_{\sigma'}$
is an injective affine map
(cf. Remark~\ref{rem:subdivisional}),
we have $f (u) \in \relin ( \sigma ') \subset \relin ( \sigma )$,
as required.
\end{Pf}

\section{Tropical triviality and density of small points}
\label{application}

\subsection{Notation, convention and remarks} \label{NT}

In the sequel, we fix the following notation and convention.
Let $k$ be an algebraically closed field,
$\mathfrak{B}$ a normal projective variety over $k$,
and let $\mathcal{H}$ be an ample line bundle on $\mathfrak{B}$.\footnote{
We assume $\mathfrak{B}$ to be a curve in \S~\ref{GBCforC}.}
Let
$K$ be the function field 
of 
$\mathfrak{B}$,
and let $\overline{K}$ be an algebraic closure of $K$.

For a finite extension $K'$ 
in
$\overline{K}$
of $K$,
let $\mathfrak{B}_{K'}$ denote the normalization of $\mathfrak{B}$ in $K'$.
Let $M_{K'}$ denote the set of points in $\mathfrak{B}_{K'}$ of codimension one.
For any $w \in M_{K'}$, the local ring $\OO_{\mathfrak{B}_{K'}, w}$ 
is a discrete valuation
ring having $K'$ as the fraction field,
and the order function $\ord_{w} : (K')^{\times} \to \ZZ$ gives an 
additive discrete valuation.
If $K''$ is a finite extension of $K'$, then
we have a canonical finite surjective morphism 
$\mathfrak{B}_{K''} \to \mathfrak{B}_{K'}$,
which induces a surjective map $M_{K''} \to M_{K'}$.
Thus we have an inverse system
$\left( M_{K'} \right)_{K'}$, where $K'$ runs through the finite extensions
of $K$ in $\overline{K}$.
We set $M_{\overline{K}} := \varprojlim_{K'} M_{K'}$,
and
call an element of $M_{\overline{K}}$ a \emph{place} of $\overline{K}$.

Each place $v = ( v_{K'} )_{K'} \in M_{\overline{K}}$  determines a unique
non-archimedean 
multiplicative absolute value $| \cdot |_{v}$ on $\overline{K}$
in such a way that
the following conditions are satisfied.
\begin{itemize}
\item
The restriction of $| \cdot |_{v}$
to $K'$ is equivalent to the absolute value
associated with the order function $\ord_{v_{K'}}$.
\item
For any $x \in K^{\times}$, $|x|_{v} = e^{ - \ord_{v_{K}} x}$.
\end{itemize}
Through this correspondence, we regard a place of $\overline{K}$ as
an absolute value of $\overline{K}$.
For a $v \in M_{\overline{K}}$,
let $\overline{K}_{v}$ denote the completion of $\overline{K}$
with respect to $v$.
It is an algebraically closed field complete with respect to
the non-archimedean absolute value $| \cdot |_{v}$.

For each $v_{K} \in M_{K}$,
let $| \cdot |_{v_{K} , \mathcal{H}}$ be the absolute value 
normalized in such a way that
\[
| x |_{v_{K} , \mathcal{H}} := e^{- (\ord_{v_{K}} x)( \deg_{\mathcal{H}} v_{K}) }
,
\]
where $\deg_{\mathcal{H}} v_{K}$
stands for the degree of the closure of $v_{K}$
in $\mathfrak{B}$ with respect to $\mathcal{H}$.
It is well known that the set 
$\mathfrak{V} :=
\{ | \cdot |_{v_{K} , \mathcal{H}} \}_{v \in M_{K}}$ 
of absolute values satisfies the product formula,
and hence we define the notion of heights
with respect to this set of absolute values,
namely, an absolute logarithmic height
with respect to $\mathfrak{V}$
(cf. \cite[Chapter~3 \S~3]{lang2}).
The ``height'' in this article
always means this height.

Let $F / k$ be any field extension.
For a scheme $X$ over $k$, we write 
$X_{F} := X \times_{\Spec k} \Spec F$.
If $\phi : X \to Y$ is a morphism of schemes over $k$,
we write $\phi_{F} : X_{F} \to Y_{F}$ for the base extension
to $F$.


Let $X$ be an algebraic scheme over $\overline{K}$.
For each place
$v$ of $\overline{K}$,
we have a Berkovich analytic space
associated to $X \times_{\Spec \overline{K}} \Spec \overline{K}_{v}$.
We write
$X_{v}$ for this analytic space.

Let $A$ be an abelian variety over $\overline{K}$
and suppose that $X$ is a subvariety 
of $A$.
Then $A$ and $X$
can be defined over a finite extension of $K$ in $\overline{K}$,
so that $A_{v}$ and $X_{v}$ can be defined over a subfield
of $\overline{K}_{v}$ 
over which the valuation
is a discrete valuation.
Thus the assumptions in \S~\ref{strictsupport}
are fulfilled for them,
and we can apply the arguments 
in \S~\ref{strictsupport}
in this setting.

\subsection{Tropically trivial subvarieties and density of small points}
\label{valuablevariety}

In this subsection,
we introduce the notion of tropically trivial subvarieties and
investigate 
in Theorem~\ref{maintheorem1}
the relationship between
tropical triviality and density of small points.

Let $A$ be an abelian variety over $\overline{K}$.
Let $L$ be an even ample line bundle on $A$.
Then we can define the
canonical height function
$\hat{h}_{L} : A ( \overline{K} ) \to \RR$
associated to $L$.
Refer to \cite{lang2} for the definition and properties.
Note that it is a non-negative function.

Let $X$ be a closed subvariety of $A$.
We say $X$ \emph{has dense small points} if 
$X (\epsilon ; L)$
is Zariski dense in $X$ for any $\epsilon > 0$
(cf. \cite[Definition~2.2]{yamaki5}),
where $X (\epsilon ; L)
:=
\{
x \in 
X ( \overline{K})
\mid
\hat{h}_{L} (x) \leq \epsilon
\}$.
Note
that
this notion does not depend on the choice of an even ample line bundle $L$.

From the viewpoint of the geometric Bogomolov conjecture
(cf. Conjecture~\ref{intGBCforAV}),
it is interesting to ask 
what properties a closed subvariety with dense small points
has.
Theorem~\ref{maintheorem1}
gives an answer to it.
For the statement, we make the following definition.

\begin{Definition} \label{tt}
Let $A$ be an abelian variety over $\overline{K}$ 
and let $X$ be a
closed subvariety of $A$.
We say $X$ is
\emph{tropically trivial}
if $\overline{\val}(X_{v} )$ consists of
a single point for any $v \in M_{\overline{K}}$,
where $\overline{\val}$ is the valuation map for $A_{v}$
(cf. \S~\ref{RtM}).
\end{Definition}

For a closed subvariety $X$ of $A$,
let $G_{X}$ be the stabilizer of $X$, i.e.,
$
G_{X} := 
\{ a \in A \mid
a + X \subset X
\}
$.
We regard it as a reduced closed subgroup scheme of $A$.

\begin{Theorem} \label{maintheorem1}
Let $A$ be an abelian variety over $\overline{K}$
and let $X$ be a closed subvariety of $A$.
If $X$ has dense small points, then
the closed subvariety
$X / G_{X}$ of $A /G_{X}$
is tropically trivial.
\end{Theorem}

\begin{Pf}
To argue by contradiction,
suppose that we have a counterexample, that is,
there exists a closed
subvariety
$X$ satisfying the following:
$X$ has dense small points,
and 
there exists a place $v \in M_{\overline{K}}$
such that $\overline{\val}((X / G_{X})_{v})$ is not a single point,
where $\overline{\val} : A_{v} \to \RR^{n} / \Lambda$
is the valuation map for $A_{v}$
(and $n$ is the torus rank of $A_{v}$)
(cf. \S~\ref{RtM}).
Then
$X / G_{X}$ is also a counterexample
by \cite[Lemma~2.1]{yamaki5},
so
we may and do assume that the stabilizer $G_{X}$ is trivial.
We consider a homomorphism $\alpha : 
A^{N} \to A^{N-1}$ defined by 
$\alpha : 
(x_{1} , \ldots , x_{N}) \mapsto (x_{2} - x_{1} , \ldots , x_{N} - x_{N-1})$.
Put $Z := X^{N} \subset A^{N}$ and $Y := \alpha (Z)$.
Since $G_{X} = 0$, the restriction
$Z \to Y$ of $\alpha$
is a generically finite surjective morphism
for a large $N$
(cf. \cite[Lemma~3.1]{zhang2}).
We fix such an $N$.
Let the same symbol $\alpha$
denote the morphism
$Z_{v} \to Y_{v}
$
between the associated analytic spaces over $\overline{K}_{v}$.

Let $\hat{h}_{A^{N}}$ and $\hat{h}_{A^{N-1}}$ be the canonical height functions
associated to 
any even ample line bundles on $A^{N}$ and $A^{N-1}$ respectively.
Since $X$ has dense small points,
so does $Z$ (cf. \cite[Lemma~2.4]{yamaki5}), and
hence $Z$ has
a generic net $(P_{m})_{m \in I}$,
where $I$ is a directed set,
such that
$\lim_{m} \hat{h}_{Z} (P_{m}) = 0$
(cf. \cite[Proof of Theorem~1.1]{gubler2}).
The image $(\alpha (P_{m}))_{m \in I}$
is also a generic net of $Y$ with
$\lim_{m} \hat{h}_{Y} (\alpha(P_{m})) = 0$.

Let $K'$ be a finite extension of $K$ in $\overline{K}$ 
over which
$A$ and $X$, and hence $Z$ and $Y$ are defined.
For a point $P$ in $Z ( \overline{K} )$ or in $Y ( \overline{K} )$,
let $O (P)$ denote the $\Gal ( \overline{K}/K')$-orbit of $P$.
Then, by the equidistribution theorem
\cite[Theorem 1.1]{gubler3},
we find
that
\[
\nu_{Z_{v},m} :=
\frac{1}{| O (P_{m}) |} \sum_{z \in O (P_{m})} \delta_{z}
\quad
\text{and}
\quad
\nu_{Y_{v},m} :=
\frac{1}{| O (\alpha(P_{m})) |} \sum_{y \in O (\alpha(P_{m}))} \delta_{y}
\]
weakly converge,
as $m \to \infty$,
to 
the canonical measures $\mu_{Z_{v}}$ on $Z_{v}$ and
to $\mu_{Y_{v}}$ on $Y_{v}$
associated to the even ample line bundles respectively.
Since $\alpha_{\ast} ( \nu_{Z_{v},m} ) = \nu_{Y_{v} , m}$,
we obtain
$
\alpha_{\ast} ( \mu_{Z_{v}} ) = \mu_{Y_{v}}
$.
Note in particular that $\alpha (S_{Z_{v}}) = S_{Y_{v}}$,
where $S_{Z_{v}}$ and $S_{Y_{v}}$ are the canonical subsets.

We take the Mumford models of $A_{v}^{N}$ and of $A_{v}^{N-1}$
associated to rational polytopal decompositions
of $(\RR^{n} / \Lambda)^{N}$ and
$(\RR^{n} / \Lambda)^{N-1}$ respectively.
Subdividing the rational polytopal decomposition
of $(\RR^{n} / \Lambda)^{N}$ if necessary,
we may and do assume that
the map $\alpha : A_{v}^{N} \to A_{v}^{N-1}$ extends to a morphism between the Mumford models (cf. \S~\ref{RtM2}).
Let $\ZZZ$ and $\YYY$ be the closure of $Z_{v}$ and $Y_{v}$
in these Mumford models of $A_{v}^{N}$ and of $A_{v}^{N-1}$
respectively.
Note that we have a morphism $\ZZZ \to \YYY$.
It follows from
Lemma~\ref{lem:semistable-alteration}
that
there exists a semistable alteration $\varphi : \ZZZ' \to \ZZZ$ 
for $\ZZZ$
as in Remark~\ref{assumption-sa}.
Note that the composite $\psi : \ZZZ' \to \ZZZ \to \YYY$ is a semistable
alteration for $\YYY$ as in Remark~\ref{assumption-sa}.
Let $g$ be the composite $(\ZZZ')^{\an} \to Z_{v} \hookrightarrow 
A_{v}^{N}$ and put $h := \alpha \circ g : (\ZZZ')^{\an} \to A_{v}^{N-1}$.

We note that $\alpha |_{S_{Z_{v}}} : S_{Z_{v}} \to  S_{Z_{v}}$
is a piecewise linear map.
Indeed,
we
consider the diagram
\[
\begin{CD}
S_{Z_{v}} @>{\alpha|_{S_{Z_{v}}}}>> S_{Y_{v}} \\
@V{\overline{\val}|_{S_{Z_{v}}}}VV @VV{\overline{\val}|_{S_{Y_{v}}}}V \\
\overline{\val}(S_{Z_{v}})
@>>{\overline{\alpha}_{\aff} |_{\overline{\val}(S_{Z_{v}})}}>
\overline{\val} (S_{Y_{v}})
,
\end{CD}
\]
where
$\overline{\alpha}_{\aff}$ is
the affine map associated to $\alpha$ (cf. \S~\ref{homo-Raynaud}).
Then $\overline{\alpha}_{\aff} |_{\overline{\val}(S_{Z_{v}})}$
is a piecewise linear map.
Since the columns are finite piecewise linear 
maps (cf. Remark~\ref{rem:ftofaff}),
it follows that the continuous map $\alpha |_{S_{Z_{v}}}$
is a piecewise linear map.

By Lemma~\ref{lem:existence-subdivisional},
there are
$\varphi$-subdivisional rational polytopal decomposition
$\Sigma_{Z_{v}}$ of $S_{Z_{v}}$
and a $\psi$-subdivisional rational polytopal decomposition
$\Sigma_{Y_{v}}$
of $S_{Y_{v}}$.
Taking a subdivision of $\Sigma_{Z_{v}}$ if necessary, we may assume 
that
for any $\sigma '\in \Sigma_{Z_{v}}$,
there exists a unique $\sigma'' \in \Sigma_{Y_{v}}$ such that 
$\relin  ( \alpha ( \sigma ')) \subset \relin ( \sigma'')$
and that $\alpha|_{\sigma'}$
is an affine map.
In particular,
for any $\sigma '\in \Sigma_{Z_{v}}$,
$\alpha (\sigma')$ is a polytope 
contained in some polytope in $\Sigma_{Y_{v}}$.

By \cite[Theorem 1.1]{gubler4},
$\overline{\val} (X_{v})$
coincides with the support of $\overline{\val}_{\ast} (\mu_{X_{v}})$.
Since 
it
contains a polytope of positive dimension
by our assumption at the beginning,
there exists a positive dimensional
rational
polytope $P$ in $\overline{\val}( X_{v})$
such that 
$\overline{\val}_{\ast} (\mu_{X_{v}}) - \epsilon \delta_{P}$ is semipositive for
a small $\epsilon > 0$.
By
\cite[Lemma 4.1 and Proposition 4.5]{yamaki5},
$\overline{\val}^{N}_{\ast} (\mu_{Z_{v}})$ is the product measure of $N$ copies of 
$\overline{\val}_{\ast} (\mu_{X_{v}})$.
Thus
$\overline{\val}^{N}_{\ast} (\mu_{Z_{v}}) - \epsilon \delta_{P^{N}}$
is semipositive
for a small $\epsilon > 0$,
where we recall that
$\overline{\val}^{m} : A_{v}^{m} \to (\RR^{n} / \Lambda)^{m}$
is the valuation map for $A_{v}^{m}$ for $m \in \NN$
(cf. \S~\ref{homo-Raynaud}).
Since $\Sigma_{Z_{v}}$ is $\varphi$-subdivisional,
we have an expression
\[
\overline{\val}^{N}_{\ast} ( \mu_{Z_{v}})
=
\sum_{\sigma'} 
r_{\sigma'}
\overline{\val}^{N}_{\ast} ( \delta_{\sigma'})
\quad
(r_{\sigma'} > 0)
,
\]
where $\sigma'$ runs through the polytopes in $\Sigma_{Z_{v}}$
that are strict supports of $\mu_{Z_{v}}$.
It follows 
that
there exists
a strict support
$\sigma \in \Sigma_{Z_{v}}$ of
$\mu_{Z_{v}}$
 such that
$\relin \left( \overline{\val}^{N} ( \sigma ) \right) 
\cap \relin ( P^{N} ) \neq \emptyset$,
where we remark that $\overline{\val}^{N} ( \sigma )$ is a polytope.
Note that
$\dim \left( \overline{\val}^{N} ( \sigma ) \right) = \dim P^{N}$.

Consider the affine map
 $\overline{\alpha}_{\aff} : (\RR^{n} / \Lambda)^{N} \to 
(\RR^{n} / \Lambda)^{N-1}$.
Since $\overline{\alpha}_{\aff} |_{P^{N}}$ 
contracts the diagonal
of $P^{N}$ to a point, it follows that
$\dim \left( \overline{\alpha}_{\aff} ( P^{N}) \right) < \dim 
\left( P^{N} \right)$,
and thus 
we have
$\dim \left( \overline{\val}^{N-1} (\alpha (\sigma)) \right) < 
\dim \left( \overline{\val}^{N} (\sigma) \right)$.
This inequality
with Remark~\ref{rem:ftofaff}
concludes
\addtocounter{Claim}{1}
\begin{align} \label{dimension-inequality}
\dim \alpha (\sigma) < \dim \sigma
.
\end{align}

We take a polytope $\tau \in \Sigma_{Y_{v}}$ 
with $\relin (\alpha (\sigma )) \subset \relin ( \tau )$.
By 
our assumption
on $\Sigma_{Z_{v}}$ and $\Sigma_{Y_{v}}$,
such a polytope $\tau$ uniquely exists
and 
is characterized by the condition that
$\relin (\alpha (\sigma )) \cap \relin \tau \neq \emptyset$.
We claim that
$\dim \tau = \dim \alpha (\sigma)$ and
that $\tau$
is a strict support of $\mu_{Y_{v}}$.
Indeed,
we take a compact subset $V \subset 
\relin (\alpha ( \sigma ))$ 
such that $V$ is a polytope
in $\alpha ( \sigma )$
with $\dim V = \dim \alpha ( \sigma )$.
Then there is an $\epsilon ' > 0$ such that
$\alpha_{\ast} \delta_{\sigma} - \epsilon' \delta_{V} \geq 0$.
Since $\sigma$ is a strict support of $\mu_{Z_{v}}$,
there exists an $\epsilon '' > 0$ such that
$\mu_{Z_{v}} - \epsilon '' \delta_{\sigma} \geq 0$.
Putting $\epsilon := \epsilon ' \epsilon ''$,
we obtain
\[
\mu_{Y_{v}} - \epsilon \delta_{V} \geq
\alpha_{\ast} \mu_{Z_{v}} - \epsilon '' \alpha_{\ast} \delta_{\sigma} 
=
\alpha_{\ast} ( \mu_{Z_{v}} - \epsilon '' \delta_{\sigma}) 
\geq 0
.
\]
Since
we can write
$
\mu_{Y_{v}} = \sum_{\sigma' \in \Sigma_{Y_{v}}} r_{\sigma'} \delta_{\sigma'}
$
with
$r_{\sigma'} \geq 0$
(cf. (\ref{canonicalmeasure-no-katachi})),
we find a strict support $\tau' \in \Sigma_{Y_{v}}$ 
of $\mu_{Y_{v}}$
such that
$\relin (V) \cap \relin (\tau') \neq \emptyset$
and $\dim V = \dim \tau'$.
Note that $\relin ( \alpha (\sigma) ) \cap \relin (\tau') \neq \emptyset$
in particular, for $\relin (V) \subset\relin( \alpha(\sigma))$.
Since $\tau$ is the unique polytope in $\Sigma_{Y_{v}}$ with 
$\relin ( \alpha (\sigma )) \cap \relin (\tau) \neq \emptyset$,
it follows that $\tau = \tau'$. 
Thus
$\tau$ is a strict support of $\mu_{Y_{v}}$.
Further, we have
$\dim \tau = \dim V= 
 \dim \alpha ( \sigma)$.

Since
$\sigma$ is a strict support of $\mu_{Z_{v}}$,
Lemma~\ref{nondeg-supp} gives us
a canonical simplex $\Delta_{S}$ of $S ( \ZZZ )$
non-degenerate with respect to $g$
such that $g ( \Delta_{S}) \supset \sigma$ and $\dim \Delta_{S} = \dim \sigma$.
Note that 
$\dim \overline{h}_{\aff} ( \Delta_{S} ) = \dim \tau$:
indeed,
we see that
\begin{align*} 
\dim \overline{h}_{\aff} ( \Delta_{S} ) 
= \dim \overline{\alpha}_{\aff} ( \overline{g}_{\aff} ( \Delta_{S} ))
= \dim \overline{\alpha}_{\aff} ( \overline{\val} ( \sigma ))
= \dim 
\overline{\val} ( \alpha (\sigma ))
=
\dim \alpha (\sigma)
=  \dim \tau
.
\end{align*}
Since $g ( \Delta_{S} ) \supset \sigma$,
we have
$h ( \Delta_{S} ) \supset \alpha (\sigma)$.
Since
$\relin ( \tau ) \supset \relin (\alpha (\sigma))$,
it follows that $h ( \Delta_{S} ) \cap \relin (\tau) \neq \emptyset$.
By the assumption that
polytopal decomposition
$\Sigma_{Y_{v}}$ is $\psi$-subdivisional,
this shows that $h ( \Delta_{S} ) \supset \tau$.
Since $\tau$ is a strict support of $\mu_{Y_{v}}$,
it follows from
Proposition~\ref{nondeg-supp2}
with Lemma~\ref{rem:assumptionofprop} that
$\dim \overline{h}_{\aff} ( \Delta_{S} ) = \dim \Delta_{S}$.
However,
the inequality (\ref{dimension-inequality}) shows
\begin{align*}
\dim
\overline{h}_{\aff} ( \Delta_{S} ) 
=
\dim \alpha (\sigma)
< 
\dim \sigma
=
\dim \Delta_{S}
.
\end{align*}
That is a contradiction. Thus
we complete the proof of the theorem.
\end{Pf}

\section{Main results} \label{mainresults}

In this section,
after recalling the geometric Bogomolov conjecture for
abelian varieties,
we 
prove 
Theorems~\ref{maintheoremint}
and \ref{maintheorem3int}
with the use of Theorem~\ref{maintheorem1}.

\subsection{Special subvarieties and  the geometric Bogomolov conjecture} \label{susecGBC}

Let $A$ be an abelian variety over $\overline{K}$.
Let
$\left( A^{\overline{K} / k}, \Tr_{A} \right)$ 
be the $\overline{K}/k$-trace of $A$,
that is,
the pair of an abelian variety
$A^{\overline{K} / k}$ over $k$
and a homomorphism
$\Tr_{A} : (A^{\overline{K} / k})_{\overline{K}} \to A $
with
the following universal property:
for any abelian variety $A'$ over $k$
and for any homomorphism $\phi : A'_{\overline{K}} \to A$,
there exists a unique homomorphism
$\Tr (\phi) : A' \to A^{\overline{K} / k}$ over $k$ such that
$\Tr_{A} \circ 
\Tr ( \phi)_{\overline{K}} = \phi$.
We refer to \cite{lang1} and \cite{lang2} for details.

We recall the notion of special subvarieties introduced in
\cite{yamaki5}.
Let $X$ be a closed subvariety of $A$.
We say that $X$ is \emph{special} if there exist a torsion point 
$\tau \in A( \overline{K} )$ and a  closed subvariety $X' \subset 
A^{\overline{K} / k}$ over $k$
such that
$
X = G_{X} + \Tr_{A} \left( X'_{\overline{K}} \right) + \tau
$,
where $G_{X}$ is the stabilizer of $X$.
We say a point $x \in A$ is a \emph{special point} if $\{ x \}$
is a special subvariety.
In the definition of special subvariety,
we may replace the condition for $\tau$ being torsion by 
that for $\tau$ being special
(cf. \cite[Remark~2.6]{yamaki5}).

We put together some properties of special subvarieties
in the following lemma.
Note that the equivalence between (a) and (b) in (1) below
shows that
the special subvariety
in the introduction is the same as the one recalled above.
Note also that 
\cite[Proposition~2.11]{yamaki5}
shows that
$X$ is a special subvariety of $A$ if and only if
$X / G_{X}$ is a special subvariety of $A /G_{X}$,
and Lemma~\ref{traceintro}~(1) generalizes this equivalence.

\begin{Lemma} \label{traceintro}
Let $A$ be an abelian variety over $\overline{K}$ and let
$X$ be a closed subvariety of $A$.
Further, let
$G'$ be a closed algebraic subgroup of $G_X$,
where $G_X$ is the stabilizer of $X$.
\begin{enumerate}
\item
The following are equivalent to each other.
\begin{enumerate}
\renewcommand{\labelenumi}{(\alph{enumi})}
\item
$X$ is a special subvariety.
\item
There exist 
an abelian variety $B$ over $k$, a closed subvariety $Y \subset B$,
a homomorphism $\phi : B_{\overline{K}} \to A$,
a torsion point $\tau \in A ( \overline{K} )$,
and
an abelian subvariety $A'$ of $A$ such that
$X = A' + \phi ( Y_{\overline{K}}) + \tau$.
\item
The quotient $X/G'$ is a special subvariety of $A/G'$.
\end{enumerate}
\item
Let $A_1$ be an abelian variety over $\overline{K}$ and
let $\psi : A \to A_1$ be a homomorphism.
Suppose that $X$ is a special subvariety of $A$.
Then $\psi (X)$ is a special subvariety of $A_1$.
\end{enumerate}
\end{Lemma}

\begin{Pf}
In this proof, we first show the equivalence between (a) and (b),
next show that (2) holds, and finally show the equivalence
between (a) and (c).

We begin by showing that (a) implies (b).
Suppose that
$X$ is a special subvariety.
Then by definition,
$
X = G_{X} + \Tr_{A} \left( X'_{\overline{K}} \right) + \tau
$
for some closed variety $X' \subset A^{\overline{K} / k}$
and some torsion point $\tau$.
Let $G_X^{\circ}$ be the connected component
of $G_X$ with $0 \in G_{X}$.
Since $X$ is irreducible,
we then have
$X = G_X^{\circ} + \Tr_{A} \left( X'_{\overline{K}} \right) + \tau$,
which shows (b).

To complete the proof of the 
equivalence between (a) and (b),
we show that (b) implies (a) next.
Write $X = A' + \phi ( Y_{\overline{K}}) + \tau$ as in (b).
Then 
$A' \subset G_X$,
and hence $X = G_{X} + \phi ( Y_{\overline{K}}) + \tau$.
Further, by the universality of the trace,
there exists a homomorphism $\Tr (\phi) : B \to A^{\overline{K} / k}$ 
such that $\phi = \Tr_{A} \circ \Tr (\phi)_{\overline{K}}$.
Putting $X' := \Tr (\phi) (Y)$,
we have
$
X = G_{X} + \Tr_{A} \left( X'_{\overline{K}} \right) + \tau
$, which concludes (a).

Here, using the equivalence between (a) and (b) shown above,
let us prove (2).
Since $X$ is special,
we write
$X = A' + \phi ( Y_{\overline{K}}) + \tau$ as in (b).
Then we have
$\psi (X) = \psi (A')
+ \psi \circ \phi ( Y_{\overline{K}})  + \psi ( \tau )$,
which is again an expression as in (b).
This shows
that 
$\psi (X)$ is special in $A_1$.

Let us move back to (1) and prove the remaining equivalences.
It follows immediately from (2) that
(a) implies (c).
Suppose that (c) holds.
Then by (2), $X/G_X$ is also special,
and hence \cite[Proposition~2.11]{yamaki5}
tells us that
$X$ is special.
\end{Pf}

\begin{Remark}
We consider the case where
$k$ is an algebraic closure of a finite field.
In \cite[Definition~2.1]{scanlon},
Scanlon defined ``special'' subvarieties,
from
which our special subvarieties are a little different.
Let
$X$ be a closed subvariety of $A$.
Then,
in the setting of this section,
one sees that $X$ is special 
{in the sense of Scanlon}
if and only if
there exist a closed subvariety $X'$ of $A^{\overline{K}/k}$ and a
(not necessarily torsion) point
$a \in A ( \overline{K})$
such that
$
X = G_{X} + \Tr_{A} \left( X'_{\overline{K}} \right) + a
$.
Thus
$X$ is special in our sense
if and only if
it is special in the sense of Scanlon and the point $a$ above
can be taken 
to be a torsion point.
\end{Remark}

As remarked in the introduction,
it is known that any special subvariety has dense small points.
The geometric Bogomolov conjecture claims that the converse 
should also hold:

\begin{Conjecture} [Geometric Bogomolov conjecture for
abelian varieties] \label{GBCforAV}
Let $A$ be an abelian variety over $\overline{K}$.
Let $X \subset A$ be a closed subvariety.
Then if $X$ has dense small points,
then
$X$ is a special subvariety.
\end{Conjecture}

\begin{Remark} \label{rem:special-dim0}
A point is special if and only if it is of height zero
(cf. \cite[(2.5.4)]{yamaki5}),
and hence $\{ x \}$ has dense small points 
if and only if $x$
is a special point.
This means that
Conjecture~\ref{GBCforAV} holds
if
$\dim X = 0$.
\end{Remark}


We make a remark on the relationship between special subvarieties
and tropically trivial subvarieties.
Since a special subvariety has dense small points,
it follows from Theorem~\ref{maintheorem1}
that,
if
$X$ is a special subvariety of $A$,
then 
$X / G_{X}$ is tropically trivial.
This assertion
itself can be shown directly,
not as a corollary of Theorem~\ref{maintheorem1}.
Indeed, taking the quotient by $G_{X}$,
we may assume $G_{X} = 0$.
Since $X$ is a special subvariety,
we translate $X$ by a special point
to
take
an abelian variety $A'$ over $k$,
a closed subvariety $Y \subset A'$ and
a homomorphism
$\alpha : A'_{\overline{K}} \to A$ such that
$\alpha (Y_{\overline{K}}) = X$.
Since $A'_{\overline{K}}$ 
has torus rank $0$ at any place,
the subvariety
$Y_{\overline{K}}$ is tropically trivial,
and hence its image $X$ is also tropically trivial
(cf. (\ref{tamanidetekuru})).

\subsection{Isogeny and the conjecture}

Let $\phi : A \to B$ be an isogeny of abelian varieties over $\overline{K}$
and let $X \subset A$ be a closed subvariety. 
By
\cite[Lemma 2.3]{yamaki5},
$X$ has dense small points if and only if the same property
holds for $\phi (X)$.
This suggests that,
if our formulation of the geometric Bogomolov conjecture is correct,
then
$X$ being special should be equivalent to $\phi (X)$ being special.
In fact, this holds true:

\begin{Proposition} \label{isog-special}
Let $\phi : A \to B$ be an isogeny of
abelian varieties over $\overline{K}$
and let
$X \subset A$ be a closed subvariety.
Then $X$ is a special subvariety if and only if 
$Y := \phi (X)$ is a special subvariety.
\end{Proposition}

\begin{Pf}
It follows from 
Lemma~\ref{traceintro}~(2) that
if $X$ is special, $Y$ is also special.
Let us show the other implication.
We note that $\phi (G_{X}) \subset G_{Y}$,
so that we have a homomorphism $\phi' : A / G_{X} \to B / G_{Y}$.
Furthermore, we see that $\phi'$ is an isogeny.
Indeed, 
let $a \in A (\overline{K})$ be a point with $\phi (a) \in G_{Y}$.
Then
for any $m \in \ZZ$, we have $\phi (ma + X ) = m \phi (a) + Y = Y$,
and thus $m a + X \subset \phi^{-1} (Y)$.
Since 
$\phi$ is an isogeny, 
the subset $ma + X$
is an irreducible component of 
$\phi^{-1} (Y)$ for any $m \in \ZZ$.
Since the number of irreducible
components of $\phi^{-1} (Y)$
is at most $\deg \phi$,
it follows that there exists a positive integer $m_{0}$ 
with $m_{0} \leq \deg \phi$
such that $m_{0} a + X = X$,
i.e., $m_{0} a \in G_{X}$.
This
shows that any element
of the kernel of $\phi' : A / G_{X} \to B / G_{Y}$ 
is a torsion point of order at most $\deg \phi$,
which concludes that 
$\phi'$ is an isogeny.

Now, suppose that $Y$ is special.
By the equivalence between (a) and (c) in 
Lemma~\ref{traceintro}~(1),
we see that $Y / G_{Y}$ is special and that
the
specialty of $X$ should follow from that of $X/G_{X}$.
Since $\phi' : A / G_{X} \to B / G_{Y}$ is an isogeny,
we may and do thus assume $G_{X}$ and $G_{Y}$ are trivial.
Further,
taking the translation by a torsion point,
we may assume $Y = \Tr_{B}
\left( Y'_{\overline{K}} \right)$
for some closed subvariety $Y' \subset B^{\overline{K}/k}$,
where $\left( B^{\overline{K}/k} , \Tr_{B} \right)$
is the $\overline{K} / k$-trace of $B$.
It follows from the universality of $\overline{K} / k$-trace
that there is
a unique homomorphism
$\Tr (\phi) : A^{\overline{K}/k} \to B^{\overline{K}/k}$
characterized by
$\Tr_{B} \circ \Tr (\phi)_{\overline{K}} =
\phi \circ \Tr_{A}$,
where $\left( A^{\overline{K}/k} , \Tr_{A} \right)$
is the $\overline{K} / k$-trace of $A$.
Since the
homomorphism $\Tr (\phi)$ is surjective by \cite[Lemma~1.5]{yamaki5},
we have $\Tr (\phi)  \left( \Tr (\phi)^{-1} (Y')) \right) = Y'$.
Note that
$
\Tr (\phi)^{-1} (Y')_{\overline{K}} 
=
\left( \Tr (\phi)_{\overline{K}} \right)^{-1} (Y'_{\overline{K}})$.
Then
we have
\[
\phi
\left(
\Tr_{A} \left( \Tr (\phi)^{-1} (Y')_{\overline{K}} \right)
\right)
=
\Tr_{B}
\left(
\Tr (\phi)_{\overline{K}}
\left(
\left(
\Tr (\phi)_{\overline{K}}
\right)^{-1}
\left(
Y'_{\overline{K}}
\right)
\right)
\right)
=
\Tr_{B}
\left(
Y'_{\overline{K}}
\right)
=
\phi (X)
.
\]
Since $\phi$ is an isogeny and $X$ is irreducible,
it follows that
there exists a torsion point $\tau \in A (\overline{K})$ such that
$Z :=X - \tau$ is an irreducible component of $\Tr_{A}
\left( \Tr (\phi)^{-1} (Y')_{\overline{K}} \right)$.
Thus there exists an irreducible component $W$ of
$\Tr (\phi)^{-1} (Y')_{\overline{K}}$
such that $\Tr_{A} (W) = Z$.
Since $k$ is algebraically closed, there exists an irreducible
component $W'$ of $\Tr (\phi)^{-1} (Y')$ such that $W'_{\overline{K}} = W$.
Then
$\Tr_{A} \left( W'_{\overline{K}} \right) + \tau 
= Z + \tau
= X$,
which concludes
that $X$ is a special subvariety.
\end{Pf}

\begin{Corollary} \label{GBC-isogeny}
Let $\phi : A \to B$ be an isogeny of abelian varieties 
over $\overline{K}$.
Then 
the geometric Bogomolov conjecture holds for $A$
if and only if that holds for $B$.
\end{Corollary}

\begin{Pf}
This assertion follows from
\cite[Lemma 2.3]{yamaki5}
and Proposition~\ref{isog-special}.
\end{Pf}

The following lemma will be used in the proof of 
Corollary~\ref{maintheorem3}

\begin{Lemma} \label{GBC-surjective2}
Let $\phi : A \to B$ be a surjective homomorphism of abelian varieties 
over $\overline{K}$.
If the geometric Bogomolov conjecture holds for $A$,
then it holds for $B$.
\end{Lemma}

\begin{Pf}
%
By Poincar\'e's complete reducibility theorem
(cf. \cite[Proposition~12.1]{milne} or \cite[\S~19, Theorem~1]{mumford}),
we
take
an abelian subvariety $A' \subset A$ such that
$\phi |_{A'} : A' \to B$ is an isogeny.
Note that
the geometric Bogomolov conjecture for $A$ implies that for
$A'$.
Thus the geometric Bogomolov conjecture for $A$ implies that for $A'$.
Furthermore,
Corollary~\ref{GBC-isogeny}
shows that the geometric Bogomolov conjecture for $A'$
implies that for $B$.
This concludes that the conjecture for $A$ implies that for $B$.
\end{Pf}

\subsection{Maximal nowhere degenerate abelian subvariety and
nowhere-degeneracy rank}  \label{ndrc}

In this subsection,
we define
the notion of maximal nowhere degenerate abelian subvarieties
and that of nowhere-degeneracy rank,
and we describe their properties.
In our main results on the geometric Bogomolov conjecture,
they play key roles.

An abelian variety $A$ over $\overline{K}$
is said to be \emph{somewhere degenerate}
if $A$
is degenerate at some place of $\overline{K}$,
and
is said to be \emph{nowhere degenerate}
if $A_{v}$ is non-degenerate 
for all $v \in M_{\overline{K}}$
(cf. \S~\ref{RtM}).


\begin{Lemma} \label{lem:ndas1}
Let $A_{1}$, $A_{2}$ and $A$ be abelian varieties over $\overline{K}$.
\begin{enumerate}
\item
Let 
$B_{1}$ and $B_{2}$ be nowhere degenerate abelian subvarieties
of $A_{1}$ and $A_{2}$
respectively.
Then $B_{1} \times B_{2}$ is a nowhere degenerate abelian subvariety
of $A_{1} \times A_{2}$.
\item
Let $\phi : A_{1} \to A_{2}$ be a homomorphism of abelian varieties.
Let
$B_{1}$ be a nowhere degenerate abelian subvariety of $A_{1}$.
Then $\phi (B_1)$ is a nowhere degenerate abelian subvariety of $A_2$.
\item
If $B$ and $B'$ are nowhere degenerate abelian subvarieties of $A$,
then so is $B + B'$.
\item
Let $B$ be an abelian subvariety of $A$.
Suppose that $A$ is nowhere degenerate.
Then $B$ is nowhere degenerate.
\end{enumerate}
\end{Lemma}

\begin{Pf}
For an abelian variety $A'$ over $\overline{K}$
and for a place $v \in M_{\overline{K}}$,
let $n ((A')_{v})$ denote the torus rank of $(A')_{v}$.

(1) 
For any $v \in M_{\overline{K}}$,
we have $n ( (B_{1} \times B_{2})_v ) = n ((B_{1})_{v}) + n ((B_{2})_{v})$
by Proposition~\ref{exactsequence-abelianrank},
and the right-hand side is $0$ by the assumption.
This shows that $B_{1} \times B_{2}$ is nowhere degenerate.

(2) 
For any $v \in M_{\overline{K}}$,
we have $n (\phi (B_{1})_{v}) \leq n ( (B_{1})_v )$
by Proposition~\ref{exactsequence-abelianrank}.
Since $n ( (B_{1})_v ) = 0$, we have $n (\phi (B_{1})_{v}) = 0$,
which shows that
$\phi (B_{1})$ is nowhere degenerate.

(3) Let $\phi : B \times B' \to A$ be the homomorphism
given by $\phi ( b , b') = b + b'$.
Then $\phi (B \times B') = B + B'$.
It follows from (1) and (2) that,
if $B$ and $B$ are nowhere degenerate, then
$B + B'$ is nowhere degenerate.

(4)
For any 
$v \in M_{\overline{K}}$,
we have $n (B_{v}) \leq n ( A_v )$
by Proposition~\ref{exactsequence-abelianrank}.
Since $n ( A_v ) = 0$, we have $n (B_{v}) = 0$,
which shows that
$B$ is nowhere degenerate.
\end{Pf}

The following lemma is used to make the key definitions.

\begin{Lemma} \label{lem:mndas}
Let $A$ be an abelian variety over $\overline{K}$.
Then there exists
a unique nowhere degenerate
abelian subvariety
$\mathfrak{m}$ of $A$ 
such that,
for any nowhere degenerate abelian subvariety
$B'$ of $A$, we have $B' \subset \mathfrak{m}$.
\end{Lemma}

\begin{Pf}
The uniqueness follows from the condition.
We show the existence.
Let $\mathfrak{m}$ 
be a nowhere degenerate abelian subvariety of $A$ of maximal dimension.
For any nowhere degenerate abelian subvariety
$B'$,
the sum
$\mathfrak{m} + B'$ is nowhere degenerate by Lemma~\ref{lem:ndas1} (3).
Since $\mathfrak{m}$ 
has maximal dimension among the nowhere degenerate abelian subvarieties,
we have $\dim (\mathfrak{m} + B') 
= \dim \mathfrak{m}$.
Thus we have $\mathfrak{m} + B' = \mathfrak{m}$ 
and hence $B' \subset \mathfrak{m}$.
\end{Pf}

By Lemma~\ref{lem:mndas},
the following definition makes sense:

\begin{Definition} \label{ndr}
Let $A$ be an abelian variety over $\overline{K}$.
\begin{enumerate}
\item
The unique abelian subvariety $\mathfrak{m}$ 
of $A$
in Lemma~\ref{lem:mndas} is called the
\emph{maximal nowhere degenerate abelian subvariety of $A$}.
\item
Let $\mathfrak{m}$ be the maximal nowhere degenerate abelian subvariety of $A$.
The \emph{nowhere-degeneracy rank} of $A$
is defined to be $\ndr (A) := \dim \mathfrak{m}$.
\end{enumerate}
\end{Definition}

Note that 
$\mathfrak{m} = 0$
if and only if $\ndr(A) = 0$.

\begin{Proposition} \label{prop:new:new-homomorphism-trivial}
Let $A_1$ and $A_2$ be abelian varieties over $\overline{K}$.
Let 
$\mathfrak{m}_{1}$ and $\mathfrak{m}_{2}$
be the maximal nowhere degenerate abelian
subvarieties of $A_{1}$ and $A_{2}$ respectively.
Let $\phi : A_{1} \to A_{2}$ be a homomorphism.
Suppose
that $\phi$ is surjective.
Then $\phi (\mathfrak{m}_{1}) = \mathfrak{m}_{2}$.
\end{Proposition}

\begin{Pf}
It follows from Lemma~\ref{lem:ndas1}~(2)
that $\phi (\mathfrak{m}_1 )$
is nowhere degenerate.
By the maximality of $\mathfrak{m}_2$,
we obtain
$\phi (\mathfrak{m}_1) \subset \mathfrak{m}_2$.

Let us prove the other inclusion.
By considering the homomorphism 
$A_{1} / \mathfrak{m}_{1} \to A_{2} / \phi (\mathfrak{m}_1)$
between the quotients 
instead
of $\phi : A_{1} \to A_{2}$,
we may and do assume $\mathfrak{m}_{1} = 0$.
Then 
our goal is to show that $\mathfrak{m}_{2} = 0$.
By the Poincar\'e complete reducibility theorem
and \cite[Remark in p.157]{mumford},
there exists
an abelian subvariety $A_{2}'$
of $A_{2}$ with an isogeny
$\alpha : A_{2} \to A_{2}' \times \mathfrak{m}_{2}$.
Let 
$p : A_{2}' \times \mathfrak{m}_{2} \to \mathfrak{m}_{2}$ be the second projection.
Let $\psi : A_{1} \to \mathfrak{m}_{2}$
be the composite of homomorphisms:
\[
\begin{CD}
A_{1}
@>{\phi}>>
A_{2}
@>{\alpha}>>
A_{2}' \times \mathfrak{m}_{2}
@>{p}>>
\mathfrak{m}_{2} 
.
\end{CD}
\]
Then $\psi$ is surjective,
and thus we are reduced
to show $\psi = 0$.

We recall that there exist non-trivial simple abelian varieties
$B_{1} , \ldots , B_{s}$ and an isogeny
$\gamma : B_{1} \times \cdots \times B_{s} \to A_{1}$.
Let $\mathfrak{n}$ be the maximal nowhere degenerate abelian
subvariety of $B_{1} \times \cdots \times B_{s}$.
We have $\gamma ( \mathfrak{n} ) \subset \mathfrak{m}_{1}$
by Lemma~\ref{lem:ndas1}~(2).
Since $\mathfrak{m}_{1} = 0$,
it follows that
$\gamma ( \mathfrak{n} ) =  0$.
Since
$\gamma$ is finite, we obtain $\mathfrak{n} = 0$.
It follows from Lemma~\ref{lem:ndas1}~(1)
that the maximal nowhere degenerate abelian subvariety of each
$B_{i}$ is trivial.
Thus
each $B_{i}$ is somewhere degenerate.

Let $\iota_{i} : B_{i} \to B_{1} \times \cdots \times B_{s}$ be the canonical
injection, and
set 
${\varphi}_{i} := \psi \circ \gamma \circ \iota_{i} : B_{i} \to  \mathfrak{m}_{2}$.
Note that $\psi = 0$
if $\psi \circ \gamma = 0$, 
and
that
$\psi \circ \gamma = 0$
if
${\varphi}_{i} = 0$ for any $i = 1, \ldots , s$.
Thus we only have to show that
${\varphi}_{i} = 0$ for any $i = 1, \ldots , s$.

\begin{Claim} \label{homomorphism-trivial}
Let $A$ be a somewhere degenerate simple abelian variety
and let $B$ be a nowhere degenerate abelian variety.
Then any homomorphism $\varphi : A \to B$
is trivial.
\end{Claim}

\begin{Pf}
Set $B' := \varphi (A)$.
By Lemma~\ref{lem:ndas1}~(4),
it is nowhere degenerate.
To argue by contradiction,
we suppose that
$B'$ is non-trivial.
Then,
since $A$ is simple,
the morphism $\varphi: A \to B'$ is an isogeny.
Therefore we have an isogeny $\varphi' : B' \to A$.
By Lemma~\ref{lem:ndas1}~(2), this shows that $A$ is nowhere degenerate,
which is a contradiction.
Thus $B' = \varphi (A)$ is trivial.
%
%
\end{Pf}

Applying
Claim~\ref{homomorphism-trivial}
to
$\varphi_{i} : B_{i} \to  \mathfrak{m}_{2}$
shows that $\varphi_{i}$ is trivial
for each $i = 1, \ldots , s$.
This
completes the proof of Proposition~\ref{prop:new:new-homomorphism-trivial}.
\end{Pf}

\begin{Corollary} \label{corollary-ndr}
Let $\phi : A_{1} \to A_{2}$ be a homomorphism
of abelian varieties.
If $\phi$ is surjective, then
$\ndr (A_{1}) \geq \ndr (A_{2})$.
\end{Corollary}

\begin{Pf}
Let 
$\mathfrak{m}_{i}$ be the maximal nowhere degenerate abelian subvariety
of $A_{i}$ for $i = 1,2$.
It follows from
Proposition~\ref{prop:new:new-homomorphism-trivial}
that
$\dim \mathfrak{m}_{1} \geq \dim \mathfrak{m}_{2}$.
Thus $\ndr (A_{1}) \geq \ndr (A_{2})$.
\end{Pf}

\begin{Proposition} \label{cor:new2}
Let $B_{1} , \ldots , B_{s}$ and $A$ be abelian varieties over $\overline{K}$.
Suppose that $A$ is isogenous to $B_{1} \times  \cdots \times B_{s}$.
Then $\ndr (A) = \ndr ( B_{1} ) +  \cdots + \ndr ( B_{s} )$.
\end{Proposition}

\begin{Pf}
We set $B := B_{1} \times  \cdots \times B_{s}$.
Since there are surjective homomorphisms
$A \to B$ and 
$B \to A$,
it follows from
Corollary~\ref{corollary-ndr}
that $\ndr (A) = \ndr ( B )$.
It remains to show that
$\ndr ( B ) = \ndr ( B_{1} ) +  \cdots + \ndr ( B_{s} )$.
Let $\mathfrak{n}$ be the 
maximal nowhere degenerate abelian subvariety of $B$
and $p_{i} : B \to B_{i}$ be the canonical projection.
Let $\mathfrak{n}_{i}$ be the maximal nowhere degenerate abelian subvariety
of $B_{i}$ for $i = 1 , \ldots , s$.
We prove $\ndr ( B ) = \ndr ( B_{1} ) +  \cdots + \ndr ( B_{s} )$
by showing
$\mathfrak{n} = \mathfrak{n}_{1} \times \cdots \times \mathfrak{n}_{s}$.
By Lemma~\ref{lem:ndas1}~(1),
$\mathfrak{n}_{1} \times \cdots \times \mathfrak{n}_{s}$
is a nowhere degenerate abelian subvariety of $B$. 
Then,  by the maximality of $\mathfrak{n}$, we have
$\mathfrak{n} \supset \mathfrak{n}_{1} \times \cdots \times \mathfrak{n}_{s}$.
On the other hand, 
Proposition~\ref{prop:new:new-homomorphism-trivial} gives us
$p_{i} (\mathfrak{n}) = \mathfrak{n}_{i}$
for all $i = 1 , \ldots , s$,
which shows
$\mathfrak{n} \subset \mathfrak{n}_{1} \times \cdots \times \mathfrak{n}_{s}$.
Thus we obtain
$\mathfrak{n} = \mathfrak{n}_{1} \times \cdots \times \mathfrak{n}_{s}$
and hence the proposition.
\end{Pf}

\subsection{
Proof of Theorem~\ref{maintheoremint}}

In this subsection, we establish Theorem~\ref{maintheoremint}
as a consequence of Theorem~\ref{prop:newmain-1}.
This theorem shows that
if a closed subvariety
of an abelian variety
has dense small points,
then it is contained
in 
the translate 
of the sum of the maximal nowhere degenerate abelian
subvariety and the identity component of its stabilizer,
by a special point.


We begin with two lemmas.
For a subvariety $X$ of an abelian variety $A$,
let $\langle X \rangle$ denote the abelian subvariety of $A$ generated by $X$,
that is, the smallest abelian subvariety containing $X$.

\begin{Lemma} \label{generated-subgroup}
Let $X$ be a closed subvariety 
of an abelian variety $A$ over $\overline{K}$.
Assume $0 \in X$.
Let $v$ be a place of $\overline{K}$.
Then $\langle X \rangle_{v}$ is the smallest analytic subgroup 
of $A_{v}$ containing 
$X_{v}$.
\end{Lemma}

\begin{Pf}
Let us consider, for each $l \in \NN$,  a morphism
$
X^{2l} \to A
$
given by
\addtocounter{Claim}{1}
\begin{align} \label{Xl}
( x_{1} , x_{2} , \ldots , x_{2l-1} , x_{2l} ) \mapsto
(x_{1} - x_{2}) + \cdots + (x_{2l-1} - x_{2l})
.
\end{align}
Let 
$X_{l}$ 
be the image of this morphism.
It is a closed subvariety
of $A$.
Note that $X_{l} \subset \langle X \rangle$.
Taking account that $0 \in X$,
we find
$
X \subset X_{1} \subset X_{2} \subset \cdots \subset X_{l} \subset \cdots
$.
Since each $X_{l}$ is an irreducible closed subset,
there exists $l_{0}$ such that $X_{l} = X_{l_{0}}$ for all $l \geq l_{0}$,
and thus
$\bigcup_{l \in \NN} X_{l} = X_{l_{0}}$.
We have
$X_{l} + X_{m} \subset X_{l+m}$, $0 \in X_{l}$, and $- X_{l} = X_{l}$ 
for all $l,m \in \NN$ by their definitions,
which tells us that $\bigcup_{l \in \NN} X_{l} = X_{l_{0}}$
is a subgroup scheme.
Since $X_{l_{0}}$ is irreducible, it follows that
$X_{l_{0}}$
is an abelian subvariety.
Since 
$X_{l_{0}} \subset \langle X \rangle$
and since
$\langle X \rangle$ is the smallest abelian subvariety containing $X$,
we obtain 
$\langle X \rangle = X_{l_{0}}$.

Let $B$ be an analytic subgroup of $A_{v}$
containing $X_{v}$.
We then have $B \supset (X_{l_{0}})_{v}$ by the definition of $X_{l_{0}}$
and hence
$B \supset \langle X \rangle_{v}$.
This shows that 
$\langle X \rangle_{v}$
is the smallest analytic subgroup containing $X_{v}$.
\end{Pf}

\begin{Lemma}\label{simple-valuable}
Let
$A$ be
an abelian variety  over $\overline{K}$
and let $X$ be a closed subvariety
of $A$.
Let $v$ be a place of $\overline{K}$.
Suppose that $A$ is simple
and is degenerate at $v$.
Then
if
$\overline{\val} (X_{v})$
consists of a single point, then
the same holds for
$X$,
where $\overline{\val}$ is the valuation map for $A_{v}$.
\end{Lemma}

\begin{Pf}
Taking the translation
of
$X$ by a point in $A(\overline{K})$,
we may assume that 
$0 \in X$.
Then
$\overline{\val} (X_{v}) = \{ \overline{\mathbf{0}} \}$
by assumption,
and hence
$X_{v} \subset
\overline{\val}^{-1} ( \overline{\mathbf{0}} )
$. 
Since $\overline{\val}^{-1} ( \overline{\mathbf{0}} )$ 
is an analytic subgroup of $A_{v}^{\circ}$
(cf. \S~\ref{RtM}), 
it follows from Lemma~\ref{generated-subgroup} that
$\langle X \rangle_{v} \subsetneq A_{v}$,
and thus
$\langle X \rangle \subsetneq A$.
Since $A$ is simple,
it follows that
$\langle X \rangle = \{ 0 \}$,
which shows
the lemma.
\end{Pf}

Now we show the following proposition.

\begin{Proposition} \label{general-valuable}
Let $X$ be a closed subvariety of 
an abelian variety $A$ over $\overline{K}$.
Then the following statements are equivalent to each other.
\begin{enumerate}
\renewcommand{\labelenumi}{(\alph{enumi})}
\item
$X $ is tropically trivial.
\item
There exists a point $a \in A(\overline{K})$ such that
$X \subset a + \mathfrak{m}$,
where $\mathfrak{m}$ is the maximal nowhere degenerate abelian subvariety
of $A$.
\end{enumerate}
\end{Proposition}

\begin{Pf}
We first show that (a) implies (b).
We may assume $\mathfrak{m} \neq A$.
By Poincar\'e's reducibility theorem with 
\cite[p.157, Remark]{mumford},
there exist an abelian subvariety $A'$ of $A$
and an isogeny $\beta : A \to A' \times \mathfrak{m}$.
Note that $\ndr (A') = 0$ by 
Proposition~\ref{cor:new2}.
Let $p : A' \times \mathfrak{m} \to A'$ be the canonical projection
and put $\psi := p \circ \beta : A \to A'$.

By the Poincar\'e complete reducibility theorem,
there exists an isogeny
$\nu : A' \to B'_{1} \times \cdots \times B'_{s}$,
where each $B_{i}'$ is a simple abelian variety.
Since $\ndr (A') = 0$,
it follows from
Proposition~\ref{cor:new2}
that $\ndr (B'_{i}) = 0$ for any $i = 1, \ldots , s$,
and thus
any $B'_{i}$
is degenerate at some place $v_{i}$.

Let $p_{i} : B'_{1} \times \cdots \times B'_{s} \to B'_{i}$
be the canonical projection for each $i$,
and put $\psi_{i} := p_{i} \circ \nu \circ \psi : A \to B'_{i}$.
Since $\overline{\val} (X_{v_{i}})$ is a single point by our assumption,
so is
$\overline{(\psi_{i})}_{\aff} 
(\overline{\val} (X_{v_{i}}))$,
where $\overline{(\psi_{i})}_{\aff}$ is the affine map
 associated to $\psi_{i}$
(cf. \S~\ref{homo-Raynaud}).
Since $\overline{\val} (\psi_{i} (X)_{v_{i}}) = 
\overline{(\psi_{i})}_{\aff} 
(\overline{\val} (X_{v_{i}}))$
(cf. (\ref{tamanidetekuru})),
we see that $\overline{\val} (\psi_{i} (X)_{v_{i}})$ is a single point.
Since $B'_{i}$ is simple and degenerate at $v_{i},$
it follows from Lemma~\ref{simple-valuable}
that $\psi_{i} (X)$ is a single point for any $i$,
which implies that
$\nu ( \psi (X))$ is a single point.
Since $\nu$ is finite and $\psi (X)$ is connected,
we see that $\psi (X)$ is a single point,
and thus we write 
$\psi (X) = \{ b' \}$ for some $b' \in B' (\overline{K})$.

We then have $p^{-1} (\psi (X)) = \{ b' \} \times \mathfrak{m}$,
and hence 
\[
\psi^{-1} (\psi (X)) = \beta^{-1} (\{ b' \} \times \mathfrak{m})
= \bigcup_{\sigma \in \beta^{-1} (b' , 0)} (\sigma + \mathfrak{m})
,
\]
which is a finite union of subvarieties.
Since $X \subset \bigcup_{\sigma \in \beta^{-1} (b' , 0)} (\sigma + \mathfrak{m})$
and $X$ is irreducible,
there exists a point $a \in \beta^{-1} (b' , 0)$ with
$X \subset a + \mathfrak{m}$.
Thus we obtain (b).

Next we show that (b) implies (a).
Since $\mathfrak{m}$ is nowhere degenerate, 
the valuation map for $\mathfrak{m}$ is trivial.
Take any $v \in M_{\overline{K}}$.
By the compatibility
of valuation maps with homomorphisms (cf. (\ref{tamanidetekuru})),
we have
$\overline{\val} (\mathfrak{m}_v) = \{ \overline{\mathbf{0}} \}$,
where $\overline{\val} : A_v \to \RR^{n} / \Lambda$
is the valuation map for $A_{v}$.
Then
$\overline{\val} (X_v) \subset \{ \overline{\val} (a) \}$,
and thus $\overline{\val} (X_v)$
is a single point.
Since $v$ is arbitrary,
this concludes that
$X$ is tropically trivial.
\end{Pf}

Using Theorem~\ref{maintheorem1}
with Proposition~\ref{general-valuable},
we show the following theorem.
Here, let $G_X^{\circ}$ denote the connected component
of the stabilizer $G_X$ with $0 \in G_X^{\circ}$.

\begin{Theorem} \label{prop:newmain-1}
Let $A$ be an abelian variety over $\overline{K}$
with maximal nowhere degenerate abelian subvariety $\mathfrak{m}$
and let $X$ be a closed subvariety of $A$.
Suppose that $X$ has dense small points.
Then there exists a special point $a_0 $ of $A$
such that $X \subset a_0 + G_X^{\circ} +  \mathfrak{m}$.
\end{Theorem}

\begin{Pf}
Let
$\phi : A \to  A/G_X$ be the quotient homomorphism.
Then $\phi (X) = X / G_{X} $
has trivial stabilizer.
Since $X$ has dense small points,
it follows that
$\phi (X)$ is tropically trivial by Theorem~\ref{maintheorem1}.
Note that by Proposition~\ref{prop:new:new-homomorphism-trivial},
$\phi ( \mathfrak{m} )$
is 
the maximal nowhere degenerate abelian subvariety
of $A / G_{X}$.
By Proposition~\ref{general-valuable},
it follows
that
there exists a point $a_0' \in A  ( \overline{K})$
such that
$\phi (X) \subset \phi (a_0') + \phi (\mathfrak{m})$,
and hence
$X \subset a_0' + G_X + \mathfrak{m}$.
Since $X$ is irreducible, 
$X \subset a_0' + G_X' + \mathfrak{m}$ for some connected component
$G_X'$
of $G_X$.
Note that $G_X' = b + G_X^{\circ}$ for some torsion point
$b \in A ( \overline{K} )$.
Setting
$a_0 := a_0' + b$,
we then have
$X \subset a_0 + G_X^{\circ} + \mathfrak{m}$.

Now it suffices to show that
$a_0$ can be chosen to be a special point.
Let $\psi : A  \to A / ( G_{X}^{\circ} + \mathfrak{m} )$ be the
quotient. Then $\psi (X) = \{ \psi (a_0) \}$ has dense small points
by \cite[Lemma~2.1]{yamaki5},
and thus $\psi (a)$ is a special point of $
A / ( G_{X}^{\circ} + \mathfrak{m} )
$ (cf. Remark~\ref{rem:special-dim0}).
By \cite[Lemma~2.10]{yamaki5},
there exists a special point $a_0'$ of $A$ with $\psi (a_0 ') = \psi (a_0)$.
Replacing $a_0$ with $a_0'$,
we have $X \subset  a_0 + G_X^{\circ} + \mathfrak{m}$ with $a_0$ special.
Thus, we complete the proof of the theorem.
\end{Pf}

As a consequence, we obtain the following two corollaries.
The first one gives us a sufficient condition
for a closed subvariety having dense small points
to be the translate of an abelian subvariety by a special point.

\begin{Corollary} \label{maintheorem2}
Let $A$ be an abelian variety over $\overline{K}$
and let $X$ be a closed subvariety of $A$.
Let $G_{X}$ be the stabilizer of $X$ in $A$.
Assume that 
$\dim X / G_{X} \geq \ndr ( A / G_{X} )$.
Then,
if $X$ has dense small points,
then
there exists a special point $x_{0} \in X (\overline{K})$ such that
$X = x_{0} + G_{X}$.
\end{Corollary}

\begin{Pf}
Let $\phi : A \to  A / G_{X}$ be the quotient homomorphism
and
let $\mathfrak{m}$ 
be the maximal nowhere degenerate
abelian subvariety of $A$.
Then $\phi (\mathfrak{m})$
is the maximal nowhere degenerate
abelian subvarieties of $A / G_X$ by
Proposition~\ref{prop:new:new-homomorphism-trivial}.
Suppose that $X$ has dense small points.
By Theorem~\ref{prop:newmain-1},
there exists a special point $x_0$ of $A$ such that
$X / G_X = \phi (X) \subset \phi (x_{0} ) + \phi (\mathfrak{m})$.
Since $\dim ( X / G_X ) \geq \ndr ( A / G_X ) = \dim ( \phi (\mathfrak{m}) )$,
it follows that $X / G_X = \phi (x_{0} ) + \phi (\mathfrak{m})$.
Since $X /G_X$ has trivial stabilizer,
we find $X / G_X = \{ \phi (x_0) \}$,
which
shows
$X = x_{0} + G_{X}$.
Further, this equality also shows that $x_{0} \in X ( \overline{K} )$.
Thus we obtain the corollary.
\end{Pf}



The second corollary shows Theorem~\ref{maintheoremint},
which gives us a  partial answer to the geometric Bogomolov conjecture.

\begin{Corollary}[cf. Theorem~\ref{maintheoremint}] \label{GBCfornondeggeq1}
Let $A$ be an abelian variety over $\overline{K}$
with
$\ndr (A) \leq 1$
and
let $X$ be a closed subvariety of $A$.
Suppose that 
$X$ has dense small points.
Then
$X$ is the translate of an abelian subvariety by a special point.
In particular,
it is a special subvariety.
\end{Corollary}

\begin{Pf}
Let $\phi : A \to A/G_X$ be the quotient.
Since $X$ has dense small points,
Theorem~\ref{prop:newmain-1} gives us a special point 
$x_0 \in X ( \overline{K} )$
such that $X \subset x_0 + G_X + \mathfrak{m}$,
where $\mathfrak{m}$ is the maximal nowhere degenerate abelian subvariety
of $A$.
We then have $
X / G_X \subset \phi (x_0) + \phi ( \mathfrak{m} )$.

Note that if we show $\dim X/G_X = 0$,
then $X/G_X = \{ \phi (x_0) \}$
and hence
$X = x_0 + G_X$.
Further
since $X$ is irreducible, $G_X$ is an abelian subvariety.
Therefore
it suffices to claim that $\dim X/G_X = 0$.
Since
$\ndr (A) \leq 1$,
$\phi ( \mathfrak{m} )$ has dimension $0$ or $1$.
If $\dim \phi ( \mathfrak{m} ) = 0$,
then we are done.
Suppose that
$\dim \phi ( \mathfrak{m} ) = 1$.
Since $X / G_X$ has trivial stabilizer, we find that 
$X/G_X \subsetneq \phi (x_0) + \phi ( \mathfrak{m} )$,
which 
implies
$\dim X/G_X = 0$.
Thus 
we obtain the corollary.
\end{Pf}


\begin{Remark} \label{rem:generalization}
Let $A$ be an abelian variety over $\overline{K}$.
\begin{enumerate}
\item
Suppose that
$b(A_{v}) \leq 1$ for some place $v$,
where $b(A_{v})$ is the abelian rank of $A_{v}$.
Since $b(A_{v}) \geq \ndr (A)$ in general,
Corollary~\ref{GBCfornondeggeq1}
generalizes
Theorem~\ref{b=1}
and hence Theorem~\ref{thm:gubler}.
\item
Suppose that $A$ is simple and somewhere degenerate.
Then $\ndr (A) = 0$ and hence the geometric Bogomolov conjecture
holds for $A$.
\end{enumerate}
\end{Remark}

\subsection{Reduction to the
nowhere degenerate case} \label{reduction}

In this subsection,
we reduce the geometric Bogomolov conjecture for an abelian variety
to the conjecture for its maximal nowhere degenerate abelian 
subvariety
(Theorem~\ref{maintheorem3int}).

We first establish
the following theorem.

\begin{Theorem} \label{prop:new-key60}
Let $A$
be an abelian variety over $\overline{K}$,
$B$ a nowhere degenerate abelian variety over $\overline{K}$,
and
let $\phi : A \to B$ be a surjective homomorphism
with $\dim B = \ndr ( A )$.
Let $X$ be a closed subvariety of $A$.
Suppose that
$X$ has dense small points and that
$\phi (X)$ is a special subvariety.
Then $X$ is a special subvariety.
\end{Theorem}

\begin{Pf}
Let $G_X^{\circ}$ be the connected component of the stabilizer
$G_X$ of $X$ with $0 \in G_X^{\circ}$.
Let $\mathfrak{m}$ be the maximal nowhere degenerate
abelian subvariety of $A$.
Note that the restriction $\phi |_{\mathfrak{m}} : \mathfrak{m}
\to B$ is an isogeny.
Indeed,
we have
$\dim \mathfrak{m} =\dim B$ by assumption,
and
$\phi ( \mathfrak{m} ) = B$
by Proposition~\ref{prop:new:new-homomorphism-trivial}.

\emph{Step 1. First we show the assertion under the
assumption that $G_X^{\circ}  = 0$.
} 
Under this assumption,
Theorem~\ref{prop:newmain-1}
gives us
a special point $a_0$ of $A$ such that
$X \subset a_0 + \mathfrak{m}$.
We put $X' = X - a_0$.
Then $X' \subset \mathfrak{m}$, and $X'$ has dense small points.
On the other hand,
since $\phi (a_0)$ is a special point of $B$,
$\phi|_{\mathfrak{m}} (X') = \phi (X) - \phi (a_0)$ is also special
(cf. \cite[Remark~2.6]{yamaki5}).
By Proposition~\ref{isog-special},
it follows that
$X'$
is special.
Since $a_0 $ is a special point,
$X = a_0 + X'$ is special
(cf. \cite[Remark~2.6]{yamaki5}).

\emph{Step 2. Next we conclude the theorem
by reducing the general case to Step~1.}
The quotient $B / \phi (G_X^{\circ})$ is nowhere degenerate
by Lemma~\ref{lem:ndas1}~(2). 
We would like to
show that
$\ndr ( A / G_X^{\circ} ) = \dim (B / \phi (G_X^{\circ}))$.

First, we
note that
\addtocounter{Claim}{1}
\begin{align} \label{eq:new60}
\dim ( G_X^{\circ} \cap \mathfrak{m} )
= \dim \phi (G_X^{\circ}).
\end{align}
Indeed, let $\mathfrak{n}$ be the maximal nowhere degenerate abelian
subvariety of $G_X^{\circ}$ and let 
$\left( G_X^{\circ} \cap \mathfrak{m} \right)^{\circ}$
be the connected component
of $G_X^{\circ} \cap \mathfrak{m}$ containing 
with $0 \in
\left( G_X^{\circ} \cap \mathfrak{m} \right)^{\circ}$.
By Lemma~\ref{lem:ndas1}~(4),
we then have $
\left( G_X^{\circ} \cap \mathfrak{m} \right)^{\circ}
\subset \mathfrak{n}$.
Further,
the maximality of $\mathfrak{m}$
tells us
$\mathfrak{n} \subset \mathfrak{m}$
and hence $\mathfrak{n} 
\subset 
\left( G_X^{\circ} \cap \mathfrak{m} \right)^{\circ}
$.
Thus $
\left( G_X^{\circ} \cap \mathfrak{m} \right)^{\circ}
= \mathfrak{n}$.
On the other hand,
$\phi (G_X^{\circ})$ is nowhere degenerate
by Lemma~\ref{lem:ndas1}~(4).
Therefore,
by Proposition~\ref{prop:new:new-homomorphism-trivial},
we have
$\phi|_{\mathfrak{m}} 
\left( \left( G_X^{\circ} \cap \mathfrak{m} \right)^{\circ} 
\right) =
\phi |_{\mathfrak{m}} ( \mathfrak{n} ) =
\phi (G_X^{\circ})$.
Since $\phi|_{\mathfrak{m}}$ is a finite morphism,
that shows (\ref{eq:new60}).

Let $\alpha : A \to A / G_X^{\circ}$ 
and $\beta : B \to B / \phi (G_X^{\circ})$
denote the quotient homomorphisms
and
let
$\phi' : A / G_X^{\circ} \to B / \phi (G_X^{\circ})$
be the homomorphism induced from $\phi$.
Remark that
Proposition~\ref{prop:new:new-homomorphism-trivial} tells us that
$\alpha ( \mathfrak{m} )$
is the maximal nowhere degenerate abelian subvariety
of $A / G_X^{\circ} $,
and hence $\ndr (A / G_X^{\circ}) = \dim \alpha ( \mathfrak{m} )$.
We
consider
the homomorphism
$\phi'|_{\alpha ( \mathfrak{m} )} 
: \alpha ( \mathfrak{m} ) \to B / \phi (G_X^{\circ})$.
Since
$\alpha ( \mathfrak{m} ) \cong 
\mathfrak{m} / ( G_X^{\circ} \cap \mathfrak{m})$, 
we have $\dim \alpha ( \mathfrak{m} )
=
\ndr (A) - \dim ( G_X^{\circ} \cap \mathfrak{m})$.
Since $\ndr (A) = \dim B$
and since we have (\ref{eq:new60}),
that
equals
$\dim B - \dim \phi (G_X^{\circ} )$.
Therefore
$\dim \alpha ( \mathfrak{m} ) = \dim B / \phi (G_X^{\circ})$,
and thus
$
\ndr (A / G_X^{\circ} ) 
=  \dim (B / \phi (G_X^{\circ}) )
$.

Now we can apply Step~1
to $\phi' : A / G_X^{\circ} \to B / \phi (G_X^{\circ}) $.
Indeed, we 
know
$
\ndr (A / G_X^{\circ} ) 
=  \dim (B / \phi (G_X^{\circ}) )
$.
Since $X$ has dense small points,
so does
$\alpha (X)$.
Since $\phi (X)$ is special,
$\phi' ( \alpha (X) ) = \beta ( \phi (X) )$ is also special
by Lemma~\ref{traceintro}~(2).
Further
$\alpha (X)$ has 
stabilizer $G_X / G_X^{\circ}$, which has dimension $0$.
Now applying Step~1,
we see that $X / G_X^{\circ} = \alpha (X)$ is a special subvariety.
By Lemma~\ref{traceintro}~(1),
we conclude that $X$ is a special subvariety.
\end{Pf}

As a consequence, we obtain the following corollary
 (Theorem~\ref{maintheorem3int}).

\begin{Corollary} [Theorem~\ref{maintheorem3int}] \label{maintheorem3}
Let $A$ be an abelian variety and let $\mathfrak{m}$ be 
the maximal nowhere degenerate abelian subvariety of $A$.
Then the following 
are equivalent to each other:
\begin{enumerate}
\renewcommand{\labelenumi}{(\alph{enumi})}
\item
The geometric Bogomolov conjecture holds for $A$;
\item
The geometric Bogomolov conjecture holds for $\mathfrak{m}$.
\end{enumerate}
\end{Corollary}

\begin{Pf}
By the Poincar\'e complete
reducibility theorem with \cite[p.157, Remark]{mumford},
we have a surjective homomorphism
$\phi : A \to \mathfrak{m}$.
Then
(a) implies (b) by
Lemma~\ref{GBC-surjective2}.
To see the other direction,
suppose (b)
and
let $X$ be a closed subvariety of $A$ having dense small points.
Then $\phi (X)$ has dense small points.
Since the conjecture holds for $\mathfrak{m}$,
$\phi (X)$ is a special subvariety of $\mathfrak{m}$.
By Theorem~\ref{prop:new-key60},
we conclude that $X$ is special.
Thus 
(b) implies (a).
\end{Pf}

We remark that
Theorem~\ref{maintheoremint} follows also from
the following corollary,
for the geometric Bogomolov conjecture holds for elliptic curves.

\begin{Corollary}
\label{general-nondeg}
For any non-negative integer $s$, 
the following are equivalent to each other:
\begin{enumerate}
\renewcommand{\labelenumi}{(\alph{enumi})}
\item
Conjecture~\ref{GBCforAV}
holds for any abelian variety $A$
with $\ndr (A) \leq s$;
\item
Conjecture~\ref{GBCforAV}
holds for any nowhere
degenerate abelian variety $B$ 
with $\dim B \leq s$.
\end{enumerate}
\end{Corollary}

\begin{Pf}
It is trivial that (a) implies (b).
Suppose that (b) holds.
Let $A$ be an abelian variety over $\overline{K}$ with 
$\ndr (A) \leq s$.
Let $\mathfrak{m}$ be the maximal nowhere degenerate abelian
subvariety of $A$.
Then $\dim \mathfrak{m} \leq s$, so that
Conjecture~\ref{GBCforAV} holds for $\mathfrak{m}$ by assumption.
It follows from Corollary~\ref{maintheorem3} that
Conjecture~\ref{GBCforAV} holds for $A$.
\end{Pf}

Finally in this section, we remark that
by
Corollary~\ref{general-nondeg},
Conjecture~\ref{GBCforAV}
is equivalent to
the following conjecture.

\begin{Conjecture} [Geometric Bogomolov conjecture for 
nowhere degenerate abelian varieties] \label{GBCforAVwithtoutdeg}
Let $A$ be a nowhere degenerate abelian variety over $\overline{K}$
and let
$X$ be a closed subvariety of $A$.
Then if
$X$ has dense small points, then
it should be a special subvariety.
\end{Conjecture}


\section{Geometric Bogomolov conjecture for curves} \label{GBCforC}

In this section, we consider
the geometric Bogomolov conjecture
for curves.
This conjecture claims that 
the set of $\overline{K}$-points of
a non-isotrivial smooth projective curve of genus at least $2$
is ``discrete''
in its Jacobian
with respect to the N\'eron--Tate semi-norm.
To be precise,
let $C$ be a smooth projective curve over $\overline{K}$
of genus $g \geq 2$
and
let $J_{C}$ be its Jacobian variety.
For a divisor $D$ of degree $1$
on $C$,
let
$
j_{D} : C \to J_{C}
$
be the embedding given by $x \mapsto [x - D]$,
where $[x - D]$ denotes the divisor class of $x - D$.
Let $|| \cdot ||_{NT}$ be the canonical N\'eron--Tate semi-norm
arising from the canonical N\'eron--Tate pairing on
$J_{C}$.
The following assertion is called the
geometric Bogomolov conjecture for curves.
Here $C$ is said to be \emph{isotrivial} if there is a curve $C'$ over $k$
such that 
$C'_{\overline{K}} \cong C$.

\begin{Conjecture} \label{GBCforcurves}
Assume that $C$ is non-isotrivial.
For any divisor $D$ of degree $1$ on $C$
and for any
$P \in J_{C} (\overline{K})$,
there should exist an $\epsilon > 0$ such that
\[
\{ x \in C ( \overline{K} )
\mid
|| j_{D} (x) - P||_{NT} \leq \epsilon
\}
\]
is a finite set.
\end{Conjecture}

A stronger version of this conjecture
is also well known as the
effective geometric Bogomolov conjecture for curves:

\begin{Conjecture} \label{EGBCforcurves}
Assume that $C$ is non-isotrivial.
Then
there should exist an $\epsilon > 0$ such that
$
\{ x \in C ( \overline{K} )
\mid
|| j_{D} (x) - P||_{NT} \leq \epsilon
\}
$
is a finite set
for any divisor
$D$ of degree $1$ on $C$
and any 
$P \in J_{C} (\overline{K})$.
Moreover, if $C$ has a stable model over $\mathfrak{B}$,
then we can describe such an $\epsilon$ effectively
in terms of geometric information of the stable model.
\end{Conjecture}

There are some results on Conjecture~\ref{EGBCforcurves}
under the setting where $\mathfrak{B}$ is a curve,
i.e., $K$ is a function field of one variable.
In
$\ch K = 0$,
after partial results by Zhang \cite{zhang1}, 
Moriwaki \cite{moriwaki0,moriwaki3,moriwakiRBI}, 
the author \cite{yamaki1,yamaki4},
and by Faber \cite{faber},
Cinkir proved Conjecture~\ref{EGBCforcurves} in \cite{cinkir}.
In positive characteristic case,
Conjecture~\ref{GBCforcurves} as well as Conjecture~\ref{EGBCforcurves}
are unsolved in full generality,
and there are only some partial results:
Conjecture~\ref{EGBCforcurves} is solved in the case where
\begin{itemize}
\item
the stable model of $C$ has only irreducible fibers (\cite{moriwakiRBI}),
\item
$C$ is a curve of genus $2$ (\cite{moriwaki0}),
\item
$C$ is a hyperelliptic curve (\cite{yamaki4}), or
\item
$C$ is non-hyperelliptic and $g=3$ (\cite{yamaki1}).
\end{itemize}

The
geometric Bogomolov conjecture
for abelian varieties
implies Conjecture~\ref{GBCforcurves}.\footnote{The argument
from here to the end of the proof of
Theorem~\ref{GBCforcurvesC} works well
even if $\mathfrak{B}$ is a
higher dimensional variety.}
Indeed,
let $\hat{h}$ be the N\'eron--Tate height such that
$\hat{h} (x) = || x||_{NT}^{2} $.
Then we have $\hat{h} (j_{D + P} (x)) = || j_{D} (x) - P||_{NT}$
for all $P \in J_{C} ( \overline{K} )$.
Thus 
Conjecture~\ref{GBCforcurves} is equivalent to saying that,
if $C$ is non-isotrivial, then
$j_{D} (C)$ does not have dense small points for any divisor $D$ on $C$
of degree $1$.
To show its contraposition,
we
assume that $j_{D} (C)$ has dense small points for some $D$.
Since Conjecture~\ref{GBCforAV} is assumed to be true,
it follows that $j_{D} ( C )$ is a special subvariety of $J_{C}$.
Since $j_{D}$ is an embedding and since
the trace homomorphism is 
a purely inseparable
finite morphism
(cf. \cite[VIII~\S~3~Corollary~2]{lang2}
or \cite[Lemma~1.4]{yamaki5}),
there are a smooth projective curve
$C_{0}$ over $k$ and a purely inseparable 
finite morphism
$\phi: (C_{0})_{\overline{K}} \to C$.
Let $C_{0} \to C_{0}^{(q)}$ be the $q$-th relative Frobenius morphism,
where $q$ is the degree of $\phi$.
Then its base-change 
to $\overline{K}$ is also the $q$-th relative Frobenius morphism.
It follows from \cite[Corollary~II~2.12]{silverman}
that
$ (C_{0}^{(q)})_{\overline{K}} \cong C$.
This concludes that $C$ is isotrivial.

The following result is obtained
as a consequence of
Theorem~\ref{maintheorem1}
without assumption on the characteristic.
Recall that
a curve $C$ over $\overline{K}$
is of \emph{compact type at $v$} if
the special fiber of the stable model of $C_{v}$
is a tree of smooth irreducible components.
Further, $C$
is said to be of \emph{non-compact type at $v$} if it is not of compact type
at $v$.
It is well known 
that $C$ is of compact type at $v$ if and only if
$C_{v}$ has non-degenerate
Jacobian variety
(cf. \cite[Chapter~9]{BLR}).

\begin{Theorem} \label{GBCforcurvesC}
Suppose that there exists a place 
at which $C$ is of non-compact type.
Then,
for any divisor $D$ of $C$ of degree $1$ and for any $P \in
J_C ( \overline{K} )$,
there exists an
$\epsilon > 0$ such that
$\{ x \in C ( \overline{K} )
\mid
|| j_{D} (x) - P||_{NT} \leq \epsilon
\}
$
is finite.
\end{Theorem}

\begin{Pf}
Let $v$ be a place at which
$C$ is of non-compact type.
Then
the Jacobian $J_{C}$ is degenerate at $v$.
It follows from the Poincar\'e
complete reducibility theorem and Proposition~\ref{cor:new2} that
there is a non-trivial simple abelian variety $A'$
over $\overline{K}$ degenerate at $v$
with a 
surjective homomorphism
$\phi : J_{C} \to A'$.

Let $D$ be any divisor on $C$ of degree $1$.
Let $x_{0} \in j_{D} (C ( \overline{K}))$ be a point.
Note that $J_{C}$ itself is the smallest abelian subvariety
of  $J_{C}$ containing
$j_{D} (C) - x_{0}$.
Since $\phi$ is surjective and since $\dim A' > 0$,
it follows that the image of $j_{D} (C) - x_{0}$ by $\phi$
cannot be a point.
By Lemma~\ref{simple-valuable},
we see that $j_{D} (C)$ is tropically non-trivial.
Since $j_{D} (C)$ has at most finite stabilizer in $J_{C}$,
it follows form Theorem~\ref{maintheorem1} that
$j_{D} (C)$ does not have dense small points.
Thus
the assertion
holds for $C$, as 
is noted above.
\end{Pf}

We end with
a couple of remarks.
In characteristic zero,
Conjecture~\ref{GBCforcurves}
is deduced 
from the combination of \cite[Theorem E]{moriwakiRBI}
and Theorem~\ref{GBCforcurvesC}.
This suggests that
we can avoid hard analysis on metric graphs
carried out
in \cite{cinkir} as far as we consider the non-effective version only.
If the inequality
of \cite[Theorem D]{moriwakiRBI} 
also holds
in positive characteristic,
then the same proof as for \cite[Theorem E]{moriwakiRBI} works,
and hence
we obtain Conjecture~\ref{GBCforcurves}.
Thus our argument makes a contribution to Conjecture~\ref{GBCforcurves}.
We should also note however that 
our approach does not say anything on the effective version
Conjecture~\ref{EGBCforcurves}.


\renewcommand{\thesection}{Appendix} 
\renewcommand{\theTheorem}{A.\arabic{Theorem}}
\renewcommand{\theClaim}{A.\arabic{Theorem}.\arabic{Claim}}
\renewcommand{\theequation}{A.\arabic{Theorem}.\arabic{Claim}}
\renewcommand{\theProposition}
{A.\arabic{Theorem}.\arabic{Proposition}}
\renewcommand{\theLemma}{A.\arabic{Theorem}.\arabic{Lemma}}
\setcounter{section}{0}
\renewcommand{\thesubsection}{A.\arabic{subsection}}



\small{

}


\begin{thebibliography}
{99}
\bibitem{berkovich1}
V.G.\! Berkovich, 
Spectral theory and analytic geometry over nonarchimedean
fields, Mathematical Surveys and Monographs 33, Providence,
RI: AMS (1990).
\bibitem{berkovich2}
V.G.\! Berkovich, \'Etale cohomology for non-archimedean 
analytic spaces,
Publ.\! Math.\! IHES 78, 5--161 (1993).
\bibitem{berkovich3}
V.G.\! Berkovich, 
Vanishing cycles for formal schemes, Invent.\! Math.\! 115,
No.3, 539--571 (1994).
\bibitem{berkovich4}
V.G.\! Berkovich, Smooth p-adic 
analytic spaces are locally contractible,
Invent.\! Math.\! 137, No.1, 1--84 (1999).
\bibitem{bosch}
S.\! Bosch,
Rigid analytische Gruppen mit guter Reduktion,
Math.\! Ann.\! 233 (1976), 193--205.
\bibitem{BL-neron}
S.\! Bosch and W.\! L\"utkebohmert,
N\'eron models from the rigid analytic viewpoint,
J.\! Reine Angew.\! Math.\! 364 (1986), 69--84. 
\bibitem{BL}
S.\! Bosch and W.\! L\"utkebohmert,
Degenerating abelian varieties,
Topology 30 (1991), 653--698.
\bibitem{BLR}
S.\! Bosch, W.\! L\"utkebohmert and
M.\! Raynaud,
N\'eron models,
Ergeb.\! Math.\! Grenzgeb.\! (3) 21,
Springer-Verlag, Berlin, (1990).
\bibitem{chambert-loir}
A.\! Chambert-Loir, 
Mesure et \'equidistribution sur les espaces de Berkovich,
J.\! Reine Angew.\! Math. 595 (2006), 215--235.
\bibitem{cinkir}
Z.\! Cinkir,
Zhang's conjecture and the effective Bogomolov conjecture 
over function fields,
 Invent.\! Math.\! 183, No.3 (2011), 517--562.
\bibitem{dejong}
A.J.\! de Jong,
Smoothness, semi-stability and alterations, 
Inst.\! Hautes \'Etudes Sci.\! Publ.\! Math.\! 83 (1996), 51--93.
\bibitem{faber}
X.W.C.\! Faber,
The geometric Bogomolov conjecture for curves of small genus,
Experiment.\! Math.\! 18 (2009), no. 3, 347--367. 
\bibitem{gubler0}
W.\! Gubler,
Local and canonical height of subvarieties,
Ann.\! Scuola Norm.\! Sup.\! Pisa Cl.\! Sci.\! (5) Vol.\! II (2003),
711--760.
\bibitem{gubler1}
W.\! Gubler,
Tropical varieties for non-archimedean analytic spaces,
Invent.\! Math.\! 169 (2007), 321--376.
\bibitem{gubler2}
W.\! Gubler,
The Bogomolov conjecture for totally degenerate abelian varieties,
Invent.\! Math.\! 169 (2007), 377--400.
\bibitem{gubler3}
W.\! Gubler,
Equidistribution over function fields,
manuscripta math.\! 127 (2008), 485--510.
\bibitem{gubler4}
W.\! Gubler,
Non-archimedean canonical measures on abelian varieties,
Compos.\! Math.\! 146 (2010), 683--730.
\bibitem{lang1}
S.\! Lang,
Abelian varieties,
Springer-Verlag (1983).
\bibitem{lang2}
S.\! Lang,
Fundamentals of Diophantine Geometry,
Springer-Verlag (1983).
\bibitem{milne}
J.\! Milne,
Abelian varieties,
\emph{in} Arithmetic geometry (Storrs, Conn., 1984), 103--150, 
Springer, New York (1986).
\bibitem{moriwaki0}
A.\! Moriwaki,
Bogomolov conjecture for curves of genus $2$
over function fields, J.\! Math.\! Kyoto Univ.\! 36 (1996), 687--695. 
\bibitem{moriwaki3}
A. Moriwaki,
Bogomolov conjecture over 
function fields for stable curves with only irreducible fibers,
Compositio Math.\! 105 (1997), 125--140. 
\bibitem{moriwakiRBI}
A. Moriwaki,
Relative Bogomolov's inequality and the cone of positive divisors 
on the moduli space of stable curves,
J.\! of AMS 11 (1998), 569--600.
\bibitem{moriwaki5}
A.\! Moriwaki,
Arithmetic height functions
over finitely generated fields,
Invent.\! Math.\! 140 (2000), 101--142.
\bibitem{mumford}
D. Mumford,
Abelian varieties
(with appendices by C. P. Ramanujam and Yuri Manin),
corrected reprint of the second (1974) edition,
TIFR Studies in Math. 5,
Hindustan Book Agency, New Delhi, 2008,
ISBN: 978-81-85931-86-9; 81-85931-86-0.
\bibitem{scanlon}
T.\! Scanlon,
A positive characteristic Manin--Mumford theorem,
Compos.\! Math.\! 141 (2005), 1351--1364.
\bibitem{silverman}
J.H.\! Silverman,
The arithmetic of elliptic curves,
GTM 106,
Springer-Verlag, New York (1986).
\bibitem{SUZ}
L.\! Szpiro, E.\! Ullmo and S.\! Zhang, 
Equir\'epartition des petits points, Invent. Math. 127 (1997),
337--347.
\bibitem{ullmo}
E.\! Ullmo,
Positivit\'e et discr\'etion des points alg\'ebriques des courbes,
Ann.\! of Math.\! 147 (1998), 167--179.
\bibitem{yamaki1}
K.\! Yamaki,
Geometric Bogomolov's conjecture
for curves of genus $3$
over function fields,
J.\! Math.\! Kyoto Univ.\! 42 (2002),
57--81.
\bibitem{yamaki4}
K.\! Yamaki,
Effective calculation of the geometric height
and the Bogomolov conjecture
for hyperelliptic curves over
function fields,
J.\! Math.\! Kyoto.\! Univ.\! 48 (2008),
401--443.
\bibitem{yamaki5}
K.\! Yamaki,
Geometric Bogomolov conjecture for abelian varieties 
and some results for those with some degeneration 
(with an appendix by Walter Gubler: 
The minimal dimension of a canonical measure),
Manuscr.\! Math.\! 
142 (2013), no.\! 3-4, 273--306.
\bibitem{zhang1}
S.\! Zhang,
Admissible pairing on a curve, 
Invent.\! Math.\! 112 (1993), 171--193.
\bibitem{zhang2}
S.\! Zhang,
Equidistribution of small points on abelian varieties,
 Ann.\! of Math.\! (2) 147 (1998), no.\! 1, 159--165.
\end{thebibliography}
\end{document}